\documentclass[12pt]{article}
\usepackage{amssymb,amsmath}
\usepackage{graphicx}
\usepackage{a4}
\numberwithin{equation}{section}

\begin{document}
\baselineskip16pt
\title{\bf{Refined Estimates on Conjectures of Woods and Minkowski}}
\author{Leetika Kathuria{\footnote{The author acknowledges the support of CSIR, India. The paper forms a part of her
Ph.D. dissertation accepted by Panjab University, Chandigarh.  }}~ and  Madhu Raka
 \\ \small{\em Centre for Advanced Study in Mathematics}\\
\small{\em Panjab University, Chandigarh-160014, INDIA}\\
\date{}\vspace{-1cm}}
\maketitle
 {\abstract{$~~$\par Let $\wedge$ be a lattice in $\mathbb{R}^n$ reduced in the sense of Korkine and Zolotareff having a  basis of the form $(A_1,0,0,\ldots,0),(a_{2,1},A_2,0,\ldots,0)$, $\ldots,(a_{n,1},a_{n,2},\ldots,a_{n,n-1},A_n)$ where $A_1, A_2,\ldots,A_n$ are all positive.  A well known conjecture of Woods in  Geometry of Numbers asserts that if $A_{1}A_{2}\cdots A_{n}=1$ and $A_{i}\leqslant A_{1}$ for each $i$ then  any closed sphere in $\mathbb{R}^n $ of  radius $ \sqrt{n}/2$ contains a point of $\wedge$. Woods' Conjecture is known to be true for $n\leq 9$. In this paper we give estimates on  the Conjecture of Woods for $10\leq n\leq33$, improving the earlier best known results of Hans-Gill et al. These lead to an improvement, for these values of  $n$,  to the estimates on the long standing classical conjecture of Minkowski  on the product of $n$ non-homogeneous linear forms.
 }}

\section{Introduction}
Let $ L_{i}=a_{i1}x_{1}+\cdots+ a_{in}x_{n},~1\leq i\leq n$ be $n$  real linear forms in $n$ variables $x_{1},\ldots, x_{n}$ and having determinant $\Delta=\mbox{det}(a_{ij})\neq0$. The following conjecture is attributed to H. Minkowski:\\
\textbf{Conjecture I:}~\emph{For any given real numbers $c_1,\ldots, c_n$, there exists integers $x_1,\ldots,x_n$ such that
\begin{equation}
\mid(L_{1}+c_{1})\cdots(L_{n}+c_{n})\mid \leqslant \frac{1}{2^{n}}\mid \Delta \mid.
\end{equation}
Equality is necessary if and only if after a suitable unimodular transformation the linear forms $L_{i}$ have the form $2c_{i}x_{i}$ for $1\leq i \leq n$.\vspace{2mm}\\}
 This result is known to be true for $n\leq9$. For a detailed history and the related results, see Bambah et al\cite{BDH},
Gruber\cite{Gr4}, Hans-Gill et al \cite{HRS1} and Kathuria and Raka \cite{KR1}. \vspace{2mm}\\
 Minkowski's Conjecture is equivalent to saying that $$M_{n}\leqslant \frac{1}{2^{n}}\mid \Delta \mid,$$ where $M_{n}$ = $M_{n}(\Delta)$ is given by
$$\begin{array}{cccccc}
M_{n}& = &{\mbox {Sup}}&\hspace{-0.9cm}{\mbox {Sup}}&\hspace{-0.9cm}{\mbox {Inf}}&\hspace{-0.3cm}{\displaystyle\prod _{i=1}^{n}}\mid L_{i}(u_1,\ldots, u_n)+c_{i}\mid.\vspace{-6mm}\\
      &  &{\scriptstyle{L_{1},\cdots,L_{n}}}&\hspace{-0.2cm}{\scriptstyle(c_{1},\cdots,c_{n})\in \mathbb{R}^{n}}&\hspace{-0.2cm}{\scriptstyle(u_{1},\cdots,u_{n})\in \mathbb{Z}^{n}}
       \end{array}$$
Chebotarev~\cite{Ch} proved the weaker inequality \begin{equation} M_n \leq {\frac{1}{2^{n/2}}}|\Delta|.\end{equation}
Since then several authors have tried to improve upon this estimate. The bounds have been obtained in the form
\begin{equation} M_n \leq {\frac{1}{\nu_n 2^{n/2}}}|\Delta|,\end{equation} where $\nu_n >1$. Clearly $\nu_n \leq 2^{n/2}$ by considering the linear forms $ L_i =~x_i$ and $ c_i~= \frac{1}{2}$ for $ 1\leq i \leq n$. During 1949-1986, many authors such as Davenport, Woods, Bombieri, Gruber, Skubenko, Andrijasjan, Il'in and Malyshev obtained $\nu_n$ for large $n$. For  details, see Gruber and Lekkerkerker~\cite{GL}, Hans-Gill et al \cite{HRS2}. In 1960, Mordell~\cite{Mo4} obtained $ \nu_n= 4 - 2(2 - 3\sqrt{2}/4)^n - 2^{-n/2}$ for all $n$.  Il'in~\cite{IL1, IL2}(1986,1991) improved Mordell's  estimates for $6 \leq n \leq 31$. Hans-Gill et al~\cite{HRS2,HRS4}(2010,2011) got improvements on the results of Il'in~\cite{IL2} for $9 \leq n \leq 31$. Since recently $\nu_9=2^{9/2}$ has been established  by the authors \cite{KR1}, we  study $\nu_n$ for $10\leq n \leq 33$ and obtain their refined values in this paper. \vspace{1mm}

 For sake of comparison, we give results by Mordell~\cite{Mo4},  Il'in~\cite{IL2}, Hans-Gill et al~\cite{HRS4}  and our improved $\nu_n$  in Table I.\vspace{2mm}

 We shall follow the Remak-Davenport approach. For the sake of convenience of the reader we give some basic results of this approach.\\
\noindent Minkowski's Conjecture can be restated in the terminology of lattices as :\\ Any lattice $\wedge$ of determinant $d(\wedge)$ in $\mathbb{R}^n$  is a covering lattice for the set
                   $$  S : | x_{1}x_{2}\dots x_{n}|\leq \frac{d (\wedge )}{2^{n}}.$$
 \newpage
\begin{center}{\small {\bf Table I}
\vspace{3mm}\\ \begin{tabular}{|l|l|l|l|l|}
  \hline
& Estimates by &  Estimates by
& Estimates by  & Our improved \\&Mordell&Il'in &Hans-Gill et al&Estimates\\
\hline
 ~~$n$ & ~~~$\nu_n$& ~~~$\nu_n$  &~~~~ $\nu_n$ &~~~~ $\nu_n$
\\
\hline
~~10 &2.8990614& 3.4798928&24.3627506&27.6034811\\
  ~~11& 2.9731018&3.5229055  &29.2801145&33.4727227\\
  ~~12& 3.0405253&3.5502417   &32.2801213&39.5919904\\
  ~~13&3.1023558& 3.5785628   &34.8475153&45.4004068\\
  ~~14& 3.1593729&3.6020935  &37.8038391&51.2623882\\
  ~~15& 3.2121798&3.6111553  &40.9051980&57.0037507\\
  ~~16& 3.2612520&3.6190753   &44.3414913&57.4701963\\
  ~~17& 3.3069717&3.6392444   &47.2339309&57.6759791\\
  ~~18& 3.3496524&3.6617581   &46.7645724&57.3887589\\
  ~~19&3.3895562& 3.6673429   &47.2575897&60.0933912\\
  ~~20&3.4269065& 3.6723611  &46.8640155&58.4859214\\
  ~~21& 3.4618973&3.6769169   &46.0522028&56.4257125\\
  ~~22&3.4946990&3.684080   &43.6612034&53.9414220\\
  ~~23& 3.5254641& 3.6863331 &37.8802374&50.9884152\\
  ~~24&3.5543297& 3.6897821  &32.5852958&47.7463213\\
~~25&3.5814208& 3.6929517 &27.8149432&42.3908768\\
  ~~26& 3.6068520& 3.6958893 &23.0801951&38.8656991\\
  ~~27&  3.6307288& 3.7001150  &17.3895105&31.9331584\\
  ~~28&  3.6531489& 3.7026271 &12.9938763&26.1066323\\
  ~~29&  3.6742031& 3.7049722 &9.5796191&19.9625412\\
  ~~30 & 3.6939760&3.7086731 &6.7664335&16.0688443\\
  ~~31 & 3.7125466& 3.7255824 &4.7459720&11.2387160\\
  ~~32 & 3.7299885&  &&8.3258788\\
~~33 & 3.746371&  &&5.4114880\\

  \hline
\end{tabular}}\end{center}\vspace{2mm}

The weaker result (1.3) is equivalent to saying that any lattice $\wedge$ of determinant $d(\wedge)$ in $\mathbb{R}^n$  is a covering lattice for the set
                   $$  S : | x_{1}x_{2}\dots x_{n}|\leq \frac{d (\wedge )}{\nu_n 2^{n/2}}.$$
Define the homogeneous minimum of $\wedge$ as
$$m_H(\wedge)~=~{\rm Inf}\{|x_1x_2\dots x_n|:~X=(x_1,x_2,\cdots,x_n)\in \wedge, X\neq O\}.$$

\noindent In 1956, Birch and Swinnerton-Dyer\cite{BS} proved \vspace{2mm}\\
{\bf Proposition 1.} Suppose that Minkowski Conjecture has been proved for dimensions $1,2,\cdots, n-1.$ Then it holds for all lattices $\wedge$ in $\mathbb{R}^n$ for which $M_H(\wedge) =0$.\vspace{2mm}

\noindent
C.T. McMullen\cite{Mc} proved \vspace{2mm} \\{\bf Proposition 2. }If  $\wedge$ is a lattice in
$\mathbb{R}^n$ for $n\geq 3$ with $M_H(\wedge) \neq 0$ then there exists  an ellipsoid having $n$
linearly independent points of $\wedge$ on its boundary and
 no point of $\wedge$ other than $O$ in its interior.\vspace{2mm}

It is well known that using these results, Minkowski's Conjecture would follow from \vspace{2mm}

{\noindent \bf Conjecture II.} If $\wedge$ is a lattice in
$\mathbb{R}^n$ of
determinant 1 and there is a sphere $|X|<R$ which contains no
point of $\wedge$ other than $O$ in its interior and has $n$ linearly independent
points of $\wedge$ on its boundary then $\wedge$ is a covering
lattice for the closed sphere of radius $ \sqrt{n/4}.$ Equivalently, every closed sphere of radius $\sqrt{n/4}$ lying in $\mathbb{R}^n$ contains a point of $\wedge$.\vspace{2mm}

%

 Woods ~\cite{W3, W4} formulated a conjecture from which Conjecture-II follows immediately. To state Woods' conjecture, we need to introduce some terminology :
\par Let $\mathbb{L}$ be a lattice in $\mathbb{R}^{n}$. By the reduction theory of quadratic forms introduced by Korkine and Zolotareff ~\cite{KZ1}, a cartesian co-ordinate system may be chosen in $\mathbb{R}^{n}$ in such a way that $\mathbb{L}$ has a basis of the form $$(A_1,0,0,\ldots,0),(a_{2,1},A_2,0,\ldots,0),\ldots,(a_{n,1},a_{n,2},\ldots,a_{n,n-1},A_n),$$ where $A_1, A_2,\ldots,A_n$ are all positive and further for each $i=1,2,\ldots,n$ any two points of the lattice in $\mathbb{R}^{n-i+1}$ with basis
$$(A_i,0,0,\ldots,0),(a_{i+1,i},A_{i+1},0,\ldots,0),\ldots,
(a_{n,i},a_{n,i+1},\ldots,a_{n,n-1},A_n)$$
 are at a distance atleast $A_i$ apart. Such a basis of $\mathbb{L}$ is called a reduced basis.\vspace{2mm}\\
\textbf{Conjecture III (Woods)}:\emph{ If $A_{1}A_{2}\cdots A_{n}=1$ and $A_{i}\leqslant A_{1}$ for each $i$ then any closed sphere in $\mathbb{R}^{n}$ of radius  $\sqrt{n}/2$ contains a point of $\mathbb{L}$.}\vspace{2mm}

Woods~\cite{W2, W3, W4} proved this conjecture for $4\leqslant n \leqslant 6$. Hans-Gill et al~\cite{HRSS} gave a unified proof of Woods' Conjecture for $n \leqslant 6$. Hans-Gill et al ~\cite{HRS2,HRS4} proved Woods' Conjecture  for $n=7$ and $n=8$ and thus completed the proof of Minkowski's Conjecture for $n=7$ and 8. Kathuria and Raka~\cite{KR1} proved Woods Conjecture and hence Minkowski's Conjecture for $n=9$.
 With the assumptions as in Conjecture III, a weaker result would be that\vspace{2mm}\\ \textbf{If $\omega_n\geq n$, any closed sphere in $\mathbb{R}^n$ of radius $\sqrt{\omega_n}/2$ contains a point of $\mathbb{L}$.}\vspace{2mm}

  Hans-Gill et al~\cite{HRS2, HRS4} obtained the estimates $\omega_n$ on  Woods' Conjecture for $n\geq 9$. As $\omega_9=9$ has been established by the authors \cite{KR1} recently, in this paper we obtain improved estimates $\omega_n$ on Woods' Conjecture for $10\leq n\leq33$.
 Together with the following result  of Hans-Gill et al.~\cite{HRS2}, we get improvements of $\omega_n$ for  $n\geq 34$ also.\vspace{2mm}

\noindent\textbf{Proposition 3.} \emph{Let $\mathbb{L}$ be a lattice in $\mathbb{R}^{n}$ with $A_{1}A_{2}\cdots A_{n}=1$ and $A_{i}\leqslant A_{1}$ for each $i$. Let $0 < l_n \leq A_n^2 \leq m_n$, where
$l_n$ and $m_n$ are real numbers. Then $\mathbb{L}$ is a covering lattice for the sphere
$|X|\leq{\sqrt{\omega_n}}/2$, where
$\omega_n$ is defined inductively by
$$\omega_n=\max\{\omega_{n-1}l_n^{-1/l_{n-1}}+l_n, \omega_{n-1}m_n^{-1/m_{n-1}}+m_n\}.$$}
Here we prove\\
\noindent{\bf Theorem 1.}
\emph{Let $10\leq n\leq33$. If $d(\mathbb{L})= A_{1}\dots A_{n} = 1 $ and $A_{i }\leq A_{1}$ for $i =2,\cdots, n$, then any closed sphere in
$\mathbb{R}^n$ of radius $\sqrt{\omega_n}/2$ contains a point of
$\mathbb{L}$, where $\omega_{n}$ are as listed in Table II.}\vspace{2mm}\\For the sake of comparison we give results by Hans-Gill et. al \cite{HRS4} and our improved $\omega_n$ in Table II.

To deduce the results on the estimates of Minkowski's Conjecture we also need
the following generalization of Proposition 1 (see Theorem 3 of \cite{HRS2}; for a proof see \cite{KR4}):\vspace{2mm}\\
\textbf{Proposition 4.}
 \emph{Suppose that we know
$$M_{j}\leqslant \frac{1}{\nu_{j}2^{j/2}\mid\Delta\mid}
~\mbox{for}~ 1 \leqslant j \leqslant n-1.$$
Let $\nu_n < \min {\nu_{k_1}\nu_{k_2}\cdots\nu_{k_s}}$, where the minimum is taken over all $(k_{1}, k_{2},\cdots, k_{s})$
such that $n = k_{1} + k_{2} + \ldots + k_{s}$, $k_{i}$ positive integers for all $i$ and $s\geqslant 2$. Then
for all lattices $\Lambda$ in $\mathbb{R}^{n}$ with homogeneous minimum $M_{H}(\Lambda) = 0$, the estimate $\nu_n$ holds for
Minkowski's Conjecture}.\vspace{2mm}\\
Since by arithmetic-geometric inequality the sphere $\{X \in \mathbb{R}^n :~  |X| \leq \frac{\sqrt{\omega_n}}{2} \}$  is a subset of $\{ X : |x_1x_2...x_n|\leq \frac{1}{2^{n/2}}(\frac{\omega_n}{2n})^{n/2}\}$, Propositions 2 and 4 immediately imply \vspace{2mm} \\
\noindent{\bf Theorem 2:} \emph{The values of $\nu_n$  for
 the estimates of Minkowski's Conjecture can be taken as $(\frac{2n}{\omega_n})^{n/2}$.\\}
For $10\leq  n \leq 33$, these values are listed in Table I.
In Section 2 we state some preliminary results and in Sections 3-9 we  prove Theorem 1 for $10\leq n\leq 33$.
\section{Preliminary Results and Plan of the Proof}
Let $\mathbb{L}$ be a lattice in $\mathbb{R}^n$ reduced in the sense of Korkine and Zolotareff. Let $\Delta (S_n)$ denotes the critical determinant of
the unit sphere $S_n$ with center $O$
in $\mathbb{R}^n$ i.e.\vspace{2mm}\\ $ ~~~\Delta (S_n)= {\mbox {Inf}}\{d(\wedge) :  ~\wedge {\mbox {~has no   point other  than }}O{\mbox {~in the interior of ~}} S_n\}.$ \vspace{1mm}\\
Let $\gamma_n$ be the Hermite's constant i.e. $\gamma_n$ is the smallest real number such that for any positive definite quadratic form $Q$ in $n$ variables of determinant $D$, there exist integers $u_1,u_2,\cdots, u_n$ not all zero satisfying $$ Q(u_1,u_2,\cdots, u_n) \leq \gamma_n D^{1/n}.$$
 It is well known that $\Delta^{2}(S_{n})=\gamma_n^{-n}$. We write $A_i^{2}=B_i$. \vspace{1mm} \\We state below some preliminary lemmas. Lemmas 1 and 2 are due to Woods~\cite{W2}, Lemma 3 is due to Korkine and Zolotareff~\cite{KZ1},
 and Lemma 4 is due to  Pendavingh and Van Zwam~\cite{Pen}.
 In Lemma 5, the cases $n=2$ and $3$ are classical results of Lagrange and
 Gauss; $n=4$ and $5$ are due to Korkine and Zolotareff~\cite{KZ1} while $n=6, 7$ and $8$ are due to Blichfeldt ~\cite{Bl}.\vspace{2mm}\\
{\noindent {\bf Lemma 1.}  If $2 \Delta (S_{n+1})A_1^n\geq d(\mathbb{L})$ then any closed sphere of radius
$$R=A_1(1-\{A_1^n\Delta (S_{n+1}) /d(\mathbb{L})\}^2)^{1/2}$$ in $\mathbb{R}^n$ contains a point of $\mathbb{L}.$\vspace{1mm}\\
{\noindent {\bf Lemma 2. }
 For a fixed integer $i$ with $1\leq i\leq
n-1,$ denote by $\mathbb{L}_1$ the lattice in $\mathbb{R}^i$ with
 reduced basis\vspace{-2mm}
$$(A_1,0,\cdots,0),(a_{2,1},A_2,0,\cdots,0),\cdots,(a_{i,1},a_{i,2},\cdots,a_{i,i-1},A_i)$$
and denote by $\mathbb{L}_2$ the lattice in $\mathbb{R}^{n-i}$ with  reduced basis \vspace{-2mm}
$$(A_{i+1},0,\cdots,0),(a_{i+2,i+1},A_{i+2},0,\cdots,0),\cdots,(a_{n,i+1},a_{n,i+2},\cdots,a_{n,n-1},A_n).$$
If any closed sphere in  $\mathbb{R}^i$ of radius $r_1$ contains a point
of $\mathbb{L}_1$ and if any  closed sphere in $\mathbb{R}^{n-i}$ of
radius $r_2$ contains a point of $\mathbb{L}_2$ then any closed sphere in
$\mathbb{R}^n$ of radius $(r_1^2+r_2^2)^{1/2}$ contains a point of
$\mathbb{L}.$\vspace{1mm}\\
{\noindent {\bf Lemma 3.} For all relevant $i$,\begin{equation}\label{2.0a} B_{i+1}\geq\frac{3}{4}B_i  ~\mbox{and}~
B_{i+2}\geq \frac{2}{3}B_i .\end{equation}
{\noindent {\bf Lemma 4. } For all relevant $i$, \begin{equation}\label{2.0b}B_{i+4}\geq (0.46873) B_i .\end{equation}
\textbf{Throughout the paper we shall denote 0.46873 by $\varepsilon$.}\vspace{2mm}\\
{\noindent {\bf Lemma 5.}$~~~\Delta (S_n) = ~\sqrt{3}/2, ~1/\sqrt{2}, ~1/2,
~1/2\sqrt{2}, ~\sqrt{3}/8,~1/8$ and $1/16$  for
 $n=~2,~3,~4,~5,~6,~7~$ and $~8 {\mbox ~{respectively}}.$\vspace{2mm}\\
\noindent {\bf Lemma 6.} For any integer $s, ~ 1\leq s\leq n-1$
\begin{equation}\label{2.1} B_1B_2\cdots B_{s-1}B_{s}^{n-s+1}\leq\gamma_{n-s+1}^{n-s+1} ~\mbox{and}\end{equation}
\begin{equation}\label{2.2} B_1B_2 \dots B_s\leq (\gamma_n^{\frac{1}{n-1}}\gamma_{n-1}^{\frac{1}{n-2}}\dots \gamma_{n-s+1}^{\frac{1}{n-s}})^{n-s}.\end{equation}
This is Lemma 4 of Hans-Gill et al~\cite{HRS2}.\vspace{1mm}\\
 \noindent {\bf Lemma 7.} \begin{equation}\label{2.3}\{(8.5337)^{\frac{1}{5}}\gamma_n^{\frac{1}{n-1}}\gamma_{n-1}^{\frac{1}{n-2}}\dots \gamma_6^{\frac{1}{5}}\}^{-1}
 \leq B_n \leq \gamma_{n-1}^{\frac{n-1}{n}}.\end{equation}
 This is Lemma 6 of Hans-Gill et al~\cite{HRS4}.\vspace{2mm} \\
\noindent{\bf Remark 1.} Let \vspace{2mm}\\
$\begin{array}{ll}\delta_n~  &= {\mbox {~ the ~best~ centre~ density~ of ~packings~ of ~unit~ spheres ~in ~}} \mathbb{R}^{n},\\ \delta_n^{\star} &= {\mbox {~ the~ best~ centre~ density ~of ~lattice~ packings~ of ~unit~ spheres ~in~}} \mathbb{R}^{n}.\end{array}$\vspace{3mm}\\ Then it is known that (see Conway and Sloane~\cite{CS}, page 20)
\begin{equation}\label{2.4} \gamma_n= 4 (\delta_n^{\star})^{\frac{2}{n}}\leq 4 (\delta_n)^{\frac{2}{n}}.\end{equation} $\delta_n^{\star}$ and hence $\gamma_n$ is known for $ n \leq 8$. Also $\gamma_{24}=4$ has been proved by  Cohn and Kumar~\cite{CK}. Using the bounds on $\delta_n$ given by Cohn and Elkies~\cite{CE} and  inequality \eqref{2.4} we find  bounds on  $\gamma_n$ for $10\leq n \leq 33$ which we have listed in Table II. Also $\gamma_9\leq 2.1326324$.\vspace{2mm}
\begin{center}{\small {\bf Table II}
\vspace{3mm}\\ \begin{tabular}{|c|c|c|c|c|}
 \hline
 $n$&
 $\gamma_n  \leq $& Estimates by Hans-Gill et al  & Our improved Estimates\\& & $\omega_n$  & $\omega_n$\\\hline
  10 & 2.2636302 & 10.5605061&10.3\\
  \hline
  11 &  2.3933470 &11.9061976&11.62\\
  \hline
  12&  2.5217871 & 13.4499927& 13\\\hline
  13&  2.6492947 &15.0562267&14.455765 \\\hline
14&  2.7758041 & 16.6646332&15.955156 \\\hline
15&  2.9014777 & 18.2901579&17.498499\\\hline
  16& 3.0263937 & 19.9204292&19.285\\\hline
  17&  3.1506793&  21.6026907&21.101\\\hline
18&  3.2743307 & 23.4831402&22.955\\\hline
  19&  3.3974439 & 25.3234826&24.691\\\hline
  20&  3.5200620 & 27.2255111&26.629\\\hline
   21&  3.6422432 & 29.1638254&28.605\\\hline
  22&  3.7640371 & 31.2142617&30.62\\\hline
23&  3.8854763&33.5354821&32.68\\\hline
  24&        4.0065998     &35.9050965&34.78\\\hline
25&4.1274438&38.3201985&37.05\\\hline
%
  26& 4.2480446&40.8449876&39.24\\\hline
 27& 4.3684312&43.7039431&41.78\\ \hline
28& 4.488631&46.6267624&44.36\\\hline 29& 4.6086676&49.6305176&47.18\\\hline
30&4.7285667&52.8194566&49.86\\\hline
31&4.8483483&56.0735184&53.04\\\hline
32&4.9680344&&56.06\\\hline
33&5.0876409&&59.58\\\hline

\end{tabular}}
\end{center}

We assume that Theorem 1  is false  and derive a contradiction.  Let $\mathbb{L}$ be a lattice
satisfying the hypothesis of the conjecture. Suppose that there exists a  closed sphere of radius $\sqrt{\omega_n}/2$
in $\mathbb{R}^n$ that contains no point of $\mathbb{L}$ in $\mathbb{R}^n$.  Since $B_i=A_i^2$
and $d(\mathbb{L})=1,$ we have $B_1B_2\dots B_n = 1.$ \vspace{2mm}\\
$~~~~$We give some examples of inequalities that arise. Let $\mathbb{L}_1$ be a lattice in $\mathbb{R}^4$ with basis
$(A_1,0,0,0)$, $(a_{2,1},A_2,0,0),$ $(a_{3,1},a_{3,2},A_{3},0),$  $(a_{4,1},a_{4,2},a_{4,3},A_{4}),$ and $\mathbb{L}_i$  for $2\leq i \leq n$ be lattices in
$\mathbb{R}^1$ with basis $(A_{i+3})$. Applying Lemma 2  repeatedly and using Lemma 1   we see that  if
$~2\Delta(S_5)A_1^4$ $ \geq A_1A_2A_3A_4$ then any closed sphere of radius
$$\left(A_1^2-\frac{A_1^{10}
\Delta(S_5)^2}{A_1^2A_2^2A_3^2A_4^2}+\frac{1}{4}A_5^2+\cdots+\frac{1}{4}A_n^2\right)^{1/2}$$
\noindent contains a point of $\mathbb{L}.$ By the initial hypothesis this radius exceeds $\sqrt{\omega_n}/2$.  Since $\Delta(S_5)=1/2\sqrt{2}$ and
$B_1B_2...B_n = 1,$ this results in the conditional inequality : ${\mbox {if}}~~ B_1^{4}B_5 B_6 \dots B_n \geq 2$ ~~{\mbox then}~~
\begin{equation}\label{3.1}
 4B_1-\frac{1}{2}B_1^5B_5B_6\dots B_n+B_5+B_6+\cdots+B_n>\omega_n.
\end{equation}
We call this inequality  $(4,1,\cdots,1),$ since it corresponds to the ordered partition $(4,1,\cdots,1)$ of $n$ for the
purpose of applying  Lemma 2.    Similarly the conditional inequality
 $(1,\cdots,1,2,1,\cdots,1)$ corresponding to the ordered partition  $(1,\cdots,1,2,1,\cdots,1)$ is : ${\mbox {if}~~} 2B_i\geq B_{i+1}~~ {\mbox {then}}$
\begin{equation}\label{3.2}
 B_1+\cdots+B_{i-1}+4B_i-\frac{2B_i^2}{B_{i+1}}+B_{i+2}+\cdots+B_n>\omega_n.
\end{equation}
Since $4B_i-\frac{2B_i^2}{B_{i+1}}\leq 2B_{i+1}$, \eqref{3.2} gives
\begin{equation}\label{3.3}
B_1+\cdots+B_{i-1}+2B_{i+1}+B_{i+2}+\cdots+B_n>\omega_n.
\end{equation}
One may remark here that the condition $2B_i\geq B_{i+1}$ is necessary only if we want to use inequality \eqref{3.2}, but it is not necessary if we want to  use the weaker inequality \eqref{3.3}. This is so because if $2B_i< B_{i+1}$, using the partition $(1,1)$ in place of $(2)$ for the relevant part, we get the upper bound $B_i+B_{i+1}$ which is clearly less than $2B_{i+1}$. We shall call inequalities of type \eqref{3.3} as \emph{weak} inequalities and denote it by $(1,\cdots,1,2,1,\cdots,1)_w$. \vspace{2mm}\\
$~~~~$ If $(\lambda_{1},\lambda_{2},\cdots,\lambda_{s})$ is an ordered partition of $n$, then the conditional inequality arising
from it, by using Lemmas 1 and 2, is also denoted by $(\lambda_{1},\lambda_{2},\cdots,\lambda_{s})$. If the conditions
in an inequality $(\lambda_{1},\lambda_{2},\cdots,\lambda_{s})$ are satisfied then we say that
$(\lambda_{1},\lambda_{2},\cdots,\lambda_{s})$ holds.\vspace{2mm}\\
$~~~~$ Sometimes, instead of Lemma 2, we are able to use induction. The
use of this is indicated by putting $(^{*})$ on the corresponding part of the partition.
For example, if for $n=10$, $B_5$ is larger than each of $B_6, B_7, \cdots, B_{10}$, and if $\frac{B_1^3}{B_2B_3B_4}>2$, the inequality $(4,6^*)$  gives \begin{equation}\label{3.4} 4B_1-\frac{1}{2}\frac{B_1^4}{B_2B_3B_4}+6(B_1B_2B_3B_4)^{-1/6}>\omega_{10}.\end{equation}

 \noindent In particular the inequality
$((n-1)^{*},1)$ always holds. This can  be written as  \begin{equation}\label{3.5}\omega_{n-1}( B_{n})^{\frac{-1}{(n-1)}}+
B_n >\omega_n.\end{equation}
 Also we have $B_1\geq 1$, because if $B_1< 1$, then  $B_i\leq B_1<1$ for each $i$ contradicting  $B_1B_
2...B_n=1$.\vspace{2mm} \\
  Using the upper bounds on $\gamma_n$  and the inequality \eqref{2.3}, we obtain numerical lower and upper bounds on $B_n$, which we denote by $l_n$ and $m_n$ respectively.  We use the approach of Hans-Gill et al~\cite{HRS4}, but our method of dealing with is somewhat different. In Sections 3-5 we give proof of Theorem 1 for $n=10, 11 ~\mbox{and}~ 12$ respectively. The proof of these cases is based on the truncation of the interval $[l_n,m_n]$ from both the sides. In Sections 6-8 we give proof of Theorem 1 for $n=13, 14 ~\mbox{and}~ 15$. The proof of these cases is based on the truncation of the interval $[l_n,m_n]$ from  one side only. (Truncation from both the sides makes the proof very complicated and it does not give any significant improvement on $\omega_n$.) For $16\leq n\leq33$ we use
the inequalities in somewhat different way and this is discussed in Section 9.\vspace{2mm} \\
In this paper we need to maximize or minimize frequently functions of several variables. When we say that a given function of several variables in $x,y,\cdots$ is an increasing/decreasing function of $x,y,\cdots$, it means that the concerned property holds when function is considered as a function of one variable at a time, all other variables being fixed.
\section{Proof of Theorem 1 for $n=10$}
\noindent{Here we have $\omega_{10}=10.3$, $B_1 \leq \gamma_{10}< 2.2636302$. Using \eqref{2.3}, we have $l_{10}=0.4007<B_{10}<1.9770808=m_{10}$. \vspace{2mm}\\
The inequality $(9^*,1)$ gives $9(B_{10})^{\frac{-1}{9}}+B_{10}<10.3$. But for $0.4398\leq B_{10}\leq 1.9378$, this inequality is not true. Hence we must have either $B_{10}<0.4398$ or $B_{10}>1.9378$.\\ We will deal with the two cases $0.4007<B_{10}<0.4398$ and $1.9378<B_{10}<1.9770808$ separately:\vspace{2mm}

\subsection{$0.4007<B_{10}<0.4398$}
Using \eqref{2.0a},\eqref{2.0b}   we have:
\begin{equation}{\label{3.6}}\left \{ \begin{array}{lll}
B_9 \leq \frac{4}{3}B_{10}<0.5864,& B_8 \leq \frac{3}{2}B_{10}<0.6597,& B_7 \leq 2B_{10}<0.8796,\vspace{1mm} \\
B_6 \leq \frac{B_{10}}{\varepsilon}<0.9383,&B_5 \leq \frac{4}{3}\frac{B_{10}}{\varepsilon} < 1.2511,& B_4 \leq \frac{3}{2}\frac{B_{10}}{\varepsilon}< 1.4075,\vspace{1mm}\\
B_3\leq \frac{2B_{10}}{\varepsilon}<1.8766, &B_2 \leq \frac{B_{10}}{(\varepsilon)^2}<2.0018&.
 \end{array}\right.{}\end{equation}
\noindent {\bf Claim(i)} $B_2>1.7046$

The  inequality $(2,2,2,2,2)_w$ gives $ 2B_2+2B_4+2B_6+2B_8+2B_{10} > 10.3$. Using (3.1), we find that this inequality is not true for $B_2\leq1.7046$. Hence we must have $B_2> 1.7046$.\vspace{2mm}\\
\noindent {\bf Claim(ii)} $B_2<1.8815$

Suppose $B_2\geq1.8815$, then using (3.1) and that $B_6\geq \varepsilon B _2$ we find that $\frac{B_2^3}{B_3B_4B_5}>2$ and $\frac{B_6^3}{B_7B_8B_9}>2$. So the inequality (1,4,4,1) holds, i.e. $B_1+4B_2-\frac{1}{2}\frac{B_2^4}{B_3B_4B_5}+4B_6-\frac{1}{2}\frac{B_6^4}{B_7B_8B_9}+B_{10}>10.3$. Applying AM-GM inequality we get $B_1+4B_2+4B_6+B_{10}-\sqrt{B_2^5B_6^5B_1B_{10}}>10.3.$ Now since ${\varepsilon}^2B_{2}\leq B_{10}<0.4398$, $B_6\geq\varepsilon B _2$, $B_1\geq B_2$ and $B_2\geq1.8815$, we find that  the left side is a decreasing function of $B_{10}$ and $B_6$. So replacing $B_{10}$ by ${\varepsilon}^2B_{2}$ and $B_6$ by $\varepsilon B _2$ we get $\phi_1=B_1+(4+4\varepsilon+{\varepsilon}^2)B_2-\sqrt{(\varepsilon)^7B_2^{11}B_1}>10.3$. Now the left side is a decreasing function of $B_2$, so replacing $B_2$ by 1.8815 we find that $\phi_1<10.3$ for $1<B_1<2.2636302$, a contradiction. Hence we must have $B_2<1.8815$.\vspace{2mm}\\
\noindent{\bf Claim(iii)} $B_3<1.5652$

Suppose $B_3\geq 1.5652$. From (3.1) we have $B_4B_5B_6<1.6524$ and $B_8B_9B_{10}\\<0.1702$, so we find that $\frac{B_3^3}{B_4B_5B_6}>2$ and $\frac{B_7^3}{B_8B_9B_{10}}\geq\frac{(\varepsilon B _3)^3}{B_8B_9B_{10}}>2$, for $B_3>1.49$.

 Applying AM-GM to inequality (2,4,4) we get $4B_1-\frac{2B_1^2}{B_2}+4B_3+4B_7-\sqrt{B_3^5B_7^5B_1B_2}>10.3$. Since $B_1\geq B_2>1.7046$, $B_7\geq\varepsilon B _3$ ~and~ $B_3\geq1.5652$, we find that left side is a decreasing function of $B_1$ and $B_7$. So we replace $B_1$ by $B_2$, ~$B_7$ by $\varepsilon B _3$ and get that
$\phi_2=2B_2+4(1+\varepsilon)B_3-\sqrt{(\varepsilon)^5B_3^{10}B_2^2}>10.3.$ But left side is a decreasing function of $B_3$, so replacing $B_3$ by 1.5652 we find that $\phi_2<10.3$ for $1.7046<B_2<1.8815$, a contradiction. Hence we must have $B_3<1.5652$.\vspace{2mm}\\
\noindent{\bf Claim(iv)} $B_1> 1.9378$

Suppose  $B_1\leq 1.9378$. Using (3.1) and that $B_3<1.5652$, $B_2>1.7046$, we find that $B_2$ is larger than each of $B_3, ~B_4,\cdots,B_{10}$. So the inequality $(1,9^*)$ holds. This gives $B_1+9(B_1)^{-1/9}>10.3$, which is not true for $B_1\leq 1.9378$. So we must have $B_1>1.9378$.\vspace{2mm}\\
\noindent{\bf Claim(v)} $B_3<1.5485$

Suppose $B_3\geq 1.5485$. We proceed as in Claim(iii) and replace $B_1$ by 1.9378 and $B_7$ by $\varepsilon B _3$  to get that
$\phi_3=4(1.9378)-\frac{2(1.9378)^2}{B_2}+4(1+\varepsilon)B_3-\sqrt{(\varepsilon)^5(1.9378)B_3^{10}B_2}>10.3.$
One easily checks that $\phi_3<10.3$ for $1.5485\leq B_3<1.5652$ and $1.7046<B_2<1.8815$. Hence we have $B_3<1.5485$.\vspace{2mm}\\
\noindent{\bf Claim(vi)} $B_1<2.0187$

Suppose $B_1\geq 2.0187$. Using (3.1) and Claims (ii), (v) we have $B_2B_3B_4<4.11$. Therefore $\frac{B_1^3}{B_2B_3B_4}>2$. As $B_5\geq \varepsilon B _1>0.9462$, we see using (3.1) that $B_5$ is larger than each of $B_6, B_7, \cdots, B_{10}$. Hence the inequality $(4,6^*)$ holds. This gives $\phi_4=4B_1-\frac{1}{2}\frac{B_1^4}{B_2B_3B_4}+6(B_1B_2B_3B_4)^{-1/6}>10.3$. Left side is an increasing function of $B_2B_3B_4$ and decreasing function of $B_1$. So we can replace $B_2B_3B_4$  by $4.11$ and $B_1$  by 2.0187 to find $\phi_4<10.3,$ a contradiction. Hence we have $B_1<2.0187$.\vspace{2mm}\\
\noindent{\bf Claim(vii)} $B_4<1.337$

Suppose $B_4\geq 1.337$, then using (3.1) we get  $\frac{B_4^3}{B_5B_6B_7}>2$. Applying AM-GM to inequality (1,2,4,2,1) we have
$B_1+4B_2-\frac{2B_2^2}{B_3}+4B_4+4B_8+B_{10}-2\sqrt{B_4^5B_8^3B_1B_2B_3B_{10}}>10.3$. Since $B_2>1.7046$, $B_3\geq\frac{3}{4}B_2$, $B_4\geq1.337$, $B_{8}\geq\varepsilon B _4$ and $B_{10}\geq\frac{2\varepsilon}{3}B_4$, we find that left side is a decreasing function of $B_2$, $B_{8}$ and $B_{10}$. So we can replace $B_2$ by 1.7046; $B_{8}$ by $\varepsilon B _4$ and $B_{10}$ by $\frac{2\varepsilon}{3}B_4$ to get
$\phi_5=B_1+4(1.7046)-\frac{2(1.7046)^2}{B_3}+(4+4\varepsilon+\frac{2\varepsilon}{3})B_4-2\sqrt{\frac{2}{3}(\varepsilon)^4(1.7046)B_4^9B_1B_3}>10.3$. Now left side is a decreasing function of $B_4$, replacing $B_4$ by $1.337$, we find that
$\phi_5 <10.3$ for $1<B_1<2.0187$ and $1<B_3<1.5485$, a contradiction. Hence we have $B_4<1.337$.\vspace{2mm}\\
\noindent {\bf Claim(viii)} $B_5<1.1492$

Suppose $B_5\geq1.1492$. Using (3.1), we get $B_6B_7B_8<0.5445.$ Therefore $\frac{B_5^3}{B_6B_7B_8}>2$. Also using \eqref{2.0a},\eqref{2.0b}, $2B_9\geq2(\varepsilon B _5)>1.077>B_{10}$. So the inequality $(4^*,4,2)$ holds, i.e. $4(\frac{1}{B_5B_6B_7B_8B_9B_{10}})^\frac{1}{4}+4B_5-\frac{1}{2}\frac{B_5^4}{B_6B_7B_8}+4B_9-\frac{2B_9^2}{B_{10}}>10.3$. Now left side is a decreasing function of $B_5$ and $B_9$. So we replace $B_5$ by 1.1492 and $B_9$ by $1.1492\varepsilon$ and get that
$\phi_6(x,B_{10})=4(\frac{1}{(\varepsilon)(1.1492)^2xB_{10}})^\frac{1}{4}+4(1+\varepsilon)(1.1492)-\frac{1}{2}\frac{(1.1492)^4}{x}
-\frac{2(1.1492\varepsilon)^2}{B_{10}}>10.3$, where $x=B_6B_7B_8$. Using \eqref{2.0a},\eqref{2.0b}   we have $x=B_6B_7B_8\geq\frac{B_5^3}{4}\geq\frac{(1.1492)^3}{4}$ and $B_{10}\geq\frac{3\varepsilon}{4}B_5\geq\frac{3\varepsilon}{4}(1.1492)$.  It can be verified that $\phi_6(x,B_{10})<10.3$ for  $\frac{(1.1492)^3}{4}\leq x<0.5445$ and $\frac{3\varepsilon}{4}(1.1492)\leq B_{10}<0.4398$, giving thereby a contradiction. Hence we must have $B_5<1.1492$.\vspace{2mm}\\
\noindent{\bf Claim(ix)} $B_2<1.766$.

 Suppose $B_2\geq 1.766.$ We have $B_3B_4B_5<2.3793$. So $\frac{B_2^3}{B_3B_4B_5}>2$. Also $B_6\geq\varepsilon B _2>0.8277$. Therefore $B_6$ is larger than each of $B_7, B_8, B_9, B_{10}$. Hence the inequality $(1,4,5^*)$ holds. This gives $B_1+4B_2-\frac{1}{2}\frac{B_2^4}{B_3B_4B_5}+5({\frac{1}{B_1B_2B_3B_4B_5}})^{\frac{1}{5}}>10.3$. Left side is an increasing function of $B_3B_4B_5$, a decreasing function of $B_2$ and an increasing function of $B_1$. One easily checks that this inequality is not true for $B_1 <2.0187$; $B_2\geq1.766$ and $B_3B_4B_5< 2.3793.$ Hence we have $B_2<1.766$.\vspace{2mm}\\
\noindent{\bf Final contradiction}\\
 As $2(B_2+B_4+B_6+B_8+B_{10})<2(1.766+1.337+0.9383+0.6597+0.4398)<10.3$, the weak inequality $(2,2,2,2,2)_w$ gives a contradiction.
 \subsection{$1.9378<B_{10}<1.9770808$}
 Here $B_1\geq B_{10}>1.9378$. And $B_2=(B_1B_3\cdots B_{10})^{-1}$\vspace{1mm}\\
  $\leq (B_1\cdot \frac{3}{4}B_2\cdot\frac{2}{3}B_2\cdot\frac{1}{2}B_2\cdot\varepsilon B_2\cdot\frac{3\varepsilon}{4} B_2\cdot\frac{2\varepsilon}{3} B_2\cdot\frac{\varepsilon}{2} B_2\cdot B_{10})^{-1}=(\frac{1}{16}{\varepsilon}^4B_2^7B_1B_{10})^{-1},$\vspace{1mm}\\ which implies $(B_2)^8\leq(\frac{1}{16}{\varepsilon}^4 (1.9378)^2)^{-1}$, i.e. $B_2<1.75076$.\\ Similarly \\
  $\begin{array}{ll}
  B_3&=(B_1B_2B_4\cdots B_{10})^{-1}\leq (\frac{3}{32}{\varepsilon}^3B_3^6B_1^2B_{10})^{-1} ;\vspace{1mm}\\
  B_4&=(B_1B_2B_3B_5\cdots B_{10})^{-1}\leq (\frac{3}{32}{\varepsilon}^2B_4^5B_1^3B_{10})^{-1};\vspace{1mm}\\
  B_6&=(B_1\cdots B_5B_7B_8B_9 B_{10})^{-1}\leq (\frac{1}{16}{\varepsilon}B_6^3B_1^5B_{10})^{-1};\vspace{1mm}\\
   B_8&=(B_1\cdots B_7B_9 B_{10})^{-1}\leq (\frac{3}{32}{\varepsilon}^3B_8B_1^7B_{10})^{-1}.\end{array}$\vspace{1mm}\\
   These respectively give $B_3<1.46138$, $B_4<1.22883$, $B_6<0.896058$ and $B_8<0.721763$.
So we have  $B_1^4B_5B_6B_7B_8B_9B_{10}=\frac{B_1^3}{B_2B_3B_4}>2$. Also $2B_5\geq2(\varepsilon B_1)>1.8166>B_6$ and $2B_7\geq2(\frac{2\varepsilon}{3}B_1)>B_8$.

 Applying AM-GM to inequality (4,2,2,1,1) we have
$4B_1+4B_5+4B_7+B_9+B_{10}-3(2B_1^5B_5^3B_7^3B_9B_{10})^{\frac{1}{3}}>10.3$. We find that left side is a decreasing function of $B_7$ and $B_5$, so can replace  $B_{7}$ by $\frac{2}{3}\varepsilon B_1$ and $B_{5}$ by $\varepsilon B _1$; then it is a decreasing function of $B_1$, so replacing $B_1$ by $B_{10}$ we have
$4(1+\varepsilon+\frac{2}{3}\varepsilon)B_{10}+B_9+B_{10}-2^{\frac{4}{3}}(\varepsilon)^2(B_{10})^{4}(B_9)^{\frac{1}{3}}>10.3$, which is not true for $(1.9378)\varepsilon^2<B_9\leq B_1<2.2636302$ and $1.9378<B_{10}<1.9770808$. Hence we get a contradiction.~~$\Box$

\section{Proof of Theorem 1 for $n=11$}
Here we have $\omega_{11}=11.62$, $B_1 \leq \gamma_{11}<2.393347$. Using \eqref{2.3}, we have $l_{11}=0.3673<B_{11}<2.1016019=m_{11}$.\\ The inequality $(10^*,1)$ gives $10.3(B_{11})^{\frac{-1}{10}}+B_{11}>11.62$. But for $0.4409\leq B_{11}\leq 2.018$ this inequality is not true. So we must have either $B_{11}<0.4409$ or $B_{11}>2.018$.
\subsection{$0.3673<B_{11}<0.4409$}
 \noindent{\bf Claim(i)}  $B_{10}<0.4692$\\
Suppose $B_{10}\geq 0.4692$, then $2B_{10}>B_{11}$, so $(9^*,2)$ holds, i.e. $9(\frac{1}{B_{10}B_{11}})^{\frac{1}{9}}+4B_{10}-\frac{2B_{10}^2}{B_{11}}>11.62$. As left side is a decreasing function of $B_{10}$, we can replace $B_{10}$ by 0.4692 and find that it is not true for $0.3673<B_{11}<0.4409$. Hence we must have $B_{10}<0.4692$.\vspace{1mm}\\
Using  \eqref{2.0a},\eqref{2.0b}   we have:
\begin{equation}\label{3.7}
\begin{array}{lll}
B_9 \leq \frac{4}{3}B_{10}<0.6256, &B_8 \leq \frac{3}{2}B_{10}<0.7038,&B_7 \leq \frac{B_{11}}{\varepsilon}<0.94063, \vspace{1mm}\\
B_6 \leq \frac{B_{10}}{\varepsilon}<1.0011,&B_5 \leq \frac{4}{3}\frac{B_{10}}{\varepsilon} < 1.3347,&B_4 \leq \frac{3}{2}\frac{B_{10}}{\varepsilon}< 1.50151,\vspace{1mm}\\
B_3\leq \frac{B_{11}}{\varepsilon^2}<2.0068,&B_2 \leq \frac{B_{10}}{\varepsilon^2}<2.13557.
\end{array}{}\end{equation}
 \noindent{\bf Claim(ii)} $B_2> 1.913$

The  inequality $(2,2,2,2,2,1)_w$ gives $ 2B_2+2B_4+2B_6+2B_8+2B_{10}+B_{11} > 11.62$. Using (4.1) we find that this inequality is not true for $B_2\leq 1.913$. So we must have $B_2> 1.913$.\vspace{2mm}\\
\noindent{\bf Claim(iii)} $B_3<1.761$

Suppose $B_3\geq 1.761$, then we have $\frac{B_3^3}{B_4B_5B_6}>2$~ and ~$\frac{B_7^3}{B_8B_9B_{10}}>\frac{(\varepsilon B_3)^3}{B_8B_9B_{10}}>2$.  Applying AM-GM to the inequality (2,4,4,1) we get $
4B_1-\frac{2B_1^2}{B_2}+4B_3+4B_7+B_{11}-\sqrt{B_3^5B_7^5B_1B_2B_{11}}>11.62$. One easily finds that it is not true for $B_1\geq B_2>1.913$, ~$B_3\geq1.761$, ~$B_7\geq \varepsilon B_3$, $B_{11}\geq \varepsilon^2B_3$,  $1.913<B_2<2.13557$ and $1.761\leq B_3<2.0068$. Hence we must have $B_3<1.761.$\vspace{2mm}\\
\noindent{\bf Claim(iv)} $B_1<2.2436$

Suppose $B_1\geq 2.2436$. As $B_2B_3B_4<2.13557\times1.761\times1.50151<5.6468$, we have $\frac{B_1^3}{B_2B_3B_4}>2$. Also $B_5\geq\varepsilon B_1>1.051$, so $B_5$ is larger than each of $B_6, B_7,\cdots,B_{11}$. Hence the inequality $(4,7^*)$ holds. This gives $4B_1-\frac{1}{2}\frac{B_1^4}{B_2B_3B_4}+7({\frac{1}{B_1B_2B_3B_4}})^{\frac{1}{7}}>11.62$. Left side is an increasing function of $B_2B_3B_4$ and decreasing function of $B_1$. One easily checks that the inequality is not true for $B_2B_3B_4<5.6468$  and $B_1\geq 2.2436$. Hence we have $B_1<2.2436$.\vspace{2mm}\\
\noindent{\bf Claim(v)} $B_4<1.4465$ and $B_2>1.9686$

Suppose $B_4\geq 1.4465$. We have $B_5B_6B_7<1.2569$ and $B_9B_{10}B_{11}<0.1295$. Therefore for $B_4>1.36$, we have $\frac{B_4^3}{B_5B_6B_7}>2$~ and ~$\frac{B_8^3}{B_9B_{10}B_{11}}>\frac{(\varepsilon B_4)^3}{B_9B_{10}B_{11}}>2$. So the inequality (1,2,4,4) holds.  Applying AM-GM to inequality (1,2,4,4), we get $
B_1+4B_2-\frac{2B_2^2}{B_3}+4B_4+4B_8-\sqrt{B_4^5B_8^5B_1B_2B_3}>11.62$. A simple calculation shows that this is not true for $B_1\geq B_2>1.913$, $B_4\geq1.4465$, ~$B_{8}\geq \varepsilon B_4$, $B_4\geq  1.4465 $, $B_1<2.2436$ and $B_3<1.761$. Hence we have $B_4<1.4465.$\vspace{1mm}

Further if $B_2\leq1.9686$, then $ 2B_2+2B_4+2B_6+2B_8+2B_{10}+B_{11} <11.62$. So the inequality $(2,2,2,2,2,1)_w$ gives a contradiction.\vspace{2mm}\\
\noindent{\bf Claim(vi)} $B_4<1.4265$ and $B_2>1.9888$

Suppose $B_4\geq 1.4265$. We proceed as in Claim(v) and get a contradiction with improved bounds on $B_2$ and $B_4$.\vspace{2mm}

\noindent{\bf Claim(vii)} $B_1<2.2056$

Suppose $B_1\geq2.2056$. We proceed as in Claim(iv) and get a contradiction with improved bounds on $B_1$, $B_2$ and $B_4$.\vspace{2mm}

\noindent{\bf Claim(viii)} $B_2<2.025$

Suppose $B_2\geq2.025$. As $B_3B_4B_5<1.761\times1.4265\times1.3347<3.3529$, we have $\frac{B_2^3}{B_3B_4B_5}>2$. Also $B_6\geq\varepsilon B_2>0.9491$, so $B_6$ is larger than each of $B_7, B_8,\cdots,B_{11}$. Hence the inequality $(1,4,6^*)$ holds, i.e. $B_1+4B_2-\frac{1}{2}\frac{B_2^4}{B_3B_4B_5}+6({\frac{1}{B_1B_2B_3B_4B_5}})^{\frac{1}{6}}>11.62$. Left side is an increasing function of $B_3B_4B_5$, a decreasing function of $B_2$ and an increasing function of $B_1$. One easily checks that this inequality is not true for $B_1 <2.2056$; $B_2\geq 2.025$ and $B_3B_4B_5< 3.3529.$ Hence we have $B_2<2.025$.\vspace{2mm}\\
\noindent{\bf Claim(ix)} $B_1<2.1669$

Suppose $B_1\geq2.1669$. We proceed as in Claim(iv) and get a contradiction with improved bounds on $B_1$, $B_2$ and $B_4$.\vspace{2mm}

\noindent{\bf Claim(x)} $B_4<1.403$ and $B_2>2.012$

Suppose $B_4\geq 1.403$. We proceed as in Claim(v) and get a contradiction with improved bounds on $B_2$ and $B_4$.\vspace{2mm}

\noindent{\bf Final Contradiction:}

As now $B_3B_4B_5<1.761\times1.403\times1.3347<3.2977$, we have $\frac{B_2^3}{B_3B_4B_5}>2$ for $B_2>2.012$. Also $B_6\geq\varepsilon B_2>0.943>$~ each of $B_7, B_8,\cdots,B_{11}$. Hence the inequality $(1,4,6^*)$ holds. Proceeding as in Claim(viii) we find that this inequality is not true for $B_1 <2.1669$; $B_2>2.012$ and $B_3B_4B_5< 3.2977,$ giving thereby a contradiction.

\subsection{$2.018<B_{11}<2.1016019$}
Here $B_1\geq B_{11}>2.018$. Therefore using \eqref{2.0a},\eqref{2.0b} we have \\$B_{10}=(B_1\cdots B_{9}B_{11})^{-1}$\vspace{1mm}\\
$\leq(B_1\cdot\frac{3}{4}B_1\cdot\frac{2}{3}B_1\cdot\frac{1}{2}B_1\cdot\varepsilon B_1\cdot\frac{3}{4}\varepsilon B_1\cdot\frac{2}{3}\varepsilon B_1\cdot\frac{1}{2}\varepsilon B_1\cdot{\varepsilon}^2B_1\cdot B_{11})^{-1}\vspace{1mm}\\=(\frac{1}{16}{\varepsilon}^6{B_1}^9B_{11})^{-1}<
(\frac{1}{16}{\varepsilon}^6(2.018)^{10})^{-1}<1.34702.$\vspace{1mm}\\
Similarly \\
$B_4=(B_1B_2B_3B_5\cdots B_{11})^{-1}\leq (\frac{1}{16}{\varepsilon}^3B_4^6B_1^3B_{11})^{-1},$  which gives $B_4<1.37661$.\vspace{2mm}\\
\noindent{\bf Claim(i)} $B_{10}<0.4402$

The inequality $(9^*,1,1)$ gives $9(\frac{1}{B_{10}B_{11}})^{\frac{1}{9}}+B_{10}+B_{11}>11.62$. But this inequality is not true for $0.4402\leq B_{10}<1.34702$ and $2.018<B_{11}<2.1016019$. Hence we must have $B_{10}<0.4402$.\vspace{1mm}\\
Now we have $B_9\leq\frac{4}{3}B_{10}<0.58694$, ~$B_{8}\leq \frac{3}{2}B_{10}<0.6603$, ~$B_7\leq2B_{10}<0.8804$ and $B_6\leq\frac{B_{10}}{\varepsilon}<0.93914$.\vspace{2mm}\\
\noindent{\bf Claim(ii)} $B_7<0.768$

Suppose $B_7\geq 0.768$. Then $\frac{B_7^3}{B_8B_9B_{10}}>2$, so $(6^*,4,1)$ holds. This gives
$\phi_7(x)=6(x)^{1/6}+4B_7-\frac{1}{2}B_7^5B_{11}x +B_{11}>11.62$, where $x=B_1B_2\dots B_6$. The function $\phi_7(x)$ has its maximum value at $x=(\frac{2}{B_7^5B_{11}})^{6/5}$. Therefore $\phi_7(x)\leq \phi_7((\frac{2}{B_7^5B_{11}})^{6/5})$,  which  is less than $11.62$ for $0.768\leq B_7<0.8804$ and $2.018<B_{11}<2.1016019$. This gives a contradiction.\vspace{1mm}\\
Now $B_5\leq\frac{3}{2}B_7<1.1521$ and $B_3\leq\frac{B_7}{\varepsilon}<1.6385.$\vspace{2mm}\\
\noindent{\bf Claim(iii)} $B_2<1.795$

Suppose $B_2\geq1.795$, then $\frac{B_2^3}{B_3B_4B_5}>2$ and $\frac{B_6^3}{B_7B_8B_9}>2$.  Applying AM-GM to the inequality (1,4,4,1,1) we get $B_1+4B_2+4B_6+B_{10}+B_{11}-\sqrt{B_2^5B_6^5B_1B_{10}B_{11}}>11.62.$  We find that  left side is a decreasing function of $B_{6}$,  so we first replace $B_6$ by $\varepsilon B_2$; then it is a decreasing function of $B_2$, so we replace $B_2$ by 1.795 and get that $\phi_8(B_{11})=B_1+4(1+\varepsilon)(1.795)+B_{10}+B_{11}-\sqrt{(\varepsilon)^{5}(1.795)^{10}B_1B_{10}B_{11}}>11.62.$ Now $\phi_8''(B_{11})>0$, so $\phi_8(B_{11})<{\mbox{max}}\{\phi_8(2.018),\phi_8(2.1016019)\}$, which can be verified to be at most 11.62 for $(\varepsilon)^2(1.795)\leq B_{10}<0.4402$ and $2.018<B_1<2.393347$, giving thereby a contradiction.

\noindent{\bf Claim(iv)} $B_5<0.98392$

 Suppose $B_5\geq0.98392$. We have  $\frac{B_1^3}{B_2B_3B_4}>2$ and $\frac{B_5^3}{B_6B_7B_8}>2.$ Also $2B_9\geq2(\varepsilon B_5)>B_{10}$.  Applying AM-GM to the inequality $(4,4,2,1)$ we get $
4B_1+4B_5+4B_9-\frac{2B_9^2}{B_{10}}+B_{11}-\sqrt{B_1^5B_5^5B_9B_{10}B_{11}}>11.62$. One can easily check that left side is a decreasing function of $B_9$  and $B_1$  so we can replace $B_9$ by $\varepsilon B_{5}$ and $B_1$ by $B_{11}$ to get $\phi_9=5B_{11}+4(1+\varepsilon)B_5-\frac{2(\varepsilon B_{5})^2}{B_{10}}-\sqrt{\varepsilon B_{11}^{6}B_5^6B_{10}}>11.62$. Now the  left side is a decreasing function of $B_5$, so replacing $B_5$ by 0.98392 we see that $\phi_9<11.62$  for $\frac{3\varepsilon}{4}(0.98392)<B_{10}<0.4409$
and $2.018<B_{11}<2.1016019$, a contradiction.\vspace{2mm}\\
\noindent{\bf Final Contradiction:}

As in Claim(iv), we have $\frac{B_1^3}{B_2B_3B_4}>2$. Also $B_5\geq\varepsilon B_1>0.9458>$~each of $B_6,B_7,\cdots,B_{10}$. Therefore the inequality $(4,6^*,1)$ holds, i.e. $\phi_{10}=4B_1-\frac{1}{2}\frac{B_1^4}{B_2B_3B_4}+6({\frac{1}{B_1B_2B_3B_4B_{11}}})^{\frac{1}{6}}+B_{11}>11.62.$ Left side is an increasing function of $B_2B_3B_4$ and $B_{11}$ and decreasing function of $B_1$. Using $B_5<0.98392$, we have $B_3\leq\frac{3}{2}B_{5}<1.47588$ and $B_4\leq\frac{4}{3}B_5<1.311894$. One easily checks that $\phi_{10}<11.62$  for $B_2B_3B_4<1.795\times1.47588\times1.311894$, $B_{11}<2.1016019$ and $B_1\geq 2.018$. Hence we have a contradiction.~$\Box$
\section{Proof of Theorem 1 for $n=12$}
\noindent Here we have $\omega_{12}=13$, $B_1 \leq \gamma_{12}<2.5217871$. Using \eqref{2.3}, we have $l_{12}=0.3376<B_{12}<2.2254706=m_{12}$ and using \eqref{2.1} we have $B_1B_2^{11}\leq\gamma_{11}^{11}$, i.e. $B_2\leq\gamma_{11}^{\frac{11}{12}}<2.2254706$. \vspace{2mm}\\ The inequality $(11^*,1)$ gives $11.62(B_{12})^{\frac{-1}{11}}+B_{12}>13$. But this is not true for $0.4165\leq B_{12}\leq 2.17$. So we must have either $B_{12}<0.4165$ or $B_{12}>2.17$.
\subsection{ $0.3376<B_{12}<0.4165$}
\noindent{\bf Claim(i)} $B_{11}<0.459$

Suppose $B_{11}\geq 0.459$, then $B_{12}\geq\frac{3}{4}B_{11}>0.34425$ and $2B_{11}>B_{12}$, so $(10^*,2)$ holds, i.e. $\phi_{11}=10.3(\frac{1}{B_{11}B_{12}})^{\frac{1}{10}}+4B_{11}-\frac{2B_{11}^2}{B_{12}}>13$. Left side is a decreasing function of $B_{11}$, so we can replace $B_{11}$ by 0.459 to find that  $\phi_{11}<13$ for $0.34425<B_{12}<0.4165$, a contradiction. Hence we have $B_{11}<0.459$.\vspace{2mm}\\
\noindent{\bf Claim(ii)} $B_{10}<0.5432$

Suppose $B_{10}\geq 0.5432$.  From \eqref{2.0a}, $B_{11}B_{12}\geq\frac{1}{2}B_{10}^2$ and $B_{10}\leq\frac{3}{2}B_{12}$. Therefore  $\frac{1}{2}(0.5432)^2\leq B_{11}B_{12}<0.1912$ and $B_{10}^2>B_{11}B_{12}$, so the inequality $(9^*,3)$ holds, i.e.
$9(\frac{1}{B_{10}B_{11}B_{12}})^{\frac{1}{9}}+4B_{10}-\frac{B_{10}^3}{B_{11}B_{12}}>13$. One easily checks that it is not true  noting that left side is a decreasing function of $B_{10}$. Hence we must have $B_{10}<0.5432$.\vspace{2mm}

\noindent{\bf Claim(iii)} $B_{9}<0.6655$

Suppose $B_9\geq0.6655$, then $\frac{B_9^3}{B_{10}B_{11}B_{12}}>2$. So the inequality $(8^*,4)$ holds. This gives $\phi_{12}(x)= 8(x)^{1/8}+4B_9-\frac{1}{2}B_9^5x>13$, where $x=B_1B_2\dots B_8$. The function $\phi_{12}(x)$ has its maximum value at $x=(\frac{2}{B_9^5})^{\frac{8}{7}}$, so $\phi_{12}(x)<\phi_{12}((\frac{2}{B_9^5})^{\frac{8}{7}})<13$ for $0.6655\leq B_9\leq\frac{3}{2}B_{11}<0.6885$. This gives a contradiction.\vspace{2mm}\\
Using \eqref{2.0a},\eqref{2.0b}   we have:
\begin{equation}\label{3.8}\begin{array}{lll}
B_8 \leq \frac{3}{2}B_{10}<0.8148,&B_7 \leq \frac{B_{11}}{\varepsilon}<0.9793,
&B_6 \leq \frac{B_{10}}{\varepsilon}<1.1589,\\ B_5 \leq \frac{B_{9}}{\varepsilon} < 1.4198,&B_4 \leq \frac{3}{2}\frac{B_{10}}{\varepsilon}< 1.7384,
&B_3\leq \frac{B_{11}}{\varepsilon^2}<2.0892.
\end{array}{}\end{equation}
\noindent{\bf Claim(iv)} $B_{2}>1.828$, $B_4>1.426$, $B_6>1.019$ and $B_8>0.715$

Suppose $B_2\leq1.828$. Then $2(B_2+B_4+B_6+B_8+B_{10}+B_{12})<2(1.828+1.7384+1.1589+0.8148+0.5432+0.4165)<13$, giving thereby a contradiction to the weak inequality $(2,2,2,2,2,2)_w$.\\
Similarly we obtain  lower bounds on $B_4, B_6$ and $B_8$ using  $(2,2,2,2,2,2)_w$.\vspace{2mm}

\noindent{\bf Claim(v)} $B_{2}>2.0299$

Suppose $B_2\leq2.0299$. Consider following two cases:\\
\noindent{\bf Case(i)} $B_{3}>B_{4}$

We have $B_3>B_4>1.426>$~each of $B_5,\cdots,B_{12}$. So the inequality $(2,10^*)$ holds, i.e. $4B_1-\frac{2B_1^2}{B_2}+10.3(\frac{1}{B_1B_2})^{\frac{1}{10}}>13.$ The left side is a decreasing function of $B_1$, so replacing $B_1$ by $B_2$ we get $2B_2+10.3(\frac{1}{B_2^2})^{\frac{1}{10}}>13$, which is not true for $B_2\leq2.0299.$\\
\noindent{\bf Case(ii)} $B_{3}\leq B_{4}$

As $B_4>1.426>$~each of $B_5,\cdots,B_{12}$, the inequality $(3,9^*)$ holds, i.e. $\phi_{13}(X)=4B_1-\frac{B_1^3}{X}+9(\frac{1}{B_1X})^{\frac{1}{9}}>13$, where $X=B_2B_3<\min\{B_1^2,(2.0299)(1.7384)\}$ $=\alpha$ say. Now $\phi_{13}(X)$ is an increasing function of $X$ for $B_1\geq B_2>1.828$ and so $\phi_{13}(X)<\phi_{13}(\alpha)$, which can be seen to be less than 13.
Hence we have $B_2>2.0299$.\vspace{2mm}\\
\noindent{\bf Claim(vi)} $B_{1}>2.17$ and $B_3<1.9517$

Using \eqref{2.1} we have $B_3\leq(\frac{\gamma_{10}^{10}}{B_1B_2})^{\frac{1}{10}}<1.9648.$ Therefore $B_2>2.0299>$~each of $B_3,\cdots,B_{12}$. So the inequality $(1,11^*)$ holds, i.e. $B_1+11.62(\frac{1}{B_{1}})^{\frac{1}{11}}>13$. But this is not true for $B_1\leq2.17$. So we must have $B_1>2.17.$ Again using \eqref{2.1} we have $B_3<(\frac{2.2636302^{10}}{2.17\times2.0299})^{\frac{1}{10}}<1.9517.$\vspace{2mm}\\
\noindent{\bf Claim(vii)} $B_4<1.646$

 Suppose $B_4\geq 1.646$. From (5.1) and Claims (i)-(iii),  we have $\frac{B_4^3}{B_5B_6B_7}>2$ ~and~ $\frac{B_8^3}{B_9B_{10}B_{11}}>\frac{(\varepsilon B_4)^3}{B_9B_{10}B_{11}}>2$.  Applying AM-GM to the  inequality (1,2,4,4,1) we get $\phi_{14}=B_1+4B_2-\frac{2B_2^2}{B_3}+4B_4+4B_8+B_{12}-\sqrt{B_4^5B_8^5B_1B_2B_{3}B_{12}}>13$. We find that left side is a decreasing function of $B_2$, ~$B_{8}$ and ~$B_{12}$. So we can replace $B_2$ by $2.0299$, ~$B_8$ by $\varepsilon B_4$ and $B_{12}$ by $\varepsilon^2B_4$. Then it turns a decreasing function of $B_4$, so can replace $B_4$ by 1.646 to find that $\phi_{14}< 13$ for $B_1<2.52178703$ and $B_{3}<1.9517$, a contradiction. Hence we have $B_4<1.646$.\vspace{2mm}

 \noindent{\bf Claim(viii)} $B_1<2.4273$

Suppose $B_1\geq2.4273$.  Consider following two cases:\\
\noindent{\bf Case(i)}  $B_5>B_6$\\
Here $B_5>$~each of $B_6,\cdots,B_{12}$ as  $B_5\geq\varepsilon B_1>1.137>$~each of $B_7,\cdots,B_{12}$. Also $B_2B_3B_4<2.2254706\times1.9517\times1.646<7.15$. So $\frac{B_1^3}{B_2B_3B_4}>2$. Hence the inequality $(4,8^*)$ holds. This gives $4B_1-\frac{1}{2}\frac{B_1^4}{B_2B_3B_4}+8(B_1B_2B_3B_4)^{-1/8}>13$. Left side is an increasing function of $B_2B_3B_4$ and decreasing function of $B_1$. So we can replace $B_2B_3B_4$ by 7.15 and $B_1$ by 2.4273 to get a contradiction.\\
\noindent{\bf Case(ii)}  $B_5\leq B_6$\\
Using (5.1) we have $B_5\leq B_6<1.1589$ and so $B_4\leq\frac{4}{3}B_5<1.5452.$ Therefore $\frac{B_2^3}{B_3B_4B_5}>2$ as $B_2>2.0299$ and $B_3<1.9517$. Also from Claim(iv), $B_6>1.019>$~each of $B_7,\cdots, B_{12}$. Hence the inequality $(1,4,7^*)$ holds. This gives $B_1+4B_2-\frac{1}{2}\frac{B_2^4}{B_3B_4B_5}+7(B_1B_2B_3B_4B_5)^{-1/7}>13.$ Left side is an increasing function of $B_3B_4B_5$ and $B_1$ and a decreasing function of $B_2$. One can check that inequality is not true for $B_3B_4B_5<1.9517\times1.5452\times1.1589$, $B_1<2.5217871$ and for $B_2>2.0299.$\\
Hence we must have $B_1<2.4273.$\vspace{2mm}\\
\noindent{\bf Claim(ix)} $B_5<1.396$

 Suppose $B_5\geq 1.396$. From (5.1), $B_6B_7B_8<0.925$ ~and~ $B_{10}B_{11}B_{12}<0.104$, so we have $\frac{B_5^3}{B_6B_7B_8}>2$ and $\frac{B_9^3}{B_{10}B_{11}B_{12}}>\frac{(\varepsilon B_5)^3}{B_{10}B_{11}B_{12}}>2$. Applying AM-GM to the inequality (1,2,1,4,4) we get $B_1+4B_2-\frac{2B_2^2}{B_3}+B_{4}+4B_5+4B_9-\sqrt{B_5^5B_9^5B_1B_2B_{3}B_{4}}>13$.  We find that left side is a decreasing function of $B_2$ and ~$B_9$. So we replace $B_2$ by $2.0299$ and ~$B_9$ by $\varepsilon B_5$. Now
   it becomes a decreasing function of $B_5$ and an increasing function of $B_1$ so replacing $B_5$ by 1.396
 and $B_1$ by 2.4273, we find that above inequality is not true for $1.522<B_3<1.9517$ and $1.426<B_4<1.646$, giving thereby a contradiction. Hence we must have  $B_5<1.396$.\vspace{2mm}\\
 \noindent{\bf Claim(x)} $B_3>1.7855$

 Suppose $B_3\leq1.7855$. We have $B_4>1.426>$~each of $B_5,B_6,\cdots,B_{12}$, hence the inequality $(1,2,9^*)$ holds. It gives $\phi_{15}=B_1+4B_2-\frac{2B_2^2}{B_3}+9(\frac{1}{B_1B_2B_3})^{\frac{1}{9}}>13.$ It is easy to check that left side of above inequality is a decreasing function of $B_2$ and an increasing function of $B_1$ and $B_3$. So replacing $B_1$ by 2.4273, $B_3$ by 1.7855 and $B_2$ by 2.0299 we get $\phi_{15}<13,$ a contradiction. Hence we have $B_3>1.7855$.\vspace{2mm}\\
  \noindent{\bf Claim(xi)} $B_2>2.0733$

  Suppose $B_2\leq2.0733$. We have $B_3>1.7855>$~each of $B_4,B_5,\cdots,B_{12}$, hence the inequality $(2,10^*)$ holds.
 It gives $\phi_{16}=4B_1-\frac{2B_1^2}{B_2}+10.3(\frac{1}{B_1B_2})^{\frac{1}{10}}>13.$ The left side is a decreasing function of $B_1$ and an increasing function of $B_2$, so replacing $B_1$ by $2.17$ and $B_2$ by 2.0733 we get $\phi_{16}<13,$ a contradiction.\vspace{2mm}\\
  \noindent{\bf Claim(xii)} $B_7<0.92$ and $B_5<1.38$

Suppose $B_7\geq0.92$. Here we have $B_4B_5B_6<2.67$ and $B_8B_9B_{10}<0.295$, so $\frac{B_3^3}{B_4B_5B_6}>2$ and $\frac{B_7^3}{B_8B_9B_{10}}>2$. Also $2B_{11}\geq2\varepsilon B_7>B_{12}$.  Applying AM-GM to the inequality (2,4,4,2) we get $\phi_{17}=4B_1-\frac{2B_1^2}{B_2}+4B_3+4B_7-\sqrt{B_3^5B_7^5B_1B_2B_{11}B_{12}}+4B_{11}-\frac{2B_{11}^2}{B_{12}}>13$. We find that left side is a decreasing function of $B_1$ and $B_{11}$. So we can replace $B_1$ by $2.17$ and $B_{11}$ by $\varepsilon B_7$. Then left side becomes a decreasing function of $B_7$ and an increasing function of $B_2$, so can replace $B_7$ by 0.92 and $B_2$ by 2.2254706 to see that $\phi_{17}<13$ for $1.7855<B_3<1.9517$ and $0.3376<B_{12}<0.4156$, a contradiction. Hence $B_7<0.92$.
Further $B_5\leq\frac{3}{2}B_7$ gives $B_5<1.38$.\vspace{2mm}

\noindent{\bf Claim(xiii)} $B_6<1.097$

Suppose $B_6\geq1.097$. Here we have $B_3B_4B_5<4.44$ and $B_7B_8B_{9}<0.5$, so $\frac{B_2^3}{B_3B_4B_5}>\frac{(2.0733)^3}{4.44}>2$ and $\frac{B_6^3}{B_7B_8B_9}>2$. Also $2B_{10}\geq2\varepsilon B_6>B_{11}$.  Applying AM-GM to the inequality (1,4,4,2,1) we get $\phi_{18}=B_1+4B_2+4B_6-\sqrt{B_2^5B_6^5B_1B_{10}B_{11}B_{12}}+4B_{10}-\frac{2B_{10}^2}{B_{11}}+B_{12}>13$. We find that left side is a decreasing function of $B_{10}$, $B_{12}$ and $B_{11}$. So we can replace $B_{10}$ by $\varepsilon B_6$ and $B_{12}$ by 0.3376 and $B_{11}$ by $\frac{3\varepsilon}{4}B_{6}$. Then left side becomes a decreasing function of $B_6$, so we can replace $B_6$ by 1.097 to find that  $\phi_{18}<13$,  for $2.17<B_1<2.4273$ and $2.0733<B_{2}<2.2254706$, a contradiction. Hence we must have $B_6<1.097$.\vspace{2mm}

\noindent{\bf Claim(xiv)} $B_5>B_6$ and $\frac{B_1^3}{B_2B_3B_4}<2$

First suppose $B_5\leq B_6$, then $B_4B_5B_6<1.646\times1.097^2<1.981$ and $\frac{B_3^3}{B_4B_5B_6}>2$. Also $B_7\geq\varepsilon B_3>0.83>$ each of $B_8,\cdots,B_{12}$. Hence the inequality $(2,4,6^*)$ holds, i.e. $4B_1-\frac{2B_1^2}{B_2}+4B_3-\frac{1}{2}\frac{B_3^4}{B_4B_5B_6}+6(\frac{1}{B_1B_2B_3B_4B_5B_6})^{\frac{1}{6}}>13$. Now the left side is a decreasing function of $B_1$ and $B_3$ as well; also it is an increasing function of $B_2$ and $B_4B_5B_6$. But one can check that this inequality is not true for $B_1>2.17$, $B_3>1.7855$, $B_2<2.2254706$ and $B_4B_5B_6<1.981$, giving thereby a contradiction.\\
Further suppose $\frac{B_1^3}{B_2B_3B_4}\geq2$, then as $B_5>B_6>1.019>$~each of $B_7,\cdots,B_{12}$, the inequality $(4,8^*)$ holds. Now working as in Case(i) of Claim(viii) we get contradiction for $B_1>2.17$ and $B_2B_3B_4<2.2254706\times1.9517\times1.646<7.14934$.\vspace{2mm}\\
 \noindent{\bf Claim(xv)} $B_3<1.9$ and $B_1<2.4056$

Suppose $B_3\geq1.9$, then for $B_4B_5B_6< 1.646\times1.38\times1.097<2.492$, $\frac{B_3^3}{B_4B_5B_6}>2$. Also $B_7\geq\varepsilon B_3>0.89>{\mbox {each~ of}}~ B_8,\cdots,B_{12}$. Hence the inequality $(2,4,6^*)$ holds. Now working as in  Claim(xiv) we get contradiction for $B_1>2.17$, $B_2<2.2254706$, $B_3>1.9$ and $B_4B_5B_6<2.492$. So $B_3<1.9$.\\ Further if $B_1\geq2.4056$, then $\frac{B_1^3}{B_2B_3B_4}>\frac{(2.4056)^3}{2.2254706\times1.9\times1.646}>2$, contradicting Claim(xiv).\vspace{2mm}\\
 \noindent{\bf Claim(xvi)} $B_4<1.58$ and $B_1<2.373$

Suppose $B_4\geq1.58$, then for $B_5B_6B_7<1.38\times1.097\times0.92<1.393$, $\frac{B_4^3}{B_5B_6B_7}>2$. Also $B_8\geq\varepsilon B_4>0.74>{\mbox {each~ of}}~ B_9,\cdots,B_{12}$. Hence the inequality $(1,2,4,5^*)$ holds, i.e. $\phi_{19}=B_1+4B_2-\frac{2B_2^2}{B_3}+4B_4-\frac{1}{2}\frac{B_4^4}{B_5B_6B_7}+5(\frac{1}{B_1B_2B_3B_4B_5B_6B_7})^{\frac{1}{5}}>13$. Left side is a decreasing function of $B_2$ and $B_4$. So we replace $B_2$ by 2.0733 and $B_4$ by 1.58. Then it becomes an increasing function of $B_1$, $B_3$ and $B_5B_6B_7$. So we replace $B_1$ by 2.4056, $B_3$ by 1.9 and $B_5B_6B_7$ by 1.393 to find that $\phi_{19}<13$, a contradiction. Further if $B_1\geq2.373$, then $\frac{B_1^3}{B_2B_3B_4}>2$, contradicting Claim(xiv).\vspace{2mm}

 \noindent{\bf Final Contradiction:}

 We have $B_3B_4B_5<1.9\times1.58\times1.38<4.15$. Therefore $\frac{B_2^3}{B_3B_4B_5}>2$. Also $B_6>1.019>$~each of $B_7,\cdots, B_{12}$. Hence the inequality $(1,4,7^*)$ holds. Now we get contradiction working as in Case(ii) of Claim(viii).
 \subsection{ $2.17<B_{12}<2.2254706$}
  Here $B_1\geq B_{12}>2.17$. Using \eqref{2.0a},\eqref{2.0b}, we have\vspace{1mm}\\
   $B_{11}=(B_1B_2\cdots B_{10}B_{12})^{-1}<(\frac{3}{64}\varepsilon^{8}B_1^{10}B_{12})^{-1}<1.8223$.\vspace{1mm}\\
   \noindent{\bf Claim(i)} Either $B_{11}<0.4307$ or $B_{11}>1.818$

  Suppose $0.4307\leq B_{11}\leq1.818$. The inequality $(10^*,1,1)$ gives $10.3(\frac{1}{B_{11}B_{12}})^{\frac{1}{10}}+B_{11}+B_{12}>13$,~ which is not true for $0.4307\leq B_{11}\leq1.818$ and $2.17<B_{12}<2.2254706$. So we must have either $B_{11}<0.4307$ or $B_{11}>1.818$.\vspace{2mm}\\
  \noindent{\bf Claim(ii)} $B_{11}<0.4307$

  Suppose $B_{11}\geq0.4307$, then using Claim(i) we have $B_{11}>1.818$. Now we have using \eqref{2.0a},\eqref{2.0b},\vspace{1mm}\\
  $B_2=(B_1B_3\cdots B_{12})^{-1}<(\frac{1}{16}\varepsilon^6B_2^8B_1B_{11}B_{12})^{-1}$. This gives $B_2<1.777$.\vspace{2mm}\\
   $B_3=(B_1B_2B_4\cdots B_{12})^{-1}<(\frac{3}{64}\varepsilon^4B_3^7B_1^2B_{11}B_{12})^{-1}$. This gives
  $B_3<1.487$\vspace{2mm}\\
    $B_4=(B_1B_2B_3B_5\cdots B_{12})^{-1}<(\frac{1}{16}\varepsilon^3B_4^6B_1^3B_{11}B_{12})^{-1}$. This gives
  $B_4<1.213$.\vspace{2mm}\\
      $B_6=(B_1\cdots B_5B_7\cdots B_{12})^{-1}<(\frac{1}{16}\varepsilon^2B_6^4B_1^5B_{11}B_{12})^{-1}$. This gives
  $B_6<0.826$.\vspace{2mm}\\
  $B_7=(B_1\cdots B_6B_8\cdots B_{12})^{-1}<(\frac{3}{64}\varepsilon^2B_7^3B_1^6B_{11}B_{12})^{-1}$. This gives
  $B_7<0.697$.\vspace{2mm}\\
    $B_8=(B_1\cdots B_7B_9\cdots B_{12})^{-1}<(\frac{1}{16}\varepsilon^3B_8^2B_1^7B_{11}B_{12})^{-1}$. This gives
  $B_8<0.559$.\vspace{2mm}\\
      $B_9=(B_1\cdots B_8B_{10}B_{11} B_{12})^{-1}<(\frac{3}{64}\varepsilon^4B_9B_1^8B_{11}B_{12})^{-1}$. This gives
  $B_9<0.478$.\vspace{2mm}\\
  $B_{10}=(B_1\cdots B_9B_{11} B_{12})^{-1}<(\frac{1}{16}\varepsilon^6B_1^9B_{11}B_{12})^{-1}<0.359$.\vspace{2mm}\\
 Therefore we have $\frac{B_1^3}{B_2B_3B_4}>2$ and $B_5\geq\varepsilon B_1>1.01>$~each of $B_6,\cdots,B_{10}$. So the inequality $(4,6^*,1,1)$ holds, i.e. $4B_1-\frac{1}{2}\frac{B_1^4}{B_2B_3B_4}+6(B_1B_2B_3B_4B_{11}B_{12})^{-1/6}$ $+B_{11}+B_{12}>13$. Now  the left side is an increasing function of $B_2B_3B_4$, $B_{11}$ and of $B_{12}$ as well. Also it is a decreasing function of $B_1$. So we replace $B_2B_3B_4$ by $1.777\times1.487\times1.213$, $B_{11}$ by 1.8223, $B_{12}$ by 2.2254706 and $B_1$ by 2.17 to arrive at a contradiction. Hence we must have $B_{11}<0.4307$.\vspace{2mm}\\
\noindent{\bf Claim(iii)} $B_{10}<0.445$

  Suppose $B_{10}\geq0.445$, then $2B_{10}>B_{11}$. So the inequality $(9^*,2,1)$ holds, i.e. $\phi_{20}=9(\frac{1}{B_{10}B_{11}B_{12}})^{\frac{1}{9}}+4B_{10}-\frac{2B_{10}^2}{B_{11}}+B_{12}>13$. Now for $B_{10}\geq0.445$, $B_{11}\geq\frac{3}{4}B_{10}$ and $B_{12}>2.2254706$, the left side is an increasing function of $B_{12}$ and a decreasing function of $B_{10}$, so replacing $B_{12}$ by 2.2254706 and $B_{10}$ by 0.445 ~we find that $\phi_{20}<13$,  for $\frac{3}{4}(0.445)< B_{11}<0.4307$, a contradiction. Hence we must have $B_{10}<0.445$.\\
  Using \eqref{2.0a},\eqref{2.0b}   we have:
\begin{equation}\label{3.9}
\begin{array}{lll}
B_9\leq\frac{4}{3}B_{10}<0.594,&B_8 \leq \frac{3}{2}B_{10}<0.67,&B_7 \leq 2B_{10}<0.89,\vspace{1mm}\\
B_6 \leq \frac{B_{10}}{\varepsilon}<0.9494,&B_5 \leq \frac{4}{3}\frac{B_{10}}{\varepsilon} <1.266,&B_4 \leq \frac{3}{2}\frac{B_{10}}{\varepsilon}< 1.4242,\vspace{1mm}\\
B_3\leq \frac{2B_{10}}{\varepsilon}<1.899,&B_{2}\leq\frac{B_{10}}{(\varepsilon)^2}<2.0255.&
\end{array}{}\end{equation}
 \noindent{\bf Claim(iv)} $B_{3}<1.62$

 Suppose $B_3\geq1.62$. From (5.2), we have $B_4B_5B_6<1.712$ and $B_8B_9B_{10}<0.178$, so $\frac{B_3^3}{B_4B_5B_6}>2$ and $\frac{B_7^3}{B_8B_9B_{10}}\geq\frac{(\varepsilon B_3)^3}{B_8B_9B_{10}}>2$.  Applying AM-GM to the inequality (2,4,4,1,1)  we get $\phi_{21}=4B_1-\frac{2B_1^2}{B_2}+4B_3+4B_7-\sqrt{B_3^5B_7^5B_1B_2B_{11}B_{12}}+B_{11}+B_{12}>13$.  We find that left side is a decreasing function of $B_1$, $B_{7}$ and ~$B_{11}$. So we can replace $B_1$ by $B_{12}$, ~$B_{7}$ by $\varepsilon B_3$ and ~$B_{11}$ by $\varepsilon^2B_3$. Then it becomes a decreasing function of $B_3$, so replacing $B_3$ by 1.62 we find that  $\phi_{21}<13,$  for $1.6275<B_2<2.0255$ and $2.17<B_{12}<2.2254706$, a contradiction. Hence we must have $B_3<1.62$.\vspace{2mm}\\
  \noindent{\bf Claim(v)} $B_{12}>2.196$

Suppose $B_{12}\leq2.196$. From(5.2), we have $B_2B_3B_4<4.674$ and $\frac{B_1^3}{B_2B_3B_4}>2$. Also $B_5\geq\varepsilon B_1>1.01>{\mbox {each~ of}}~B_6,\cdots, B_{11}$. Therefore the inequality $(4,7^*,1)$ holds, i.e.
 $\phi_{22}=4B_1-\frac{1}{2}\frac{B_1^4}{B_2B_3B_4}+7(B_1B_2B_3B_4B_{12})^{-1/7}+B_{12}>13$. Left side is an increasing function of $B_2B_3B_4$ and of $B_{12}$ as well. Also it is a decreasing function of $B_1$. So we can replace $B_2B_3B_4$ by 4.674, $B_{12}$ by 2.196 and $B_1$ by 2.17 to get  $\phi_{22}<13$, a contradiction. Hence we must have $B_{12}>2.196$.\vspace{2mm}\\
\noindent{\bf Final Contradiction:}

Now we have $B_1\geq B_{12}>2.196$. We  proceed as in Claim(v) and use $(4,7^*,1)$. Here we replace $B_2B_3B_4$ by 4.674, $B_{12}$ by 2.2254706 and $B_1$ by 2.196 to get  $\phi_{22}<13$, a contradiction.~~$\Box$
 \section{Proof of Theorem 1 for $n=13$}
\noindent Here we have $\omega_{13}=14.455765$, $B_1\leq\gamma_{13}<2.6492947$. Using \eqref{2.3}, we have $l_{13}=0.3113<B_{13}<2.348593=m_{13}$ and using \eqref{2.1} we have $B_2\leq \gamma_{12}^{\frac{12}{13}}<2.348593$.\vspace{2mm}\\
\noindent{\bf Claim(i)} $B_{13}<0.3878$

Suppose $B_{13}\geq0.3878$. The inequality $(12^*,1)$ gives $13{(B_{13})}^{\frac{-1}{12}}+B_{13}>14.455765$. But this inequality is not true for $0.3878\leq B_{13}<2.348593$. So we must have $B_{13}<0.3878$.\vspace{2mm}\\
\noindent{\bf Claim(ii)} $B_{12}<0.4353$ and $B_{11}<0.5804$

Suppose $B_{12}\geq 0.4353$, then $B_{13}\geq\frac{3}{4}B_{12}>0.3264$ and $2B_{12}>B_{13}$, so $(11^*,2)$ holds, i.e. $11.62(\frac{1}{B_{12}B_{13}})^{\frac{1}{11}}+4B_{12}-\frac{2B_{12}^2}{B_{13}}>14.455765$, which is not true for $B_{12}\geq 0.4353$ and $0.3264<B_{13}<0.3878$. Hence we have $B_{12}<0.4353$.\\
Further $B_{11}\leq\frac{4}{3}B_{12}<\frac{4}{3}(0.4353)<0.5804.$\vspace{2mm}\\
\noindent{\bf Claim(iii)} $B_{10}<0.5942$; $B_{8}<0.8913$; $B_6<1.2677$ and $B_4<1.9016$

Suppose $B_{10}\geq0.5942$, then $\frac{B_{10}^3}{B_{11}B_{12}B_{13}}>2$. So the inequality $(9^*,4)$ holds. This gives $\psi_1(x)=9x^{1/9}+4B_{10}-\frac{1}{2}B_{10}^5x>14.455765$, where $x=B_1B_2\dots B_9$. The function $\psi_1(x)$ has its maximum value at $x=(\frac{2}{B_{10}^{5}})^{\frac{9}{8}}$, so $\psi_1(x)<\psi_1((\frac{2}{B_{10}^{5}})^{\frac{9}{8}})<14.455765$ for $0.5942\leq B_{10}\leq\frac{3}{2}B_{12}<\frac{3}{2}(0.4353)<0.653.$ This gives a contradiction. Hence we have $B_{10}<0.5942$.\\
Further $B_{8}\leq\frac{3}{2}B_{10}<0.8913$, ~$B_{6}\leq\frac{B_{10}}{\varepsilon}<1.2677$ and $B_4\leq\frac{3}{2}\frac{B_{10}}{\varepsilon}<1.9016$.\vspace{2mm}\\
\noindent{\bf Claim(iv)} $B_{9}<0.74$ and $B_5<1.5788$

Suppose $B_{9}\geq0.74$, then $\frac{B_{9}^3}{B_{10}B_{11}B_{12}}>2$. So the inequality $(8^*,4,1)$ holds, i.e. $ \psi_2(x)=8x^{1/8}+4B_{9}-\frac{1}{2}B_{9}^5B_{13}x+B_{13}>14.455765$, where $x=B_1B_2\dots B_8$. The function $\psi_2(x)$ has its maximum value at $x=(\frac{2}{B_{9}^{5}B_{13}})^{\frac{8}{7}}$, so $\psi_2(x)<\psi_2((\frac{2}{B_{9}^{5}B_{13}})^{\frac{8}{7}})<14.455765$ for $0.74\leq B_{9}\leq\frac{4}{3}B_{10}<\frac{4}{3}(0.5942)<0.793$ and $\varepsilon B_9\leq B_{13}<0.3878.$ This gives a contradiction. Hence $B_9<0.74$.
and $B_5\leq\frac{B_9}{\varepsilon}<1.5788$.\vspace{2mm}\\
\noindent{\bf Claim(v)} $B_{7}<1.088$

Suppose $B_{7}\geq1.088$, then $\frac{B_{7}^3}{B_{8}B_{9}B_{10}}>2$. Also $B_{11}\geq\varepsilon B_7>0.5$ and $B_{12}B_{13}<0.4353\times0.3878<0.169$, we find $B_{11}^2>B_{12}B_{13}.$ So the inequality $(6^*,4,3)$ holds. This gives $ \psi_3(x)=6x^{1/6}+4B_{7}-\frac{1}{2}B_{7}^5B_{11}B_{12}B_{13}x+4B_{11}-\frac{B_{11}^3}{B_{12}B_{13}}>14.455765$, where $x=B_1B_2\dots B_6$. The function $\psi_3(x)$ has its maximum value at $x=(\frac{2}{B_{7}^{5}B_{11}B_{12}B_{13}})^{\frac{6}{5}}$, so $\psi_3(x)<\psi_3((\frac{2}{B_{7}^{5}B_{11}B_{12}B_{13}})^{\frac{6}{5}})=4B_{7}+5(\frac{2}{B_{7}^{5}B_{11}B_{12}B_{13}})^{\frac{1}{5}}+4B_{11}-\frac{B_{11}^3}{B_{12}B_{13}}
=\chi(B_{11})$, say. Now for  $B_{11}\geq\varepsilon B_7>0.5$ and $B_{12}B_{13}<0.169$, we find that $\chi'(B_{11})<0$, so $\chi(B_{11})\leq\chi(\varepsilon B_7)<14.455765$ for $1.088\leq B_{7}\leq\frac{3}{2}B_{9}<\frac{3}{2}(0.74)<1.11$ and $\frac{1}{2}(\varepsilon B_7)^2\leq B_{12}B_{13}<0.4353\times0.3878.$ This gives a contradiction.\vspace{2mm}\\
\noindent{\bf Claim(vi)} $B_{2}>1.942$, $B_4>1.538$, $B_6>1.103$

Suppose $B_2\leq1.942$. Then $2(B_2+B_4+B_6+B_8+B_{10}+B_{12})+B_{13}<2(1.942+1.9016+1.2677+0.8913+0.5942+0.4353)+0.3878<14.455765$, giving thereby a contradiction to the weak inequality $(2,2,2,2,2,2,1)_w$ .\\
Similarly we obtain lower bounds on $B_4$ and $B_6$ using $(2,2,2,2,2,2,1)_w$.\vspace{2mm}\\
\noindent{\bf Claim(vii)} $B_{2}>2.12$

Suppose $B_2\leq2.12$. Now we can take $B_4\geq1.72$, for if $B_4<1.72$, then $2(B_2+B_4+B_6+B_8+B_{10}+B_{12})+B_{13}<2(2.12+1.72+1.2677+0.8913+0.5942+0.4353)+0.3878<14.455765$, giving thereby a contradiction to $(2,2,2,2,2,2,1)_w$. So we have $B_4\geq1.72>$~each of $B_5,\cdots,B_{13}$. Consider following cases :\\
\noindent{\bf Case(i)} $B_{3}>B_4$

Here $B_3>B_4>$~each of $B_5,\cdots,B_{13}$. So the inequality $(2,11^*)$ holds, i.e. $4B_1-\frac{2B_1^2}{B_2}+11.62(\frac{1}{B_1B_2})^{\frac{1}{11}}>14.455765.$ The left side is a decreasing function of $B_1$, so replacing $B_1$ by $B_2$ we get $2B_2+11.62(\frac{1}{B_2^2})^{\frac{1}{11}}>14.455765$, which is not true for $1.942<B_2\leq2.185.$\\
\noindent{\bf Case(ii)} $B_{3}\leq B_4$

Here $B_{3}\leq B_4<1.9016$. As $B_4>$ each of $B_5,\cdots,B_{13}$, the inequality $(3,10^*)$ holds, i.e. $\psi_4(X)=4B_1-\frac{B_1^3}{X}+10.3(\frac{1}{B_1X})^{\frac{1}{10}}>14.455765$, where $X=B_2B_3<\alpha=\min\{B_1^2,(2.12)(1.9016)\}$. Now $\psi_4(X)$ is an increasing function of $X$  for $1.942<B_1<2.6492947$ and $X<\alpha$. Therefore $\psi_4(X)<\psi_4(\alpha)$, which can be seen to be less than 14.455765, a contradiction. \\Hence we have $B_2>2.12$.\vspace{2mm}\\
\noindent{\bf Claim(viii)} $B_{1}>2.348593$

Using \eqref{2.1}, we have $B_3<(\frac{\gamma_{11}^{11}}{B_1B_2})^{\frac{1}{11}}<2.088.$ So $B_2>2.12>$~each of $B_3,\cdots,B_{13}$, which implies that the inequality $(1,12^*)$ holds, i.e. $B_1+13(\frac{1}{B_{1}})^{\frac{1}{12}}>14.455765$. But this is not true for $B_1\leq2.348593$. So we must have $B_1>2.348593.$\vspace{2mm}\\
\noindent{\bf Claim(ix)} $B_{3}>1.865$\\
Suppose $B_3\leq1.865$. Consider following cases :\\
\noindent{\bf Case(i)} $B_{4}>B_5$\\
As $B_4>1.538>$ each of $B_6,\cdots,B_{13}$, the inequality $(1,2,10^*)$ holds, i.e. $\psi_5=B_1+4B_2-\frac{2B_2^2}{B_3}+10.3(\frac{1}{B_1B_2B_3})^{\frac{1}{10}}>14.455765.$ It is easy to check that $\psi_5$ is a decreasing function of $B_2$ and an increasing of $B_1$. So we replace $B_2$ by 2.12 and $B_1$ by 2.6492947 and find that $\psi_5<14.455765$ for $B_3\leq1.865$, a contradiction.\\
\noindent{\bf Case(ii)} $B_{4}\leq B_5$\\
We have $B_4\leq B_5<1.5788$ and $B_3<1.865$. So $\frac{B_2^3}{B_3B_4B_5}>2$. Also $B_6>1.103>$~each of $B_7,\cdots,B_{13}$. So the inequality $(1,4,8^*)$ holds, i.e. $\psi_6=B_1+4B_2-\frac{1}{2}\frac{B_2^4}{B_3B_4B_5}+8(B_1B_2B_3B_4B_5)^{-1/8}>14.455765.$ $\psi_6$ is an increasing function of $B_3B_4B_5$ and $B_1$ as well and a decreasing function of $B_2$. A simple calculation gives  $\psi_6<14.455765$ for $B_3B_4B_5<1.865\times1.5788^2$, $B_1<2.6492947$ and for $B_2>2.12,$ a contradiction.\\ Hence we have $B_3>1.865$.\vspace{2mm}\\
Using \eqref{2.1} we find $B_4<(\frac{\gamma_{10}^{10}}{B_1B_2B_3})^{\frac{1}{10}}<(\frac{2.2636302^{10}}{(2.348593)(2.12)(1.865)})^{\frac{1}{10}}<1.812.$\vspace{2mm}\\
\noindent{\bf Claim(x)} $B_{2}>2.2366$

Suppose $B_2\leq2.2366$. As $B_3>1.865>$~each of $B_4,\cdots,B_{13}$, the inequality $(2,11^*)$ holds, i.e. $4B_1-\frac{2B_1^2}{B_2}+11.62(\frac{1}{B_1B_2})^{\frac{1}{11}}>14.455765,$ which is not true for $B_1\geq B_2$ and $2.12<B_2\leq2.2366.$\vspace{2mm}\\
\noindent{\bf Claim(xi)} $B_{3}>1.917$ and $B_5>1.278$

Suppose $B_3\leq1.917$. We work as in Claim(ix) and get a contradiction.\\
 Hence we have $B_3>1.917$ and $B_5\geq\frac{2}{3}B_3>1.278$.\\
Also using \eqref{2.1}, we find $B_4<(\frac{(2.2636302)^{10}}{(2.348593)(2.2366)(1.917)})^{\frac{1}{10}}<1.7969$.\vspace{2mm}\\
\noindent{\bf Claim(xii)} $B_{1}<2.57$

Suppose $B_1\geq2.57$. Using \eqref{2.1} we get $B_2<(\frac{\gamma_{12}^{12}}{B_1})^{\frac{1}{12}}<2.3311$, $B_3<(\frac{\gamma_{11}^{11}}{B_1B_2})^{\frac{1}{11}}<2.042$ and $B_4<(\frac{\gamma_{10}^{10}}{B_1B_2B_3})^{\frac{1}{10}}<1.7808$. So $\frac{B_1^3}{B_2B_3B_4}>2$. Also $B_5>1.278>$~each of $B_6,\cdots,B_{13}$. So the inequality $(4,9^*)$ holds, i.e. $4B_1-\frac{1}{2}\frac{B_1^4}{B_2B_3B_4}+9(B_1B_2B_3B_4)^{-1/9}>14.455765.$ Left side is an increasing function of $B_2B_3B_4$ and a decreasing function of $B_1$. One can check that inequality is not true for $B_2B_3B_4<2.3311\times2.042\times1.7808$ and for $B_1\geq2.57.$ Hence we have $B_{1}<2.57$.\vspace{2mm}\\
\noindent{\bf Claim(xiii)} $B_{5}<1.522$

Suppose $B_5\geq1.522$. As $\frac{B_5^3}{B_6B_7B_8}>2$ and $\frac{B_9^3}{B_{10}B_{11}B_{12}}\geq\frac{\varepsilon B_5^3}{B_{10}B_{11}B_{12}}>2$, the inequality (2,2,4,4,1) holds. This gives $4B_1-\frac{2B_1^2}{B_2}+4B_3-\frac{2B_3^2}{B_4}+4B_5-\frac{1}{2}\frac{B_5^4}{B_6B_7B_8}+4B_9-\frac{1}{2}\frac{B_9^4}{B_{10}B_{11}B_{12}}+B_{13}>14.455765$. Applying AM-GM inequality we get $\psi_7=4B_1-\frac{2B_1^2}{B_2}+4B_3-\frac{2B_3^2}{B_4}+4B_5+4B_9+B_{13}-\sqrt{B_5^5B_9^5B_1B_2B_{3}B_{4}B_{13}}>14.455765$. As $B_1>2.348593$, $B_2>2.2366$, $B_3>1.917$, ~$1.538<B_4<1.7969$, ~$B_5\geq1.522$, ~$B_{9}\geq\varepsilon B_5$ and $\varepsilon^2B_5\leq B_{13}<0.3878$, we find that $\psi_7$ is a decreasing function of $B_1$, $B_3$, $B_9$ and $B_{13}$. So we can replace $B_1$ by $2.348593$, $B_3$ by 1.917, $B_9$ by $\varepsilon B_5$ and $B_{13}$ by $\varepsilon^2B_5$. Now left side becomes a decreasing function of $B_5$, so replacing $B_5$ by 1.522 we find that $\psi_7<14.455765$ for
 $2.2366<B_2<2.348593$ and $1.538<B_{4}<1.7969$. Hence $B_5<1.522$.\vspace{2mm}\\
 \noindent{\bf Claim(xiv)} $B_{2}<2.278$

 Suppose $B_2\geq2.278$. Using \eqref{2.1} we get $B_3<(\frac{\gamma_{11}^{11}}{B_1B_2})^{\frac{1}{11}}<2.055$. Also $B_5<1.522$. So $\frac{B_2^3}{B_3B_4B_5}>2$. Also $B_6>1.103>$~each of $B_7,\cdots,B_{13}$. So the inequality $(1,4,8^*)$ holds. Proceeding as in Case(ii) of Claim(ix) we find that the inequality is not true for $B_3B_4B_5<2.055\times1.7969\times1.522$, $B_1<2.57$ and for $B_2\geq2.278.$ Hence we have $B_2<2.278$.\vspace{2mm}\\
 \noindent{\bf Claim(xv)} $B_{3}<1.96$

 Suppose $B_3\geq1.96$, then $\frac{B_3^3}{B_4B_5B_6}>2$. Also $B_7\geq\varepsilon B_3>0.9187>$~each of $B_8,\cdots,B_{13}$. So the inequality $(2,4,7^*)$ holds, i.e. $4B_1-\frac{2B_1^2}{B_2}+4B_3-\frac{1}{2}\frac{B_3^4}{B_4B_5B_6}+7(B_1B_2B_3B_4B_5B_6)^{-1/7}>14.455765.$ Left side is an increasing function of $B_4B_5B_6$ and  $B_2$ as well. Also it is decreasing function of $B_1$ and $B_3$ as well. One can check that inequality is not true for $B_4B_5B_6<1.7969\times1.522\times1.2677$, $B_2<2.278$, $B_1>2.348593$ and for $B_3\geq1.96.$ Hence we have $B_3<1.96$.\vspace{2mm}\\
 \noindent{\bf Final Contradiction:}
 As $\frac{B_2^3}{B_3B_4B_5}>\frac{(2.2366)^3}{1.96\times1.7969\times1.522}>2$, we get contradiction using the inequality $(1,4,8^*)$ and working as in Claim(xiv).~$\Box$
\section{Proof of Theorem 1 for $n=14$}
 \noindent Here we have $\omega_{14}=15.955156$, $B_1 \leq \gamma_{14}<2.7758041$. Using \eqref{2.3}, we have $l_{14}=0.2878<B_{14}<2.4711931=m_{14}$ and using \eqref{2.1}, we have $B_2\leq\gamma_{13}^{\frac{13}{14}}<2.4711931$.\vspace{2mm}\\
 \noindent{\bf Claim(i)} $B_{14}<0.3789$

The inequality $(13^*,1)$ gives $14.455765(B_{14})^{\frac{-1}{13}}+B_{14}>15.955156$, which is not true for $0.3789\leq B_{14}<2.4711931$. So we must have $B_{14}<0.3789$.\vspace{2mm}\\
\noindent{\bf Claim(ii)} $B_{13}<0.4183$ and $B_{11}<0.62745$

Suppose $B_{13}\geq 0.4183$, then $B_{14}\geq\frac{3}{4}B_{13}>0.3137$ and $2B_{13}>B_{14}$, so $(12^*,2)$ holds, i.e. $13(\frac{1}{B_{13}B_{14}})^{\frac{1}{12}}+4B_{13}-\frac{2B_{13}^2}{B_{14}}>15.955156$, which is not true for $B_{13}\geq 0.4183$ and $0.3137<B_{14}<0.3789$. Hence we have $B_{13}<0.4183$ and $B_{11}\leq\frac{3}{2}B_{13}<0.62745.$\vspace{2mm}\\
\noindent{\bf Claim(iii)} $B_{12}<0.4994$

Suppose $B_{12}\geq 0.4994$, then $B_{12}^2>B_{13}B_{14}$, so $(11^*,3)$ holds, i.e. \\ $11.62(\frac{1}{B_{12}B_{13}B_{14}})^{\frac{1}{11}}+4B_{12}-\frac{B_{12}^3}{B_{13}B_{14}}>15.955156$. Left side is a decreasing function of $B_{12}$. Replacing $B_{12}$ by 0.4994 we get $11.62(\frac{1}{(0.4994)B_{13}B_{14}})^{\frac{1}{11}}+4(0.4994)-\frac{(0.4994)^3}{B_{13}B_{14}}>15.955156$, which is not true for $\frac{3}{4}(0.4994)\leq\frac{3}{4}B_{12}\leq B_{13}<0.4183$ and $\frac{2}{3}(0.4994)\leq\frac{2}{3}B_{12}\leq B_{14}<0.3789$. Hence we must have $B_{12}<0.4994$.\vspace{2mm}\\
\noindent{\bf Claim(iv)} $B_{10}<0.669$; $B_8<1.0035$; $B_6<1.4273$; $B_4<2.1409$

Suppose $B_{10}\geq0.669$, then $\frac{B_{10}^3}{B_{11}B_{12}B_{13}}>2$. So the inequality $(9^*,4,1)$ holds. This gives $ \psi_8(x)=9x^{1/9}+4B_{10}-\frac{1}{2}B_{10}^5B_{14}x+B_{14}>15.955156$, where $x=B_1B_2\dots B_9$. The function $\psi_8(x)$ has its maximum value at $x=(\frac{2}{B_{10}^{5}B_{14}})^{\frac{9}{8}}$, so $\psi_8(x)<\psi_8((\frac{2}{B_{10}^{5}B_{14}})^{\frac{9}{8}})=4B_{10}+8(\frac{2}{B_{10}^{5}B_{14}})^{\frac{1}{8}}+B_{14}<15.955156$ for $0.669\leq B_{10}\leq\frac{3}{2}B_{12}<\frac{3}{2}(0.4994)<0.75$ and $\varepsilon B_{10}\leq B_{14}<0.3789.$ This gives a contradiction.\\
Further $B_8\leq\frac{3}{2}B_{10}<1.0035$, $B_6\leq\frac{B_{10}}{\varepsilon}<1.4273$ and $B_4\leq\frac{3}{2}\frac{B_{10}}{\varepsilon}<2.1409$.\vspace{2mm}\\
\noindent{\bf Claim(v)} $B_{9}<0.8233$; $B_7<1.23495$; $B_5<1.7565$

Suppose $B_{9}\geq0.8233$, then $\frac{B_{9}^3}{B_{10}B_{11}B_{12}}>2$. Also $2B_{13}\geq2(\varepsilon B_9)>0.77>B_{14}.$ So the inequality $(8^*,4,2)$ holds. This gives $ \psi_9(x)=8x^{1/8}+4B_{9}-\frac{1}{2}B_{9}^5B_{13}B_{14}x+4B_{13}-\frac{2B_{13}^2}{B_{14}}>15.955156$, where $x=B_1B_2\dots B_8$. The function $\psi_9(x)$ has its maximum value at $x=(\frac{2}{B_{9}^{5}B_{13}B_{14}})^{\frac{8}{7}}$, so $\psi_9(x)<\psi_9((\frac{2}{B_{9}^{5}B_{13}B_{14}})^{\frac{8}{7}})=4B_{9}+7(\frac{2}{B_{9}^{5}B_{13}B_{14}})^{\frac{1}{7}}+4B_{13}-\frac{2B_{13}^2}{B_{14}}$, which is a decreasing function of $B_{13}$, so for $B_{13}\geq\varepsilon B_9$ we have $\psi_9(x)<4(1+\varepsilon)B_{9}+7(\frac{2}{\varepsilon B_{9}^{6}B_{14}})^{\frac{1}{7}}-\frac{2(\varepsilon B_{9})^2}{B_{14}}<15.955156$ for $0.8233\leq B_{9}\leq\frac{4}{3}B_{10}<\frac{4}{3}(0.669)<0.892$ and $\frac{3}{4}\varepsilon B_9\leq B_{14}<0.3789.$ This gives a contradiction.\\
Further $B_7\leq\frac{3}{2}B_{9}<1.23495$ and $B_5\leq\frac{B_{9}}{\varepsilon}<1.7565$.\vspace{2mm}\\
\noindent{\bf Claim(vi)} $B_{2}>1.858$, $B_4>1.599$, $B_6>1.182$

Suppose $B_2\leq1.858$. Then $2(B_2+B_4+B_6+B_8+B_{10}+B_{12}+B_{14})<2(1.858+2.1409+1.4273+1.0035+0.669+0.4994+0.3789)<15.955156$, giving thereby a contradiction to the weak inequality $(2,2,2,2,2,2,2)_w$.\\
Similarly we obtain lower bounds on $B_4$ and $B_6$ using $(2,2,2,2,2,2,2)_w$. \vspace{2mm}\\
\noindent{\bf Claim(vii)} $B_{2}>2.24$

Suppose $B_2\leq2.24$. Here we can take $B_4\geq1.759$, because if $B_4<1.759$, then $2(B_2+B_4+B_6+B_8+B_{10}+B_{12}+B_{14})<2(2.24+1.759+1.4273+1.0035+0.669+0.4994+0.3789)<15.955156$, giving thereby a contradiction to $(2,2,2,2,2,2,2)_w$. So we have $B_4\geq1.759>$~each of $B_5,\cdots,B_{14}$. Also using \eqref{2.1} we find $B_4\leq(\frac{\gamma_{11}^{11}}{B_2B_3})^{\frac{1}{12}}\leq(\frac{\gamma_{11}^{11}}{(3/4)B_2^2})^{\frac{1}{12}}<
(\frac{(2.393347)^{11}}{(3/4)(1.858)^2})^{\frac{1}{12}}<2.05586.$ Consider following cases:\\
\noindent{\bf Case(i)} $B_{3}>B_4$

Here $B_3>B_4>$~each of $B_5,\cdots,B_{14}$. So the inequality $(2,12^*)$ holds, i.e. $4B_1-\frac{2B_1^2}{B_2}+13(\frac{1}{B_1B_2})^{\frac{1}{12}}>15.955156.$ The left side is a decreasing function of $B_1$, so replacing $B_1$ by $B_2$ we get $2B_2+13(\frac{1}{B_2^2})^{\frac{1}{12}}>15.955156$, which is not true for $1.858<B_2\leq2.33.$\\
\noindent{\bf Case(ii)} $B_{3}\leq B_4$

As $B_4>$~each of $B_5,\cdots,B_{14}$, the inequality $(3,11^*)$ holds, i.e. $\psi_{10}(X)=4B_1-\frac{B_1^3}{X}+11.62(\frac{1}{B_1X})^{\frac{1}{11}}>15.955156,$ where $X=B_2B_3<\alpha=\min\{B_1^2,\\(2.24)(2.05586)\}$. Now $\psi_{10}'(X)=\frac{B_1^3}{X^2}(1-\frac{11.62}{11}(\frac{X^{10}}{B_1^{34}})^{\frac{1}{11}})>0$ for $B_1\geq B_2>1.858$ and $X<\alpha$. Therefore we have $\psi_{10}(X)<\psi_{10}(\alpha)<15.955156$.\\
Hence we have $B_2>2.24$.\vspace{2mm}\\
\noindent{\bf Claim(viii)} $B_{1}>2.471194$

Using \eqref{2.1} we find $B_3<2.20463.$ So $B_2>2.24>$~each of $B_3,\cdots,B_{14}$, which implies that the inequality $(1,13^*)$ holds, i.e. $B_1+14.455765(\frac{1}{B_{1}})^{\frac{1}{13}}>15.955156$. But this is not true for $B_1\leq2.471194$. Therefore $B_1>2.471194.$\vspace{2mm}\\
\noindent{\bf Claim(ix)} $B_{3}>1.998$

Suppose $B_3\leq1.998$. We have $B_4>1.599>$ each of $B_6,\cdots,B_{14}$. Consider following cases:\\
\noindent{\bf Case(i)} $B_{4}>B_5$\\
Here the inequality $(1,2,11^*)$ holds, i.e. $B_1+4B_2-\frac{2B_2^2}{B_3}+11.62(\frac{1}{B_1B_2B_3})^{\frac{1}{11}}>15.955156.$ It is easy to check that the left side is a decreasing function of $B_2$ and an increasing of $B_1$ and $B_3$ as well. So we replace $B_2$ by 2.24, $B_1$ by 2.7758041 and $B_3$ by 1.998 and get a contradiction.\\
\noindent{\bf Case(ii)} $B_{4}\leq B_5$ and $B_1\leq2.66$\\
Here $B_4\leq B_5<1.7565$ and $B_5\geq B_4>$~each of $B_6,\cdots,B_{14}$. So the inequality $(1,2,1,10^*)$ holds, i.e. $B_1+4B_2-\frac{2B_2^2}{B_3}+B_4+10.3(B_1B_2B_3B_4)^{-1/10}>15.955156.$ Left side is an increasing function of $B_4$, $B_3$ and $B_1$. Also it is a decreasing function of $B_2$. So we replace $B_4$ by 1.7565, $B_3$ by 1.998, $B_1$ by 2.66 and $B_2$ by 2.24 to get a contradiction.\\
\noindent{\bf Case(iii)} $B_{4}\leq B_5$ and $B_1>2.66$\\
Here $\frac{B_1^3}{B_2B_3B_4}>2$ and $B_5\geq B_4>$~each of $B_6,\cdots,B_{14}$. So the inequality $(4,10^*)$ holds, i.e. $4B_1-\frac{1}{2}\frac{B_1^4}{B_2B_3B_4}+10.3(B_1B_2B_3B_4)^{-1/10}>15.955156.$ Left side is an increasing function of $B_2B_3B_4$ and a decreasing function of $B_1$. But one can check that the inequality is not true for $B_1>2.66$ and $B_2B_3B_4<2.471194\times1.998\times1.7565$.\\
Hence we must have $B_3>1.998$.\vspace{2mm}\\
Using \eqref{2.1} we find $B_3<2.18665$, $B_4<1.92362$ and $B_5<1.69842$.\vspace{2mm}\\
\noindent{\bf Claim(x)} $B_{2}<2.4266$

 Suppose $B_2\geq2.4266$. We have $B_6>1.182>$~each of $B_8,\cdots,B_{14}$. Consider following cases:\\
 \noindent{\bf Case(i)} $B_{6}>B_7$\\
Here $B_6>$~each of $B_7,\cdots,B_{14}$ and $\frac{B_2^3}{B_3B_4B_5}>2$. So the inequality $(1,4,9^*)$ holds, i.e. $B_1+4B_2-\frac{1}{2}\frac{B_2^4}{B_3B_4B_5}+9(B_1B_2B_3B_4B_5)^{-1/9}>15.955156.$ Left side is an increasing function of $B_3B_4B_5$ and $B_1$ as well and a decreasing function of $B_2$. One can check that inequality is not true for $B_3B_4B_5<2.18665\times1.92362\times1.69842$, $B_1<2.7758041$ and for $B_2\geq2.4266.$\\
 \noindent{\bf Case(ii)} $B_{6}\leq B_7$\\
 Here $B_6\leq B_7<1.23495$, $B_5\leq\frac{4}{3}B_6<1.6466$ and $B_4\leq\frac{3}{2}B_6<1.8525$. So $\frac{B_3^3}{B_4B_5B_6}>2$. Also $B_7\geq B_6>$~each of $B_8,\cdots,B_{14}$. So the inequality $(2,4,8^*)$ holds, i.e. $\psi_{11}=4B_1-\frac{2B_1^2}{B_2}+4B_3-\frac{1}{2}\frac{B_3^4}{B_4B_5B_6}+8(B_1B_2B_3B_4B_5B_6)^{-1/8}>15.955156.$ Left side is an increasing function of $B_4B_5B_6$. Also it is decreasing function of $B_1$ and $B_3$ as well. We replace $B_4B_5B_6$ by $1.8525\times1.6466\times1.23495$, $B_3$ by 1.998 and $B_1$ by $B_2$ to find that $\psi_{11}<15.955156$  for $2.24<B_2<2.4711931$, a contradiction.\vspace{2mm}\\
 \noindent{\bf Claim(xi)} $B_2>2.372$ and $B_1<2.635$

 First suppose $B_2\leq2.372$. As $B_3>1.998>$~each of $B_4,\cdots,B_{14}$, the inequality $(2,12^*)$ holds, i.e. $\psi_{12}(B_1)=4B_1-\frac{2B_1^2}{B_2}+13(\frac{1}{B_1B_2})^{\frac{1}{12}}>15.955156.$ $\psi_{12}(B_1)$ is a decreasing function of $B_1$, so for $B_1>2.471194$, $\psi_{12}(B_1)<\psi_{12}(2.471194)<15.955156$ for $2.24<B_2\leq2.372.$ So we must have $B_2>2.372$.\\
 Further if $B_1\geq2.635$, then $\psi_{12}(B_1)\leq\psi_{12}(2.635)<15.955156$ for $2.24<B_2<2.4266.$ So we must have $B_1<2.635$.\vspace{2mm}\\
  \noindent{\bf Claim(xii)} $B_{3}<2.097$

 Suppose $B_3\geq2.097$. We have $B_7\geq\varepsilon B_3>0.9829>$~each of $B_9,\cdots,B_{14}$. Consider following cases:\\
  \noindent{\bf Case(i)} $B_{7}>B_8$\\
 Here $B_7>$~each of $B_8,\cdots,B_{14}$. Using \eqref{2.1} we have $B_4<1.90524$ and $B_5<1.68058$. Also $B_6<1.4273$. So $\frac{B_3^3}{B_4B_5B_6}>2$, which implies that the inequality $(2,4,8^*)$ holds. Now proceeding as in Case(ii) of Claim(x) we get contradiction for $B_4B_5B_6<1.90524\times1.68058\times1.4273$, $B_2<2.4266$, $B_1>2.471194$ and for $B_3\geq2.097.$\\
  \noindent{\bf Case(ii)} $B_{7}\leq B_8$\\
 Here $B_7\leq B_8<1.0035$, $B_6\leq\frac{4}{3}B_7<1.338$ and $B_5\leq\frac{3}{2}B_7<1.50525$. So $\frac{B_4^3}{B_5B_6B_7}>2$. Also $B_8\geq B_7>$~each of $B_9,\cdots,B_{14}$. So the inequality $(1,2,4,7^*)$ holds. This gives $B_1+4B_2-\frac{2B_2^2}{B_3}+4B_4-\frac{1}{2}\frac{B_4^4}{B_5B_6B_7}+7(B_1\cdots B_7)^{-1/7}>15.955156.$ It is easy to check that left side of this inequality is a decreasing function of $B_2$ and $B_4$. Also it is increasing function of $B_5B_6B_7$, $B_1$ and $B_3$. But the inequality is not true for $B_2>2.372$, $B_4>1.599$, $B_1<2.635$, $B_3<2.18665$ and $B_5B_6B_7<1.50525\times1.338\times1.0035$.\vspace{2mm}\\
 \noindent{\bf Claim (xiii)} $B_{4}<B_5$

 Suppose $B_4\geq B_5$. Also $B_4>1.599>$~each of $B_6,\cdots,B_{14}$. So the inequality $(1,2,11^*)$ holds, i.e. $B_1+4B_2-\frac{2B_2^2}{B_3}+11.62(\frac{1}{B_1B_2B_3})^{\frac{1}{11}}>15.955156.$ It is easy to check that the left side is a decreasing function of $B_2$ and an increasing of $B_1$ and $B_3$ as well. So we replace $B_2$ by 2.372, $B_1$ by 2.635 and $B_3$ by 2.097 and get a contradiction. Hence we have $B_4<B_5$.\vspace{2mm}\\
   \noindent{\bf Final Contradiction:}

   Using \eqref{2.1} we have $B_5<1.68872$. Also $B_4<B_5$, i.e. $1.599<B_4<B_5<1.68872$. It also gives $B_5>1.599>$~each of $B_6,\cdots,B_{14}$. So the inequality $(2,2,10^*)$ holds, i.e. $4B_1-\frac{2B_1^2}{B_2}+4B_3-\frac{2B_3^2}{B_4}+10.3(\frac{1}{B_1B_2B_3B_4})^{\frac{1}{10}}>15.955156.$ Now the left side is a decreasing function of $B_1$ and $B_3$. Also it is an increasing function of $B_2$ and $B_4$. But this inequality is not true for $B_1>2.471194$, $B_3>1.998$, $B_2<2.4266$ and $B_4<1.68872$.~$\Box$
\section{Proof of Theorem 1 for $n=15$}
 \noindent Here we have $\omega_{15}=17.498499$, $B_1 \leq \gamma_{15}<2.90147763$. Using \eqref{2.3} we have $l_{15}=0.2667<B_{15}<2.5931615=m_{15}$. Using \eqref{2.1} we have $B_2\leq\gamma_{14}^{\frac{14}{15}}<2.5931615$ and $B_4^{12}\leq\frac{\gamma_{12}^{12}}{B_1B_2B_3}\leq\frac{\gamma_{12}^{12}}{B_1.(3/4)B_1.(2/3)B_1}\leq\frac{2\gamma_{12}^{12}}{B_4^3}$, i.e. $B_4\leq(2\gamma_{12}^{12})^{\frac{1}{15}}<(2(2.52178703)^{12})^{\frac{1}{15}}<2.195$. \vspace{2mm}\\
 \noindent{\bf Claim(i)}  $B_{15}<0.3705$

 The inequality $(14^*,1)$ gives $15.955156(B_{15})^{\frac{-1}{14}}+B_{15}>17.498499$, which is not true for $0.3705\leq B_{15}< 2.5931615$. So we must have $B_{15}<0.3705$.\vspace{2mm}\\
 \noindent{\bf Claim(ii)} $B_{14}<0.4101$

Suppose $B_{14}\geq 0.4101$, then $B_{15}\geq\frac{3}{4}B_{14}>0.3075$ and $2B_{14}>B_{15}$, so $(13^*,2)$ holds, i.e. $14.455765(\frac{1}{B_{14}B_{15}})^{\frac{1}{13}}+4B_{14}-\frac{2B_{14}^2}{B_{15}}>17.498499$. But this is not true for $B_{14}\geq 0.4101$ and $0.3075<B_{15}<0.3705$. Hence $B_{14}<0.4101$.\vspace{2mm}\\
\noindent{\bf Claim(iii)} $B_{13}<0.4813$

Suppose $B_{13}\geq 0.4813$, then $B_{13}^2>B_{14}B_{15}$, so $(12^*,3)$ holds, i.e.\\ $13(\frac{1}{B_{13}B_{14}B_{15}})^{\frac{1}{12}}+4B_{13}-\frac{B_{13}^3}{B_{14}B_{15}}>17.498499$. Left side is a decreasing function of $B_{13}$. Replacing $B_{13}$ by 0.4813 we get $13(\frac{1}{(0.4813)B_{14}B_{15}})^{\frac{1}{12}}+4(0.4813)-\frac{(0.4813)^3}{B_{14}B_{15}}>17.498499$, which is not true for $\frac{3}{4}(0.4813)\leq B_{14}<0.4101$ and $\frac{2}{3}(0.4813)\leq B_{15}<0.3705$. Hence we must have $B_{13}<0.4813$.\vspace{2mm}\\
\noindent{\bf Claim(iv)} $B_{12}<0.5553$

Suppose $B_{12}\geq0.5553$, then $\frac{B_{12}^3}{B_{13}B_{14}B_{15}}>2$. So the inequality $(11^*,4)$ holds. This gives $ \psi_{13}(x)=11.62(x)^{1/11}+4B_{12}-\frac{1}{2}B_{12}^5x>17.498499$, where $x=B_1B_2\dots B_{11}$. The function $\psi_{13}(x)$ has its maximum value at $x=(\frac{11.62}{11}\times\frac{2}{B_{12}^{5}})^{\frac{11}{10}}$, so $\psi_{13}(x)<\psi_{13}((\frac{11.62}{11}\times\frac{2}{B_{12}^{5}})^{\frac{11}{10}})<17.498499$ for $0.5553\leq B_{12}\leq\frac{3}{2}B_{14}<0.616$.\vspace{2mm}\\
\noindent{\bf Claim(v)} $B_{11}<0.6471$; $B_{9}<0.97065$; $B_7<1.38054$

Suppose $B_{11}\geq0.6471$, then $\frac{B_{11}^3}{B_{12}B_{13}B_{14}}>2$. So the inequality $(10^*,4,1)$ holds. This gives $ \psi_{14}(x)=10.3x^{1/10}+4B_{11}-\frac{1}{2}B_{11}^5B_{15}x+B_{15}>17.498499$, where $x=B_1B_2\dots B_{10}$. The function $\psi_{14}(x)$ has its maximum value at $x=(\frac{10.3}{10}\times\frac{2}{B_{11}^{5}B_{15}})^{\frac{10}{9}}$, so $\psi_{14}(x)<\psi_{14}((\frac{10.3}{10}\times\frac{2}{B_{11}^{5}B_{15}})^{\frac{10}{9}})<17.498499$ for $0.6471\leq B_{11}\leq\frac{3}{2}B_{13}<0.72195$ and $\varepsilon B_{11}\leq B_{15}<0.3705.$ This gives a contradiction.\\
Further $B_9\leq\frac{3}{2}B_{11}<0.97065$ and $B_7\leq\frac{B_{11}}{\varepsilon}<1.38054$.\vspace{2mm}\\
\noindent{\bf Claim(vi)} $B_{10}<0.7525$; $B_8<1.12875$; $B_6<1.6055$

Suppose $B_{10}\geq0.7525$, then $\frac{B_{10}^3}{B_{11}B_{12}B_{13}}>\frac{(0.7525)^3}{0.6471\times0.5553\times0.4813}>2$. Also $2B_{15}\geq2(\varepsilon B_{10})>0.7>B_{15}$. So the inequality $(9^*,4,2)$ holds. This gives $ \psi_{15}(x)=9x^{1/9}+4B_{10}-\frac{1}{2}B_{10}^5B_{14}B_{15}x+4B_{14}-\frac{2B_{14}^2}{B_{15}}>17.498499$, where $x=B_1B_2\dots B_9$. The function $\psi_{15}(x)$ has its maximum value at $x=(\frac{2}{B_{10}^{5}B_{14}B_{15}})^{\frac{9}{8}}$, so $\psi_{15}(x)<\psi_{15}((\frac{2}{B_{10}^{5}B_{14}B_{15}})^{\frac{9}{8}})=4B_{10}+8(\frac{2}{B_{10}^{5}B_{14}B_{15}})^{\frac{1}{8}}+4B_{14}-\frac{2B_{14}^2}{B_{15}}$, which is a decreasing function of $B_{14}$, so for $B_{14}\geq\varepsilon B_{10}$ we have $\psi_{15}(x)<4(1+\varepsilon)B_{10}+8(\frac{2}{\varepsilon B_{10}^{6}B_{15}})^{\frac{1}{8}}-\frac{2(\varepsilon B_{10})^2}{B_{15}}$, which is less than 17.498499 for $0.7525\leq B_{10}<\frac{3}{2}(0.5553)<0.833$ and $\frac{3}{4}\varepsilon B_{10}\leq B_{15}<0.3705.$ This gives a contradiction.\\
Further $B_8\leq\frac{3}{2}B_{10}<1.12875$ and $B_6\leq\frac{B_{10}}{\varepsilon}<1.6055$.\vspace{2mm}\\
\noindent{\bf Claim(vii)} $B_{2}>1.916$, $B_4>1.644$, $B_6>1.249$ and $B_5<1.96235$

Suppose $B_2\leq1.916$. Then $2(B_2+B_4+B_6+B_8+B_{10}+B_{12}+B_{14})+B_{15}<2(1.916+2.195+1.6055+1.12875+0.7525+0.5553+0.4101)+0.3705<17.498499$, a contradiction to the weak inequality $(2,2,2,2,2,2,2,1)_w$.\\
Similarly we obtain lower bounds on $B_4$ and $B_6$ using $(2,2,2,2,2,2,2,1)_w$.
Also using \eqref{2.1}, we find  $B_5\leq(\frac{\gamma_{11}^{11}}{B_2B_3B_4})^{\frac{1}{12}}<(\frac{(2.393347)^{11}}{(3/4)B_2^2B_4})^{\frac{1}{12}}<1.96235$.\vspace{2mm}\\
\noindent{\bf Claim(viii)} $B_3>1.6$; $B_5<1.9449$; $B_4<2.15477$

Suppose $B_3\leq1.6$, then $B_2\leq\frac{4}{3}B_3<2.134$.  Therefore $2B_2+B_3+2(B_5+B_7+B_9+B_{11}+B_{13}+B_{15})<2(2.134)+1.6+2(1.96235+1.38054+0.97065+0.6471+0.4813+0.3705)<17.498499$, a contradiction to $(2,1,2,2,2,2,2,2)_w$.\\
Now using \eqref{2.1} we find $B_5<1.9449$ and $B_4<2.1547$.\vspace{2mm}\\
\noindent{\bf Claim(ix)} $B_{2}>2.384$

Suppose $B_2\leq2.384$. Consider following cases:\\
\noindent{\bf Case(i)} $B_{3}>2.15477$

Here $B_3$ is larger than each of $B_4,\cdots,B_{15}$. So the inequality $(2,13^*)$ holds, i.e. $4B_1-\frac{2B_1^2}{B_2}+14.455765(\frac{1}{B_1B_2})^{\frac{1}{13}}>17.498499,$ which is not true for $B_1\geq B_2$ and $1.916<B_2\leq2.45.$\\
\noindent{\bf Case(ii)} $B_{3}\leq 2.15477$ and $B_4\geq B_5$

 As $B_4\geq B_5$ and $B_4>1.644>$ each of $B_6,\cdots,B_{15}$, the inequality $(3,12^*)$ holds, i.e. $\psi_{16}(X)=4B_1-\frac{B_1^3}{X}+13(\frac{1}{B_1X})^{\frac{1}{12}}>17.498499,$ where $X=B_2B_3<\alpha=\min\{B_1^2,(2.384)(2.15477)\}$. Now $\psi_{16}'(X)=\frac{B_1^3}{X^2}-\frac{13}{12}(\frac{1}{B_1X^{13}}))^{\frac{1}{12}}=\frac{B_1^3}{X^2}(1-\frac{13}{12}(\frac{X^{11}}{B_1^{37}})^{\frac{1}{12}}>0$ for $B_1\geq B_2>1.916$ and $X<\alpha$. Therefore $\psi_{16}(X)<\psi_{16}(\alpha)$, which can be seen to be less than 17.498499.\\
\noindent{\bf Case(iii)} $B_{3}\leq 2.15477$ and $B_4<B_5$

Here $B_5>B_4>$ each of $B_6,\cdots,B_{15}$. Therefore the inequality $(2,2,11^*)$ holds, i.e. $4B_1-\frac{2B_1^2}{B_2}+4B_3-\frac{2B_3^2}{B_4}+11.62(\frac{1}{B_1B_2B_3B_4})^{\frac{1}{11}}>17.498499.$ Left side is a decreasing function of $B_1$, so replacing $B_1$ by $B_2$ we get $2B_2+4B_3-\frac{2B_3^2}{B_4}+11.62(\frac{1}{B_2^2B_3B_4})^{\frac{1}{11}}>17.498499.$ Now the left side is an increasing function of $B_2$, so replacing $B_2$ by 2.384 we find that the inequality is not true for $1.6<B_3\leq2.15477$ and $1.644<B_4<B_5<1.9449$.\\
Hence we must have $B_2>2.384$.\vspace{2mm}\\
\noindent{\bf Claim(x)} $B_{1}>2.5931615$ and $B_5<1.8575$

Using \eqref{2.1} we find $B_3<2.318.$ So $B_2>2.384>$~each of $B_3,\cdots,B_{15}$, which implies that the inequality $(1,14^*)$ holds, i.e. $B_1+15.955156(\frac{1}{B_{1}})^{\frac{1}{14}}>17.498499$. But this is not true for $B_1\leq2.5931615$. So we must have $B_1>2.5931615.$\\
Now using \eqref{2.1} we find $B_5<1.8575$.\vspace{2mm}\\
\noindent{\bf Claim(xi)} $B_{3}>2.133$\\
Suppose $B_3\leq2.133$. Consider following cases:\vspace{1mm}\\
\noindent{\bf Case(I)} $B_{4}>B_5$\vspace{1mm}\\
As $B_4>1.644>$ each of $B_6,\cdots,B_{15}$, the inequality $(1,2,12^*)$ holds, i.e. $B_1+4B_2-\frac{2B_2^2}{B_3}+13(\frac{1}{B_1B_2B_3})^{\frac{1}{12}}>17.498499.$ The left side is a decreasing function of $B_2$ and an increasing of $B_1$ and $B_3$ as well. So we replace $B_2$ by 2.384, $B_1$ by 2.90147763 and $B_3$ by 2.133 and get a contradiction.\vspace{1mm}\\
\noindent{\bf Case(II)} $B_{4}\leq B_5$ and $B_1\leq2.8$\vspace{1mm}\\
We have $1.644<B_4\leq B_5<1.8575$. It implies $B_5>1.644>$~each of $B_6,\cdots,B_{15}$. So the inequality $(1,2,1,11^*)$ holds, i.e. $\psi_{17}=B_1+4B_2-\frac{2B_2^2}{B_3}+B_4+11.62(B_1B_2B_3B_4)^{-1/11}>17.498499.$ $\psi_{17}$ is an increasing function of $B_4$ and $B_3$. Also it is a decreasing function of $B_2$. So we replace $B_4$ by 1.8575, $B_3$ by 2.133, $B_2$ by 2.384  and find that $\psi_{17}<17.498499$ for $2.5931615<B_1\leq2.8$.\vspace{1mm}\\
\noindent{\bf Case(III)} $B_{4}\leq B_5$ and $B_1>2.8$\vspace{1mm}\\
Here again $B_4\leq B_5<1.8575$ and $B_5\geq B_4>1.644>$~each of $B_6,\cdots,B_{14}$. Also $\frac{B_1^3}{B_2B_3B_4}>\frac{(2.8)^3}{2.594\times2.133\times1.8575}>2$. So the inequality $(4,11^*)$ holds, i.e. $4B_1-\frac{1}{2}\frac{B_1^4}{B_2B_3B_4}+11.62(B_1B_2B_3B_4)^{-1/11}>17.498499.$ Left side is an increasing function of $B_2B_3B_4$ and a decreasing function of $B_1$. But one can check that the inequality is not true for $B_1>2.8$ and $B_2B_3B_4<2.594\times2.133\times1.8575$.\\
Hence we must have $B_3>2.133$.\vspace{2mm}\\
Using \eqref{2.1} we find $B_3<2.30289$, $B_4<2.03406$ and $B_5<1.80946$.\vspace{2mm}\\
\noindent{\bf Claim(xii)} $B_2>2.49$

 Suppose $B_2\leq2.49$. As $B_3>2.133>$~each of $B_4,\cdots,B_{15}$, the inequality $(2,13^*)$ holds, i.e. $\psi_{18}(B_1)=4B_1-\frac{2B_1^2}{B_2}+14.455765(\frac{1}{B_1B_2})^{\frac{1}{13}}>17.498449.$ But $\psi_{18}(B_1)<\psi_{18}(2.5931615)<17.498499$ for $2.384<B_2\leq2.49.$ Hence $B_2>2.49$.\vspace{2mm}\\
\noindent{\bf Claim(xiii)} $B_4>1.883$; $B_5<1.78611$

Suppose $B_4\leq1.883$. We have $B_5\geq\frac{2}{3}B_3>\frac{2}{3}(2.133)>1.422>$~each of $B_7,\cdots,B_{15}$. Consider following cases:\vspace{1mm}\\
\noindent{\bf Case(I)} $B_5>B_6$\vspace{1mm}\\
Here the inequality $(2,2,11^*)$ holds, i.e. $\psi_{19}=4B_1-\frac{2B_1^2}{B_2}+4B_3-\frac{2B_3^2}{B_4}+11.62(\frac{1}{B_1B_2B_3B_4})^{\frac{1}{11}}>17.498499$. Now the left side is a decreasing function of $B_1$ and $B_3$, so we replace $B_1$ by $B_2$ and $B_3$ by 2.133 and then find that $\psi_{19}<17.498499$ for $2.49<B_2<2.5931615$ and $1.644<B_4<1.883$.\vspace{1mm}\\
\noindent{\bf Case(II)} $B_5\leq B_6$\vspace{1mm}\\
Here $B_6\geq B_5>$~each of $B_7,\cdots,B_{15}$. Also $B_5\leq B_6<1.6055$. Therefore $\frac{B_2^3}{B_3B_4B_5}>2$. So the inequality $(1,4,10^*)$ holds, which gives $B_1+4B_2-\frac{1}{2}\frac{B_2^4}{B_3B_4B_5}+10.3(\frac{1}{B_1B_2B_3B_4B_5})^{\frac{1}{10}}>17.498499$. Now the left side is a decreasing function of $B_2$ and an increasing function of $B_1$ and $B_3B_4B_5$. One can check that inequality is not true for $B_2>2.49$, $B_1<2.90147763$ and $B_3B_4B_5<2.30289\times1.883\times1.6055$. Hence we must have $B_4>1.883$.\vspace{1mm}\\ Further using \eqref{2.1} we find $B_5<1.781$.\vspace{2mm}\\
\noindent{\bf Claim(xiv)} $B_2<2.5585$

Suppose $B_2\geq2.5585$. We have $B_6>1.249>$~each of $B_8,\cdots,B_{15}$. Consider following cases:\vspace{1mm}\\
\noindent{\bf Case(I)} $B_6>B_7$\\
We have $\frac{B_2^3}{B_3B_4B_5}>2$, so the inequality $(1,4,10^*)$ holds. Working as in Case(II) of Claim(xiii) we here get contradiction for $B_2\geq2.5585$, $B_1<2.90147763$ and $B_3B_4B_5<2.30289\times2.03406\times1.781$.\vspace{1mm}\\
\noindent{\bf Case(II)} $B_6\leq B_7$ and $B_3\leq2.2$\\
As $B_4>1.883>$~each of $B_5,\cdots,B_{15}$, the inequality $(1,2,12^*)$ holds. Here working as in Case(I) of Claim(xi) we get contradiction for $B_2>2.5585$, $B_1<2.90147763$ and $B_3\leq2.2$.\vspace{1mm}\\
\noindent{\bf Case(III)} $B_6\leq B_7$ and $B_3>2.2$\vspace{1mm}\\
Here $B_6\leq B_7<1.38054$. Also $B_7\geq B_6>$~each of $B_8,\cdots,B_{15}$ and $\frac{B_3^3}{B_4B_5B_6}>\frac{(2.2)^3}{2.03406\times1.781\times1.38054}>2$, so the inequality $(2,4,9^*)$ holds. This gives $\psi_{20}=4B_1-\frac{2B_1^2}{B_2}+4B_3-\frac{1}{2}\frac{B_3^4}{B_4B_5B_6}+9(\frac{1}{B_1B_2B_3B_4B_5B_6})^{\frac{1}{9}}>17.498499$. Here $\psi_{20}$ is a decreasing function of $B_3$ and $B_1$ and an increasing function of $B_4B_5B_6$. So we replace $B_1$ by $B_2$, $B_3$ by 2.2 and $B_4B_5B_6$ by $2.03406\times1.78727\times1.38054$ and then find that $\psi_{20}<17.498499$ for $B_2<2.5931615$.
Hence we must have $B_2<2.5585$.\vspace{2mm}\\
\noindent{\bf Claim(xv)} $B_1<2.797$

Suppose $B_1\geq2.797$. As $B_3>2.133>$~each of $B_4,\cdots,B_{15}$, the inequality $(2,13^*)$ holds. Here working as in Case(I) of Claim(ix) we get contradiction for $B_1\geq2.797$ and $B_2<2.5585$.\vspace{2mm}\\
\noindent{\bf Claim(xvi)} $B_3<2.2398$

Suppose $B_3\geq2.2398$, then $B_7\geq\varepsilon B_3>1.0498>$~each of $B_9,\cdots,B_{15}$. Using \eqref{2.1} we get $B_4<2.0185$; $B_5<1.7724$ and $B_6<1.564$. Now consider following cases:\vspace{1mm}\\
\noindent{\bf Case(I)} $B_7>B_8$\vspace{1mm}\\
Here $B_7>$~each of $B_8,\cdots,B_{15}$. Also $\frac{B_3^3}{B_4B_5B_6}>2$. Therefore the inequality $(2,4,9^*)$ holds. Now working as in Case(III) of Claim(xiv) we get contradiction for $B_1>B_2$, $B_3\geq2.2398$, $B_4B_5B_6<2.0185\times1.7724\times1.564$ and $B_2<2.5585$.\vspace{1mm}\\
\noindent{\bf Case(II)} $B_7\leq B_8$\vspace{1mm}\\
Here $B_8\geq B_7>$~each of $B_9,\cdots,B_{15}$ and $B_7\leq B_8<1.12875$. Also $\frac{B_4^3}{B_5B_6B_7}>2$. So the inequality $(1,2,4,8^*)$ holds. This gives $B_1+4B_2-\frac{2B_2^2}{B_3}+4B_4-\frac{1}{2}\frac{B_4^4}{B_5B_6B_7}+8(\frac{1}{B_1B_2B_3B_4B_5B_6B_7})^{\frac{1}{8}}>17.498499$. Now the left side is a decreasing function of $B_2$ and $B_4$. This is also an increasing function of $B_1$, $B_3$ and $B_5B_6B_7$. One can check that inequality is not true for $B_2>2.49$, $B_4>1.883$, $B_1<2.797$, $B_3<2.30289$ and $B_5B_6B_7<1.7724\times1.564\times1.12875$. Hence we must have $B_3<2.2398$.\vspace{2mm}\\
\noindent{\bf Claim(xvii)} $B_1>2.705$

Suppose $B_1\leq2.705$. As $B_4>1.883>$~each of $B_5,\cdots,B_{15}$, the inequality $(1,2,12^*)$ holds. Now working as in Case(I) of Claim(xi) we get contradiction for $B_2>2.49$, $B_1\leq2.705$ and $B_3<2.2398$.\vspace{2mm}\\
\noindent{\bf Claim(xviii)} $B_2>2.525$

Suppose $B_2\leq2.525$. As $B_3>2.133>$~each of $B_4,\cdots,B_{15}$, the inequality $(2,13^*)$ holds. Now working as in Claim(xv) we get contradiction for $B_1>2.705$ and $B_2\leq2.525$.\vspace{2mm}\\
\noindent{\bf Final Contradiction:}

Using \eqref{2.1} we find $B_4<2.0173$ and $B_5\leq(\frac{\gamma_{11}^{11}}{B_1B_2B_3B_4})^{\frac{1}{11}}<1.7712$. Also $B_6>1.249>$~each of $B_8,\cdots,B_{15}$. Consider following cases:\vspace{1mm}\\
\noindent{\bf Case(I)} $B_6>B_7$\vspace{1mm}\\
We have $\frac{B_2^3}{B_3B_4B_5}>2$, so the inequality $(1,4,10^*)$ holds. Working as in Case(II) of Claim(xiii) we here get contradiction for $B_2>2.525$, $B_1<2.797$ and $B_3B_4B_5<2.2398\times2.0173\times1.7712$.\vspace{1mm}\\
\noindent{\bf Case(II)} $B_6\leq B_7$ and $B_3\leq2.22$\vspace{1mm}\\
As $B_4>1.883>$~each of $B_5,\cdots,B_{15}$, the inequality $(1,2,12^*)$ holds, i.e. $B_1+4B_2-\frac{2B_2^2}{B_3}+13(\frac{1}{B_1B_2B_3})^{\frac{1}{12}}>17.498499$. Here working as in Case(I) of Claim(xi) we get contradiction for $B_2>2.525$, $B_1<2.797$ and $B_3\leq2.22$.\vspace{1mm}\\
\noindent{\bf Case(III)} $B_6\leq B_7$ and $B_3>2.22$\vspace{1mm}\\
Here $B_7\geq B_6>$~each of $B_8,\cdots,B_{15}$ and $B_6\leq B_7<1.38054$. Also $\frac{B_3^3}{B_4B_5B_6}>2$, so the inequality $(2,4,9^*)$ holds. Here working as in Case(III) of Claim(xiv) we get contradiction for $B_1>2.705$, $B_3>2.22$, $B_4B_5B_6< 2.0173\times1.7712\times1.38054$ and  $B_2<2.5585$.~~$\Box$
\section{Proof of Theorem 1 for $16\leq n\leq33$}

 In addition of Lemmas 1-7, we shall use the following lemmas also:\vspace{2mm}

   \noindent {\bf Lemma 8.} For any integer $s, ~ 1\leq s\leq n-1$
\begin{equation*} B_1B_2 \cdots B_s \geq \left
\{ \begin{array}{lll}\frac{(0.46873)^{k(2k-2)}B_1^s}{4^k}& {\mbox {if}}&s=4k\\\frac{(0.46873)^{k(2k-1)}B_1^s}{4^k }& {\mbox {if}}&s=4k+1\\ \frac{3(0.46873)^{k(2k)}B_1^s}{4\times4^k }& {\mbox {if}}&s=4k+2\\\frac{(0.46873)^{k(2k+1)}B_1^s}{2\times 4^k }& {\mbox {if}}&s=4k+3. \end{array}\right.{}\end{equation*}

\noindent  This is Lemma 8 of Hans-Gill et al\cite{HRS4}.\vspace{2mm}

\noindent {\bf Lemma 9.} For any integer $s, ~ 1\leq s\leq n-1$
$$ B_1B_2 \dots B_s \geq \left
\{ \begin{array}{lll}\frac{(0.46873)^{k(2k-2)}}{4^k B_{n}^{n-s}}& {\mbox {if}}&n-s=4k\\\frac{(0.46873)^{k(2k-1)}}{4^k B_{n}^{n-s}}& {\mbox {if}}&n-s=4k+1\\\frac{3(0.46873)^{k(2k)}}{4\times4^k B_{n}^{n-s}}& {\mbox {if}}&n-s=4k+2\\\frac{(0.46873)^{k(2k+1)}}{2\times 4^k B_{n}^{n-s}}& {\mbox {if}}&n-s=4k+3. \end{array}\right.{}$$
This is Lemma 9 of Hans-Gill et al\cite{HRS4}.\vspace{2mm}

 If for some $s$, $1\leq s \leq n-1$, $B_{s+1}\geq   B_{s+j}$ for all $j,~ 2\leq j\leq n-s$ then the inequality $(s^*,(n-s)^*)$ holds, which gives
\begin{equation}\label{4.3.4}\phi_{s,n-s}(B_1B_2\dots B_s)=\omega_s (B_1B_2\dots B_s)^{\frac{1}{s}}+\omega_{n-s}\left({\frac{1}{B_1B_2\dots B_s}}\right)^{\frac{1}{n-s}}> \omega_n.\end{equation} Let $\lambda_s^{(n)}$  be the larger of  the lower bounds of $B_1B_2\dots B_s$ given in Lemmas 8 and 9 and let $~\mu_s^{(n)}$ be the upper bound of $B_1B_2\dots B_s$ given in \eqref{2.2}.\\

\noindent The following two lemmas  are respectively the Lemmas 11 \& 12 of Hans-Gill et al \cite{HRS4}.\vspace{2mm}\\
\noindent{\bf Lemma 10 :}
If for some $s,~ 1\leq s \leq n-1$, $\phi_{s,n-s}(\lambda_s^{(n)})\leq \omega_n$ and $ \phi_{s,n-s}(\mu_s^{(n)})\leq \omega_n$,  then we must have $B_{s+1}<max\{B_{s+2},\cdots, B_n\}$.\\

\noindent{\bf Proof :} Suppose $B_{s+1}\geq max\{B_{s+2},\cdots, B_n\}$, then the inequality $(s^*, (n-s)^*)$ holds which gives the inequality \eqref{4.3.4}.
Also the function  $\phi_{s,n-s}(x)$ has maximum at one of the end points of the interval in which $x$ lies. For $x= B_1B_2\dots B_s$ and $\lambda_s^{(n)}\leq B_1B_2\dots B_s\leq \mu_s^{(n)}$, this  contradicts the hypothesis.\\

\noindent {\bf Remark 2:} In all the cases we find that $max\left\{\phi_{s,n-s}(\lambda_s^{(n)}),~\phi_{s,n-s}(\mu_s^{(n)})\right\}$ is $\phi_{s,n-s}(\mu_s^{(n)})$.\\

\noindent{\bf Lemma 11 :}
Suppose that for some  $s,~ 1\leq s \leq n-1$, $\phi_{s,n-s}(\lambda_s^{(n)})\leq \omega_n$, but $ \phi_{s,n-s}(\mu_s^{(n)})> \omega_n$. Let a real number  $\sigma_s^{(n)}$ be such that  $\lambda_s^{(n)}< \sigma_s^{(n)}< \mu_s^{(n)}$ and $\phi_{s,n-s}(\sigma_s^{(n)})\leq \omega_n$.  \vspace{2mm} \\
 $~~~~~~~~~$ (i) $~~~$ If  $~~B_1B_2\dots B_s < \sigma_s^{(n)}$~ then ~
 $B_{s+1}<{max}\{B_{s+2},\cdots, B_n\}$, \vspace{2mm}\\
$~~~~~~~~~$ (ii) $~~~$ If $~~B_1B_2\dots B_s \geq \sigma_s^{(n)}$~ then ~
 $B_{s+1}\leq \frac{\mu_{s+1}^{(n)}}{\sigma_s^{(n)}}.$
 \vspace{2mm}\\
 \noindent{\bf Proof :} In Case (i) Lemma 10  gives the result. In Case (ii), since  $B_{s+1}=\frac{B_1\dots B_{s+1}}{B_1\dots B_s}$ we use Lemma 6 to get the desired result.\vspace{2mm}\\
 In Sections 9.1-9.18, we give proof of Theorem 1 for $16\leq n\leq33$. For the proof, we obtain upper bounds on $B_s$ for different $s$ by applying various inequalities and using Lemmas 10 and 11 and get a final contradiction
by applying weak inequality $(2,2,\cdots,2)$ for even $n$ and $(1,2,2,\cdots,2)$ or $(2,1,2,\cdots,2)$ or $(2,2,\cdots,2,1)$ for odd $n$.
\subsection{$n=16$}

Here we have $\omega_{16}=19.285$, $B_1\leq\gamma_{16}<3.0263937$. Using \eqref{2.3}, we have $l_{16}=0.2477<B_{16}<2.7145981=m_{16}$. Using \eqref{2.1} we have
\begin{equation}\label{2.1.1}
\begin{array}{lll}
 B_2\leq\gamma_{15}^{\frac{15}{16}}<2.7145981,&~B_4\leq(2\gamma_{13}^{13})^{\frac{1}{16}}<2.3047,
~B_5\leq(4\gamma_{12}^{12})^{\frac{1}{16}}<2.1823.
\end{array}
\end{equation}
\noindent{\bf Claim(i)} $B_{16}<0.2928$

Suppose $B_{16}\geq 0.2928$. The inequality $(15^*,1)$ gives $17.498499(B_{16})^{\frac{-1}{15}}+B_{16}>19.285$. But this is not true for $0.2928\leq B_{16}< 2.7145981$. So we have $B_{16}<0.2928$. \vspace{2mm}\\
Using \eqref{2.0a},\eqref{2.0b} and Claim(i) we find
\begin{equation}\label{2.1.2}max\{B_6,B_7,\cdots,B_{16}\}<B_{6}\leq\frac{3}{2}\frac{B_{16}}{\varepsilon^2}<1.9991\end{equation}
\noindent{\bf Claim(ii)} $B_{3}<max\{B_4,B_5,\cdots,B_{16}\}<2.3047$

 If $B_3\geq max\{B_4,B_5,\cdots,B_{16}\}$, then the inequality $(2,14^*)$ holds, i.e. $4B_1-\frac{2B_1^2}{B_2}+\omega_{14}(\frac{1}{B_1B_2})^{\frac{1}{14}}>19.285$. The left side is a decreasing function of $B_1$, so replacing $B_1$ by $B_2$ we get $2B_2+(15.955156)(\frac{1}{B_2^2})^{\frac{1}{14}}>19.285$, which is not true for $\frac{3}{4}<\frac{3}{4}B_1\leq B_2<2.7145981$. Hence $B_{3}<max\{B_4,B_5,\cdots,B_{16}\}$, which is $<2.3047$ (using \eqref{2.1.1},\eqref{2.1.2}). So we have $B_3<2.3047$.\vspace{2mm}\\
 \noindent{\bf Claim(iii)} $B_{2}<max\{B_3,B_4,\cdots,B_{16}\}<2.3047$

As $1=\lambda_1^{(16)}<B_1<\mu_1^{(16)}=3.0263937$, we find that \\ ${max}\{\phi_{1,15}(\lambda_1^{(16)}),\phi_{1,15}(\mu_1^{(16)})\}=\phi_{1,15}(\mu_1^{(16)})<\omega_{16}$, therefore using Lemma 10 we have $B_{2}<max\{B_3,B_4,\cdots,B_{16}\}$, which is $<2.3047$ (by  Claim (ii) ). \vspace{2mm}\\
\noindent{\bf Claim(iv)} $B_{4}<max\{B_5,B_6,\cdots,B_{16}\}<2.1823$; $B_{2}, B_3<2.1823$

As $\frac{1}{2}=\lambda_3^{(16)}<B_1B_2B_3<\mu_3^{(16)}=(3.0264)(2.3047)^2$ and \\${\max}\{\phi_{3,13}(\lambda_3^{(16)}), \phi_{3,13}(\mu_3^{(16)})\}=\phi_{3,13}(\mu_3^{(16)})<\omega_{16}$, therefore using Lemma 10 we have $B_{4}<max\{B_5,B_6,\cdots,B_{16}\}$, which is $<2.1823$ (using \eqref{2.1.1},\eqref{2.1.2}). So $max\{B_5,\cdots,B_{16}\}<2.1823$. Now by Claims (ii) and (iii), we have $B_3<2.1823$ and $B_2<2.1823$ respectively.\vspace{2mm}\\
\noindent{\bf Claim(v)} $B_{2}, B_3, B_4, B_5 <\frac{3}{2}\frac{B_{16}}{\varepsilon^2}$

As $\frac{1}{4}=\lambda_4^{(16)}<B_1B_2B_3B_4<\mu_4^{(16)}=(3.0264)(2.1823)^3$ and \\${\max}\{\phi_{4,12}(\lambda_4^{(16)}), \phi_{4,12}(\mu_4^{(16)})\}=\phi_{4,12}(\mu_4^{(16)})<\omega_{16}$, therefore using Lemma 10 we have $B_{5}<max\{B_6,B_7,\cdots,B_{16}\}\leq\frac{3}{2}\frac{B_{16}}{\varepsilon^2}$. So $B_5<\frac{3}{2}\frac{B_{16}}{\varepsilon^2}$. Now again using Claims (ii), (iii) and (iv)  we get each of $B_5$, $B_4$, $B_3$ and $B_2$ is $<\frac{3}{2}\frac{B_{16}}{\varepsilon^2}$.\vspace{2mm}\\
\noindent{\bf Final Contradiction:}

\noindent Now \vspace{2mm}\\
$\begin{array}{l} 2B_2+2B_4+2B_6+2B_8+2B_{10}+2B_{12}+2B_{14}+2B_{16}\vspace{2mm}\\ \leq 2(\frac{3/2}{\varepsilon^2}+\frac{3/2}{\varepsilon^2}+\frac{3/2}{\varepsilon^2}+\frac{1}{\varepsilon^2}
+\frac{3/2}{\varepsilon}+\frac{1}{\varepsilon}+\frac{3}{2}+1)B_{16}<19.285 {\rm ~~for~} B_{16}<0.2928.
\end{array}$\vspace{2mm}\\ This gives a contradiction to the weak inequality $(2,2,2,2,2,2,2,2)_w$.~$\Box$
\subsection{$n=17$}
Here we have $\omega_{17}=21.101$, $B_1\leq\gamma_{17}<3.1506793$. Using \eqref{2.3} we have $l_{17}=0.2306<B_{17}<2.8355395=m_{17}$. Using \eqref{2.1} we have
\begin{equation}\label{2.2.1}
\begin{array}{lll}B_2\leq\gamma_{16}^{\frac{16}{17}}<2.8355395,&~ B_4\leq(2\gamma_{14}^{14})^{\frac{1}{17}}<2.4147,& ~B_5\leq(4\gamma_{13}^{13})^{\frac{1}{17}}<2.2856 \\ B_6\leq(\frac{4\gamma_{12}^{12}}{\varepsilon})^{\frac{1}{17}}<2.1794.&
\end{array}
\end{equation}
\noindent{\bf Claim(i)} $B_{17}<0.298$

Suppose $B_{17}\geq0.298$. The inequality $(16^*,1)$ gives $19.285(B_{17})^{\frac{-1}{16}}+B_{17}>21.101$. But this is not true for $0.298\leq B_{17}<2.8355395$. So we must have $B_{17}<0.298$.\vspace{2mm}\\
\noindent{\bf Claim(ii)} $B_{16}<0.3328$

Suppose $B_{16}\geq0.3328$. The inequality $(15^*,2)$ gives $17.498499(B_{16}B_{17})^{\frac{-1}{15}}+4B_{16}-\frac{2B_{16}^2}{B_{17}}>21.101$. But this is not true for $0.3328\leq B_{16}\leq \frac{4}{3}B_{17}<\frac{4}{3}(0.298)$ and $\frac{3}{4}(0.3328)<B_{17}<0.298$. So we must have $B_{16}<0.3328$.\vspace{2mm}\\
Using \eqref{2.0a},\eqref{2.0b} and Claims(i), (ii) we find
\begin{equation}\label{2.2.2}max\{B_7,B_8,\cdots,B_{17}\}<B_7<\frac{4}{3}\frac{B_{16}}{\varepsilon^2}<2.0197\end{equation}

\noindent{\bf Claim(iii)} $B_{3}<max\{B_4,B_5,\cdots,B_{17}\}<2.4147$

Suppose $B_3\geq max\{B_4,B_5,\cdots,B_{17}\}$, then the inequality $(2,15^*)$ gives $4B_1-\frac{2B_1^2}{B_2}+\omega_{15}(\frac{1}{B_1B_2})^{\frac{1}{15}}>21.101$. The left side is a decreasing function of $B_1$, so replacing $B_1$ by $B_2$ we get $2B_2+\omega_{15}(\frac{1}{B_2^2})^{\frac{1}{15}}>21.101$, which is not true for $\frac{3}{4}<B_2<2.8356$. Hence $B_{3}<max\{B_4,B_5,\cdots,B_{17}\}$, which is $<2.4147$~(using \eqref{2.2.1},\eqref{2.2.2}). So $B_3<2.4147$.\vspace{2mm}\\
\noindent{\bf Claim(iv)} $B_{2}<max\{B_3,B_4,\cdots,B_{17}\}<2.4147$

As $1=\lambda_1^{(17)}<B_1<\mu_1^{(17)}=3.1507$, we find that ${max}\{\phi_{1,16}(\lambda_1^{(17)}),\phi_{1,16}(\mu_1^{(17)})\}=\phi_{1,16}(\mu_1^{(17)})<\omega_{17}$, therefore using Lemma 10 we have $B_{2}<max\{B_3,B_4,\cdots,B_{17}\}$, which is $<2.4147$~(using Claim (iii)). So $B_2<2.4147$.\vspace{2mm}\\
\noindent{\bf Claim(v)} $B_{4}<max\{B_5,B_6,\cdots,B_{17}\}<2.2856$ ; $B_3, B_2<2.2856$

Now $\frac{1}{2}=\lambda_3^{(17)}<B_1B_2B_3<\mu_3^{(17)}=(3.1507)(2.4147)^2$ and ${\max}\{\phi_{3,14}(\lambda_3^{(17)}),\\ \phi_{3,14}(\mu_3^{(17)})\}=\phi_{3,14}(\mu_3^{(17)})<\omega_{17}$, therefore using Lemma 10 we have $B_{4}<max\{B_5,B_6,\cdots,B_{17}\}$, which is $<2.2856$(using \eqref{2.2.1},\eqref{2.2.2}). So $B_4<2.2856$. Now again using Claims (iii), (iv) respectively, we have $B_3<2.2856$ and $B_2<2.2856$.\vspace{2mm}\\
\noindent{\bf Claim(vi)} $B_{5}<max\{B_6,B_7,\cdots,B_{17}\}<2.1794$; ~$B_4, B_3, B_2<2.1794$

As $\frac{1}{4}=\lambda_4^{(17)}<B_1B_2B_3B_4<\mu_4^{(17)}=(3.1507)(2.2856)^3$ and ${\max}\{\phi_{4,13}(\lambda_4^{(17)}),\\ \phi_{4,13}(\mu_4^{(17)})\}=\phi_{4,13}(\mu_4^{(17)})<\omega_{17}$, therefore using Lemma 10 we have $B_{5}<max\{B_6,B_7,\cdots,B_{17}\}$, which is $<2.1794$ (using \eqref{2.2.1},\eqref{2.2.2}). So $B_5<2.1794$. Now again using Claims (iii), (iv) and (v), we get each of $B_4$, $B_3$ and $B_2$ is $<2.1794$.\vspace{2mm}\\
\noindent{\bf Claim(vii)} $B_6, B_5$, $B_4$, $B_3$ and $B_2$ $<\frac{4}{3}\frac{B_{16}}{\varepsilon^2}$

As $\frac{\varepsilon}{4}=\lambda_5^{(17)}<B_1B_2B_3B_4B_5<\mu_5^{(17)}=(3.1507)(2.1794)^4$ and ${\max}\{\phi_{5,12}(\lambda_5^{(17)}),\\ \phi_{5,12}(\mu_5^{(17)})\}=\phi_{5,12}(\mu_5^{(17)})<\omega_{17}$, therefore using Lemma 10 we have $B_{6}<max\{B_7,B_8,\cdots,B_{17}\}<\frac{4}{3}\frac{B_{16}}{\varepsilon^2}$. Now again using Claims (iii)-(vi), we get each of $B_5$, $B_4$, $B_3$ and $B_2$ is $<\frac{4}{3}\frac{B_{16}}{\varepsilon^2}$.\vspace{3mm}\\
\noindent{\bf Final Contradiction}

 Now $2B_2+2B_4+2B_6+2B_8+2B_{10}+2B_{12}+2B_{14}+2B_{16}+B_{17}<2(\frac{4/3}{\varepsilon^2}+\frac{4/3}{\varepsilon^2}+\frac{4/3}{\varepsilon^2}+\frac{1}{\varepsilon^2}
+\frac{3/2}{\varepsilon}+\frac{1}{\varepsilon}+\frac{3}{2}+1)B_{16}+B_{17}<21.101$ for $B_{16}<0.3328$ and $B_{17}<0.298$. This gives a contradiction to the weak inequality $(2,2,2,2,2,2,2,2,1)_w$.~$\Box$
\subsection{$n=18$}
Here we have $\omega_{18}=22.955$, $B_1\leq\gamma_{18}<3.2743307$. Using \eqref{2.3}, we have $l_{18}=0.2150<B_{18}<2.9560725=m_{18}$.\vspace{2mm}\\
\noindent{\bf Claim(i)} $B_{18}<0.3$

Suppose $B_{18}\geq0.3$. The inequality $(17^*,1)$ gives $21.101(B_{18})^{\frac{-1}{17}}+B_{18}>22.955$. But this is not true for $0.3\leq B_{18}<2.9560725$. So we must have $B_{18}<0.3$.\vspace{2mm}\\
\noindent{\bf Claim(ii)} $B_{16}<0.385$

The inequality $(15^*,3)$ gives $17.498499({B_{16}B_{17}B_{18}})^{\frac{-1}{15}}+4B_{16}-\frac{B_{16}^3}{B_{17}B_{18}}>22.955$. But this is not true for $0.385\leq B_{16}\leq \frac{3}{2}B_{18}<0.45$ and $\frac{1}{2}B_{16}^2\leq B_{17}B_{18}\leq\frac{4}{3}B_{18}^2<\frac{4}{3}(0.3)^2$. So we must have $B_{16}<0.385$.\vspace{2mm}\\
\noindent{\bf Claim(iii)} $B_9, B_{10}, B_{11}\leq\frac{4}{3}\frac{B_{16}}{\varepsilon}$ and $B_2, B_3, B_4, B_5<\frac{4}{3}\frac{B_{16}}{(\varepsilon)^2}$

 Using \eqref{2.2} and Lemmas 8, 9 we have\vspace{2mm}\\
 $~~~~~~~1=\lambda_1^{(18)}<B_1<\mu_1^{(18)}=3.2743307,$\\
 $~~~~~~~\frac{3}{4}=\lambda_2^{(18)}<B_1B_2<\mu_2^{(18)}=(1.1519977)^{16},$\\
 $~~~~~~~\frac{1}{2}=\lambda_3^{(18)}<B_1B_2B_3<\mu_3^{(18)}=(1.2402616)^{15},$\\
 $~~~~~~~\frac{3\varepsilon^8}{4^3B_{18}^{10}}<\lambda_8^{(18)}<B_1B_2\cdots B_8\leq \mu_8^{(18)}<(1.8635658)^{10},$\\
 $~~~~~~~\frac{\varepsilon^6}{4^2B_{18}^{9}}<\lambda_9^{(18)}<B_1B_2\cdots B_9\leq \mu_9^{(18)}<(2.0406455)^9.$\vspace{2mm}\\
We find that $max\{\phi_{s,n-s}(\lambda_s^{(18)}),\phi_{s,n-s}(\mu_s^{(18)}), {\mbox {for}}~s=1,2,3,8,9\}=\phi_{1,n-1}(\mu_s^{(18)})$, which is $<22.955$. Therefore using Lemma 10 we have \begin{equation}\label{2.8}B_{i}<max\{B_{i+1},B_{i+2},\cdots,B_{18}\}, {\mbox {for}}~i=2,3,4,9,10\end{equation}
Using \eqref{2.0a},\eqref{2.0b}, we have $B_{11}\leq\frac{4}{3}\frac{B_{16}}{\varepsilon}$. So $B_9$, $B_{10}$ $<\frac{4}{3}\frac{B_{16}}{\varepsilon}$ (by \eqref{2.8}). Further $B_5\leq\frac{B_9}{\varepsilon}\leq\frac{4}{3}\frac{B_{16}}{(\varepsilon)^2}$. So each of $B_2, B_3$ and $B_4$ is $<\frac{4}{3}\frac{B_{16}}{(\varepsilon)^2}$.\vspace{2mm}\\
\noindent{\bf Claim(iv)} $B_2, B_3$, $B_4, B_5, B_6<\frac{8}{3}\frac{B_{16}}{\varepsilon}$

 As $B_1B_2B_3B_4<(3.2745)(\frac{4}{3}\frac{B_{16}}{(\varepsilon)^2})^3<41.765$, i.e. $0.25=\lambda_4^{(18)}<B_1B_2B_3B_4<\mu_4^{(18)}=41.765$, we find that $max\{\phi_{4,14}(\lambda_4^{(18)}),\phi_{4,14}(\mu_4^{(18)})\}=\phi_{4,14}(\mu_4^{(18)})$, which is $<22.955$. Therefore $B_{5}<max\{B_{6},B_{7},\cdots,B_{18}\}$. We have $B_6\leq2B_9\leq\frac{8}{3}\frac{B_{16}}{\varepsilon}$. Therefore using \eqref{2.8} we get $B_2, B_3$, $B_4$ and $B_5$ is $<\frac{8}{3}\frac{B_{16}}{\varepsilon}$.\vspace{2mm}\\
 \noindent{\bf Final Contradiction}

Now $2B_2+2B_4+2B_6+\cdots+2B_{18}<2(\frac{8/3}{\varepsilon}+\frac{8/3}{\varepsilon}+\frac{8/3}{\varepsilon}+\frac{(4/3)^2}{\varepsilon}+\frac{4/3}{\varepsilon}+\frac{1}{\varepsilon}+
\frac{3}{2}+1)B_{16}+2B_{18}<22.955$ for $B_{16}<0.385$ and $B_{18}<0.3$. This gives a contradiction to the weak inequality $(2,2,\cdots,2,2)_w$.~$\Box$
\subsection{$n=19$}
Here we have $\omega_{19}=24.691$, $B_1\leq\gamma_{19}<3.3974439$. Using \eqref{2.3}, we have $l_{19}=0.2009<B_{19}<3.0761736=m_{19}$ and using \eqref{2.1} we have
\begin{equation}\label{2.4.1}\begin{array}{l}
B_2\leq\gamma_{18}^{\frac{18}{19}}<3.0761736, ~~B_4\leq(2\gamma_{16}^{16})^{\frac{1}{19}}<2.635321,\\ B_5\leq(4\gamma_{15}^{15})^{\frac{1}{19}}<2.4940956,  B_6\leq(\frac{4\gamma_{14}^{14}}{\varepsilon})^{\frac{1}{19}}<2.3752798.
\end{array}\end{equation}

\noindent{\bf Claim(i)} $B_{19}<0.348$

Suppose $B_{19}\geq 0.348$.
The inequality $(18^*,1)$ gives $22.955(B_{19})^{\frac{-1}{18}}+B_{19}>24.691$. But this is not true for $0.348\leq B_{19}<3.0761736$. So we have  $B_{19}<0.348$. \vspace{2mm}\\
\noindent{\bf Claim(ii)} $B_{17}<0.403$

Suppose $B_{17}\geq0.403$. Then $B_{18}B_{19}\leq\frac{4}{3}B_{19}^2<\frac{4}{3}(0.348)^2<B_{17}^2$. So the inequality $(16^*,3)$ holds, i.e. $19.285({B_{17}B_{18}B_{19}})^{\frac{-1}{16}}+4B_{17}-\frac{B_{17}^3}{B_{18}B_{19}}>24.691$. But this is not true for $0.403\leq B_{17}\leq \frac{3}{2}B_{19}<0.522$ and $\frac{1}{2}B_{17}^2\leq B_{18}B_{19}<\frac{4}{3}(0.348)^2$. So we must have $B_{17}<0.403$.\vspace{2mm}\\
\noindent{\bf Claim(iii)} $B_{10}<max\{B_{11},B_{12},\cdots,B_{19}\}<1.2897$\vspace{1mm}\\

\noindent Using \eqref{2.2} and Lemma 9 we have\\
 $~~~~~~~\frac{\varepsilon^{10}}{2\times4^2B_{19}^{11}}<\lambda_8^{(19)}\leq B_1B_2\cdots B_8<\mu_8^{(19)}=(1.827901)^{11},$\\$~~~~~~~
  \frac{3\varepsilon^{8}}{4^3B_{19}^{10}}<\lambda_9^{(19)}\leq B_1B_2\cdots B_9<\mu_9^{(19)}=(1.994589)^{10}.$\\
We find that $max\{\phi_{9,10}(\lambda_9^{(19)}),\phi_{9,10}(\mu_9^{(19)})\}=\phi_{9,10}(\mu_9^{(19)})$, which is $<24.691$. Therefore using Lemma 10 we have $B_{10}<max\{B_{11},B_{12},\cdots,B_{19}\}$. From \eqref{2.0a},\eqref{2.0b}, we find that $B_{11}\leq\frac{3}{2}\frac{B_{17}}{\varepsilon}$ and each of $B_{12},\cdots, B_{17}$ is less than $\frac{3}{2}\frac{B_{17}}{\varepsilon}<1.2897$. Also $B_{19}<0.348$, $B_{18}<\frac{4}{3}B_{19}$. Therefore $max\{B_{11},B_{12},\cdots,B_{19}\}<1.2897$.\vspace{2mm}\\
\noindent{\bf Claim(iv)} $B_9<1.31015$

 We see that $\phi_{8,11}(\lambda_8^{(19)})<\omega_{19}$ but $\phi_{8,11}(\mu_8^{(19)})>\omega_{19}$. So we apply Lemma 11 with $\sigma_8^{(19)}=(1.8278)^{11}.$ Here $\phi_{8,11}(\sigma_8^{(19)})<\omega_{19}$.\vspace{1mm}\\
In Case(i), when $B_1B_2\cdots B_8<(1.8278)^{11}$, we have $B_{9}<max\{B_{10},B_{11},\cdots,B_{19}\}\\<1.2897$.\\
In Case(ii), when $B_1B_2\cdots B_8\geq(1.8278)^{11}$, then $B_9<\frac{\mu_9^{(19)}}{\sigma_8^{(19)}}<\frac{(1.994589)^{10}}{(1.8278)^{11}}\\<1.31015.$\\Hence $B_9<1.31015$.
\vspace{2mm}\\
\noindent{\bf Claim(v)} $B_{3}<max\{B_{4},B_{5},\cdots,B_{19}\}<2.635321$

Now suppose $B_3\geq max\{B_4,B_5,\cdots,B_{19}\}$, then the inequality $(2,17^*)$ gives $4B_1-\frac{2B_1^2}{B_2}+21.101(\frac{1}{B_1B_2})^{\frac{1}{17}}>24.691$. The left side is decreasing function of $B_1$, so replacing $B_1$ by $B_2$ we find that the inequality is not true for $\frac{3}{4}\leq B_2<3.0761736$. Hence we must have  $B_3<max\{B_4,B_5,\cdots,B_{19}\}$, which is $<2.635321$ as $B_{8}<\frac{4}{3}B_{9}$, $B_{7}<\frac{3}{2}B_{9}$, $B_{9}<1.31015$ and $B_{10}<max\{B_{11},B_{12},\cdots,B_{19}\}<1.2897$.\vspace{2mm}\\
\noindent{\bf Claim(vi)} $B_1>1.29$

For if $B_1\leq1.29$, then $B_i\leq B_1\leq1.29$ for each $i=2,3,\cdots,19$ and then the inequality $(1,1,\cdots,1)$ gives a contradiction.\vspace{2mm}\\
\noindent{\bf Claim(vii)} $B_{4}<max\{B_{5},B_{6},\cdots,B_{19}\}<2.4940956$, $B_3<2.4940956$

 Suppose $B_4\geq max\{B_5,B_6,\cdots,B_{19}\}$, then the inequality $(3,16^*)$ holds, i.e. $\eta(X)=4B_1-\frac{B_1^3}{X}+19.285(\frac{1}{B_1 X})^{\frac{1}{16}}>24.691$, where $X=B_2B_3$. $\eta'(X)=\frac{B_1^3}{X^2}(1-\frac{19.285}{16}(\frac{X^{15}}{B_1^{49}})^{\frac{1}{16}})>0$ for $1.29<B_1<3.3974439$ and $0.5\leq X\leq\alpha$, where $\alpha=min\{B_1^2,(3.0761736)(2.635321)\}$. Therefore $\eta(X)<\eta(\alpha)$, which is $<24.691$ for $1<B_1<3.3974439$, giving thereby a contradiction. Hence we have $B_4<max\{B_5,B_6,\cdots,B_{19}\}$, which is $<2.4940956$. So $B_4<2.4940956$ and hence $B_3<2.4940956$ using Claim(v).\vspace{2mm}\\
 \noindent{\bf Claim(viii)} $B_5<max\{B_6,B_7,\cdots,B_{19}\}<2.3752798$; $B_3$, $B_4<2.3752798$

  As $0.25\leq\lambda_4^{(19)}<B_1B_2B_3B_4<\mu_4^{(19)}=(3.3974439)(3.0761736)(2.4940956)^2$, we find that $max\{\phi_{4,15}(\lambda_4^{(19)}),\phi_{4,15}(\mu_4^{(19)})\}=\phi_{4,15}(\mu_4^{(19)})$, which is $<24.691$. Therefore using Lemma 10 we have $B_5<max\{B_6,B_7,\cdots,B_{19}\}$, which is $<2.3752798$. Therefore each of $B_3$, $B_4$ and $B_5$ is $<2.3752798$.\vspace{2mm}\\
  \noindent{\bf Claim(ix)} $B_6<max\{B_7,B_8,\cdots,B_{19}\}<1.9654$; $B_3$, $B_4, B_5<1.9654$

As $(0.25)(\varepsilon)=\lambda_5^{(19)}<B_1B_2B_3B_4B_5<\mu_5^{(19)}=(3.3974439)(3.0761736)(2.3752798)^3$, we find that $max\{\phi_{5,14}(\lambda_5^{(19)}),\phi_{5,14}(\mu_5^{(19)})\}=\phi_{5,14}(\mu_5^{(19)})$, which is $<24.691$. Therefore using Lemma 10 we have $B_6<max\{B_7,B_8,\cdots,B_{19}\}$, which is $\leq\frac{3}{2}B_9<\frac{3}{2}(1.31015)$. Therefore each of $B_6$, $B_5$, $B_4$ and $B_3$ is $ <\frac{3}{2}(1.31015)<1.9654$.\vspace{2mm}\\
\noindent{\bf Claim(x)} $B_7<max\{B_8,B_9,\cdots,B_{19}\} <1.747$; $B_6$, $B_5$, $B_4$ and $B_3$ is $ <1.747$

Now $\frac{3}{16}(\varepsilon^2)=\lambda_6^{(19)}<B_1B_2B_3B_4B_5B_6<\mu_6^{(19)}=(3.3974439)(3.0761736)(1.9654)^4$. We find that $max\{\phi_{6,13}(\lambda_6^{(19)}),\phi_{6,13}(\mu_6^{(19)})\}=\phi_{6,13}(\mu_6^{(19)})$, which is $<24.691$. Therefore using Lemma 10 we have $B_7<max\{B_8,B_9,\cdots,B_{19}\}$, which is $\leq\frac{4}{3}B_9<\frac{4}{3}(1.31015)<1.747$. Therefore each of $B_7$, $B_6$, $B_5$, $B_4$ and $B_3$ is $ <1.747$.\vspace{2mm}\\
\noindent{\bf Final Contradiction}

Now \vspace{1mm}\\$B_1+2B_3+2B_5+2B_7+\cdots+2B_{17}+2B_{19}<3.3974439+2(3\times1.747)+2(1.31015)+2(\frac{3/2}{\varepsilon}
+\frac{1}{\varepsilon}+\frac{3}{2}+1)B_{17}+2B_{19}<24.691$ for $B_{17}<0.403$ and $B_{19}<0.348$, giving thereby a contradiction to the weak inequality $(1,2,\cdots,2,2)_w$.~$\Box$
\subsection{$n=20$}
Here we have $\omega_{20}=26.629$, $B_1\leq\gamma_{20}<3.520062$. Using \eqref{2.3}, we have $l_{20}=0.1880<B_{20}<3.195912=m_{20}$ and using \eqref{2.1} we have $B_2\leq\gamma_{19}^{\frac{19}{20}}<3.195912$.\vspace{2mm}\\
\noindent{\bf Claim(i)} $B_{20}<0.294$

Suppose $B_{20}\geq0.294$. The inequality $(19^*,1)$ gives $24.691(B_{20})^{\frac{-1}{19}}+B_{20}>26.629$. But this is not true for $0.294\leq B_{20}<3.195912$. So we must have $B_{20}<0.294$.\vspace{2mm}\\
\noindent{\bf Claim(ii)} $B_i<max\{B_{i+1},B_{i+2},\cdots,B_{20}\}$, \mbox{for} $i=2,6,7,8,9,10$

 Using \eqref{2.2} and Lemmas 8, 9 we have\vspace{2mm}\\
  $~~~~~1=\lambda_1^{(20)}<B_1<\mu_1^{(20)}=3.520062,$\\
 $~~~~~\frac{1}{4}\varepsilon=\lambda_5^{(20)}<B_1B_2B_3B_4B_5<\mu_5^{(20)}=(1.4183645)^{15},$\\
 $~~~~~\frac{3}{8}\varepsilon^2=\lambda_6^{(20)}<B_1B_2B_3B_4B_5B_6<\mu_6^{(20)}=(1.530496)^{14},$\\
 $~~~~~\frac{\varepsilon^{15}}{4^3B_{20}^{13}}=\lambda_{7}^{(20)}<B_1\cdots B_7<\mu_7^{(20)}=(1.6555371)^{13},$\\
 $~~~~~\frac{\varepsilon^{12}}{4^3B_{20}^{12}}=\lambda_{8}^{(20)}<B_1\cdots B_8<\mu_8^{(20)}=(1.7955594)^{12},$\\
 $~~~~~\frac{\varepsilon^{10}}{2\times 4^2B_{20}^{11}}=\lambda_{9}^{(20)}<B_1\cdots B_9<\mu_9^{(20)}=(1.9530741)^{11}.$\vspace{2mm}\\
 We find that $max\{\phi_{s,n-s}(\lambda_{s}^{(20)}),\phi_{s,n-s}(\mu_{s}^{(20)}), {\mbox {for}}~s=1,5,6,7,8,9\}=\phi_{1,19}(\mu_{1}^{(20)})$, which is $<\omega_{20}$. Therefore using Lemma 10 we have
 $B_i<max\{B_{i+1},B_{i+2},\cdots,B_{20}\}$, \mbox{for} $i=2,6,7,8,9,10.$\vspace{2mm}\\
 \noindent{\bf Claim(iii)} $B_3<max\{B_{4},B_{5},\cdots,B_{20}\}$

  Suppose $B_3\geq max\{B_{4},B_{5},\cdots,B_{20}\}$, then the inequality $(2,18^*)$ holds, i.e. $4B_1-\frac{2B_1^2}{B_2}+22.955(\frac{1}{B_1B_2})^{\frac{1}{18}}>26.629$. The left side is decreasing function of $B_1$, so replacing $B_1$ by $B_2$ we find that the inequality is not true for $\frac{3}{4}\leq B_2<3.195912$. Hence $B_3<max\{B_{4},B_{5},\cdots,B_{20}\}$.

  Now using Claim(ii) and (iii), \begin{equation}\label{4.6.2}B_i<max\{B_{i+1},B_{i+2},\cdots,B_{20}\}, {\mbox {for}}~ i=2,3,6,7,8,9,10.\end{equation}  From  \eqref{2.0a},\eqref{2.0b} and Claim (i) we find that, $max\{B_{11},B_{12},\cdots,B_{20}\}<B_{11}\leq\frac{4}{3}\frac{B_{20}}{\varepsilon^2}$. So using \eqref{4.6.2} we get each of $B_{10}, B_9, B_8, B_7, B_{6}$ is $<\frac{4}{3}\frac{B_{20}}{\varepsilon^2}$ and so $B_5\leq\frac{4}{3}B_6<\frac{16}{9}\frac{B_{20}}{\varepsilon^2}$ and $B_4\leq\frac{3}{2}B_6<\frac{2B_{20}}{\varepsilon^2}$. Using \eqref{4.6.2} we have $B_3<\frac{2B_{20}}{\varepsilon^2}$ and hence $B_2<\frac{2B_{20}}{\varepsilon^2}$.\\
 We have now \vspace{2mm}\\
 $~~~~~\frac{1}{2}=\lambda_3^{(20)}<B_1B_2B_3<\mu_3^{(20)}=(3.520062)\left(\frac{2(0.294)}{\varepsilon^2}\right)^{2} {\mbox {and}}$\\
 $~~~~~\frac{1}{4}=\lambda_4^{(20)}<B_1B_2B_3B_4<\mu_4^{(20)}=(3.520062)\left(\frac{2(0.294)}{\varepsilon^2}\right)^{3}.$\vspace{2mm}\\
 We find that $max\{\phi_{s,n-s}(\lambda_{s}^{(20)}),\phi_{s,n-s}(\mu_{s}^{(20)}), {\mbox {for}}~s=3,4\}=\phi_{4,16}(\mu_{4}^{(20)})$, which is $<\omega_{20}$. Therefore using Lemma 10 we have
 $B_i<max\{B_{i+1},B_{i+2},\cdots,B_{20}\}$, for $i=4,5$. Using it together with \eqref{4.6.2} we get each of $B_2, B_3, \cdots, B_{10}$ is $<\frac{4}{3}\frac{B_{20}}{\varepsilon^2}$.\vspace{2mm}\\
 \noindent{\bf Final Contradiction}

 Now $2B_2+2B_4+\cdots+2B_{20}<2(5\times\frac{4/3}{\varepsilon^2}+\frac{1}{\varepsilon^2}+\frac{3/2}{\varepsilon}+\frac{1}{\varepsilon}+\frac{3}{2}+1)B_{20}<26.629$ for $B_{20}<0.294$, giving thereby a contradiction to the weak inequality $(2,2,\cdots,2,2)_w$.~$\Box$
\subsection{$n=21$}
Here we have $\omega_{21}=28.605$, $B_1\leq\gamma_{21}<3.6422432$. Using \eqref{2.3}, we have $l_{21}=0.1762<B_{21}<3.3153098=m_{21}$.\vspace{2mm}\\
\noindent{\bf Claim(i)} $B_{21}<0.2938$

Suppose $B_{21}\geq0.2938$. The inequality $(20^*,1)$ gives $26.629(B_{21})^{\frac{-1}{21}}+B_{21}>28.605$. But this is not true for $0.2938\leq B_{21}\leq 3.3153098$. So we must have $B_{21}<0.2938$.\vspace{2mm}\\
\noindent{\bf Claim(ii)} $B_{19}<0.4$

Suppose $B_{19}\geq0.4$. The inequality $(18^*,3)$ gives $22.955({B_{19}B_{20}B_{21}})^{\frac{-1}{18}}+4B_{19}-\frac{B_{19}^3}{B_{20}B_{21}}>28.605$. But this is not true for $0.4\leq B_{19}\leq \frac{3}{2}B_{21}<0.4407$ and $\frac{1}{2}B_{19}^2\leq B_{20}B_{21}\leq\frac{4}{3}B_{21}^2<\frac{4}{3}(0.2938)^2$. So we must have $B_{19}<0.4$.\vspace{2mm}\\
\noindent{\bf Claim(iii)} $B_2, B_3, \cdots, B_{10}<\frac{B_{19}}{\varepsilon^2}$

 Using \eqref{2.2} and Lemmas 8, 9 we have\vspace{2mm}\\
 $~~~~~1=\lambda_1^{(21)}<B_1<\mu_1^{(21)}=3.6422432,$\\
 $~~~~~\frac{3}{4}=\lambda_2^{(21)}<B_1B_2<\mu_2^{(21)}=(1.1398163)^{19},$\\
 $~~~~~\frac{1}{4}\varepsilon=\lambda_5^{(21)}<B_1\cdots B_5<\mu_5^{(21)}=(1.4053838)^{16},$\\
 $~~~~~\frac{3}{16}\varepsilon^2=\lambda_6^{(21)}<B_1\cdots B_6<\mu_6^{(21)}=(1.5130617)^{15},$\\
 $~~~~~\frac{3\varepsilon^{18}}{4^4B_{21}^{14}}=\lambda_{7}^{(21)}<B_1\cdots B_7<\mu_7^{(21)}=(1.6326792)^{14},$\\
 $~~~~~\frac{\varepsilon^{15}}{4^3B_{21}^{13}}=\lambda_{8}^{(21)}<B_1\cdots B_8<\mu_8^{(21)}=(1.7660691)^{13},$\\
 $~~~~~\frac{\varepsilon^{12}}{4^3B_{21}^{12}}=\lambda_{9}^{(21)}<B_1\cdots B_9<\mu_9^{(21)}=(1.915440)^{12}.$\vspace{2mm}\\
We find that $max\{\phi_{s,n-s}(\lambda_{s}^{(21)}),\phi_{s,n-s}(\mu_{s}^{(21)}), {\mbox {for}}~s=1,2,5,6,7,8,9\}=\phi_{1,20}(\mu_{1}^{(21)})$, which is $<\omega_{21}$. Therefore using Lemma 10 we have
 \begin{equation}\label{4.6.3}B_i<max\{B_{i+1},B_{i+2},\cdots,B_{21}\},~ {\mbox {for}}~ i=2,3,6,7,8,9,10.\end{equation}
From  \eqref{2.0a},\eqref{2.0b} and Claims (i), (ii) we find that, $max\{B_{11},B_{12},\cdots,B_{21}\}<B_{11}\leq\frac{B_{19}}{\varepsilon^2}$. So using \eqref{4.6.3} we get each of $B_{10}, B_9, B_8, B_7, B_{6}$ is $<\frac{B_{19}}{\varepsilon^2}$ and so $B_5\leq\frac{4}{3}B_6<\frac{4}{3}\frac{B_{19}}{\varepsilon^2}$ and $B_4\leq\frac{3}{2}B_6<\frac{3}{2}\frac{B_{19}}{\varepsilon^2}$. Using \eqref{4.6.3} we have $B_3\leq\frac{3}{2}\frac{B_{19}}{\varepsilon^2}$ and hence $B_2\leq\frac{3}{2}\frac{B_{19}}{\varepsilon^2}$.\\
 We have now \vspace{2mm}\\$~~~\frac{1}{2}=\lambda_3^{(21)}<B_1B_2B_3<\mu_3^{(21)}=(3.6422432)\left(\frac{3}{2}\frac{0.4}{\varepsilon^2}\right)^{2} {\mbox {and}}$\\
 $~~~\frac{1}{4}=\lambda_4^{(21)}<B_1B_2B_3B_4<\mu_4^{(21)}=(3.6422432)\left(\frac{3}{2}\frac{0.4}{\varepsilon^2}\right)^{3}.$\vspace{2mm}\\
Now we find that $max\{\phi_{s,n-s}(\lambda_{s}^{(21)}),\phi_{s,n-s}(\mu_{s}^{(21)}), {\mbox {for}}~s=3,4\}=\phi_{3,18}(\mu_{3}^{(21)})$, which is $<\omega_{21}$. Therefore using Lemma 10 we have
 $B_i<max\{B_{i+1},B_{i+2},\cdots,B_{21}\}$, for $i=4,5$. Using it together with \eqref{4.6.3} we get each of $B_2, B_3, \cdots, B_{10}$ is $<max\{B_{11},B_{12},\cdots,B_{21}\}< B_{11}\leq\frac{B_{19}}{\varepsilon^2}$.\vspace{2mm}\\
 \noindent{\bf Final Contradiction}

 \noindent Now $2B_2+B_{3}+2B_5+2B_7\cdots+2B_{19}+2B_{21}<\{3(\frac{1}{\varepsilon^2})+2(4\times\frac{1}{\varepsilon^2}+\frac{3/2}{\varepsilon}+\frac{1}{\varepsilon}+
 \frac{3}{2}+1)\}B_{19}+2B_{21}<28.605$ for $B_{21}<0.2938$ and $B_{19}<0.4$, giving thereby a contradiction to the weak inequality $(2,1,2,2,\cdots,2)_w$.~$\Box$
\subsection{$n=22$}
Here we have $\omega_{22}=30.62$, $B_1\leq\gamma_{22}<3.7640371$. Using \eqref{2.3}, we have $l_{22}=0.1655<B_{22}<3.4344103=m_{22}$. \vspace{2mm}\\
\noindent{\bf Claim(i)} $B_{22}<0.295$

Suppose $B_{22}\geq 0.295$. The inequality $(22^*,1)$ gives $28.605(B_{22})^{\frac{-1}{21}}+B_{22}>30.62$. But this is not true for $0.295\leq B_{22}\leq 3.4344103$. \vspace{2mm}\\
\noindent{\bf Claim(ii)} $B_{20}<0.378$

Suppose $B_{20}\geq0.378$. Then $B_{21}B_{22}\leq\frac{4}{3}B_{22}^2<\frac{4}{3}(0.295)^2<B_{20}^2$. Therefore the inequality $(19^*,3)$ holds, i.e. $24.691({B_{20}B_{21}B_{22}})^{\frac{-1}{19}}+4B_{20}-\frac{B_{20}^3}{B_{21}B_{22}}>30.62$. But this is not true for $0.378\leq B_{20}\leq \frac{3}{2}B_{22}<0.4425$ and $\frac{1}{2}B_{20}^2\leq B_{21}B_{22}<\frac{4}{3}(0.295)^2$. So we must have $B_{20}<0.378$.\vspace{2mm}\\
\noindent{\bf Claim(iii)} $B_2, B_3, \cdots, B_{11}$ is $<\frac{B_{20}}{\varepsilon^2}$

 Using \eqref{2.2} and Lemmas 8, 9 we have\vspace{2mm}\\
 $~~~~~1=\lambda_1^{(22)}<B_1<\mu_1^{(22)}=3.7640371,$\\
 $~~~~~\frac{3}{4}=\lambda_2^{(22)}<B_1B_2<\mu_2^{(22)}=(1.1362692)^{20},$\\
 $~~~~~\frac{1}{2}=\lambda_3^{(22)}<B_1B_2B_3<\mu_3^{(22)}=(1.21407992)^{19},$\\
 $~~~~~\frac{\varepsilon}{4}=\lambda_5^{(22)}<B_1\cdots B_5<\mu_5^{(22)}=(1.39334)^{17},$\\
 $~~~~~\frac{3\varepsilon^2}{16}=\lambda_6^{(22)}<B_1\cdots B_6<\mu_6^{(22)}=(1.496951)^{16},$\\
 $~~~~~\frac{\varepsilon^{21}}{2\times 4^3B_{22}^{15}}=\lambda_{7}^{(22)}<B_1\cdots B_7<\mu_7^{(22)}=(1.611638)^{15},$\\
 $~~~~~\frac{3\varepsilon^{18}}{4^4B_{22}^{14}}=\lambda_{8}^{(22)}<B_1\cdots B_8<\mu_8^{(22)}=(1.739055)^{14},$\\
 $~~~~~\frac{\varepsilon^{15}}{4^3B_{22}^{13}}=\lambda_{9}^{(22)}<B_1\cdots B_9<\mu_9^{(22)}=(1.8811357)^{13},$\\
 $~~~~~\frac{\varepsilon^{12}}{4^3B_{22}^{12}}=\lambda_{10}^{(22)}<B_1\cdots B_{10}<\mu_{10}^{(22)}=(2.0402387)^{12}.$\vspace{2mm}\\
We find that $max\{\phi_{s,n-s}(\lambda_{s}^{(22)}),\phi_{s,n-s}(\mu_{s}^{(22)}), {\mbox {for}}~s=1,2,3,5,6,7,8,9,10\}=\phi_{1,21}(\mu_{1}^{(22)})$, which is $<\omega_{22}$. Therefore using Lemma 10 we have
 \begin{equation}\label{4.6.4}B_i<max\{B_{i+1},B_{i+2},\cdots,B_{22}\},~ {\mbox {for}}~ i=2,3,4,6,7,8,9,10,11.\end{equation}
From  \eqref{2.0a},\eqref{2.0b} and Claims(i), (ii) we find that, $max\{B_{12},B_{13},\cdots,B_{22}\}<B_{12}\leq\frac{B_{20}}{\varepsilon^2}$. So using \eqref{4.6.4} we get each of $B_{11}, B_{10}, B_9, B_8, B_7, B_{6}$ is $<\frac{B_{20}}{\varepsilon^2}$ and so $B_5\leq\frac{4}{3}B_6<\frac{4}{3}\frac{B_{20}}{\varepsilon^2}$. Again using \eqref{4.6.4} we get each of $B_4$, $B_3$ and $B_2$ is $<\frac{4}{3}\frac{B_{20}}{\varepsilon^2}$. \\
 We have now
 \vspace{2mm}\\
 $~~~\frac{1}{4}=\lambda_4^{(22)}<B_1B_2B_3B_4<\mu_4^{(22)}=(3.7640371)\left(\frac{4}{3}\frac{0.378}{\varepsilon^2}\right)^{3}.$\vspace{2mm}\\
 We find that $max\{\phi_{4,18}(\lambda_{4}^{(22)}),\phi_{4,18}(\mu_{4}^{(22)})\}=\phi_{4,18}(\mu_{4}^{(22)})$, which is $<\omega_{22}$. Therefore using Lemma 10 we have
 $B_5<max\{B_{6},B_{7},\cdots,B_{22}\}$. Using it together with \eqref{4.6.4} we get each of $B_2, B_3, \cdots, B_{11}$ is $<\frac{B_{20}}{\varepsilon^2}$.\vspace{2mm}\\
 \noindent{\bf Final Contradiction}

\noindent Now $2B_2+2B_{4}+\cdots+2B_{20}+2B_{22}<2\{6\times\frac{1}{\varepsilon^2}+\frac{3/2}{\varepsilon}+\frac{1}{\varepsilon}+
 \frac{3}{2}+1\}B_{20}+2B_{22}<30.62$ for $B_{22}<0.295$ and $B_{20}<0.378$, giving thereby a contradiction to the weak inequality $(2,2,\cdots,2,2)_w$.~$\Box$
\subsection{$n=23$}
Here we have $\omega_{23}=32.68$, $B_1\leq\gamma_{23}<3.8854763$. Using \eqref{2.3}, we have $l_{23}=0.1556<B_{23}<3.5532476=m_{22}$.\vspace{2mm}\\
\noindent{\bf Claim(i)} $B_{23}<0.293$

Suppose $B_{23}\geq0.293$. The inequality $(22^*,1)$ gives $30.62(B_{23})^{\frac{-1}{22}}+B_{23}>32.68$. But this is not true for $0.293\leq B_{23}<3.5532476$. So we must have $B_{23}<0.293$.\vspace{2mm}\\
\noindent{\bf Claim(ii)} $B_{21}<0.376$

Suppose $B_{21}\geq0.376$. The inequality $(20^*,3)$ gives $26.629({B_{21}B_{22}B_{23}})^{\frac{-1}{20}}+4B_{21}-\frac{B_{21}^3}{B_{22}B_{23}}>32.68$. But this is not true for $0.376\leq B_{21}\leq \frac{3}{2}B_{23}<0.4395$ and $\frac{1}{2}B_{21}^2\leq B_{22}B_{23}\leq\frac{4}{3}B_{22}^2<\frac{4}{3}(0.293)^2$. So we must have $B_{21}<0.376$.\vspace{2mm}\\
\noindent{\bf Claim(iii)} $B_2,\cdots,B_{11}$ $<max\{B_{12},\frac{B_{21}}{\varepsilon^2}\}$

 Using \eqref{2.2} and Lemmas 8, 9 we have\vspace{2mm}\\
 $~~~~~1=\lambda_1^{(23)}<B_1<\mu_1^{(23)}=3.8854763,$\\
 $~~~~~\frac{3}{4}=\lambda_2^{(23)}<B_1B_2<\mu_2^{(23)}=(1.1329363)^{21},$\\
 $~~~~~\frac{1}{2}=\lambda_3^{(23)}<B_1B_2B_3<\mu_3^{(23)}=(1.2085769)^{20},$\\
  $~~~~~\frac{1}{4}=\lambda_4^{(23)}<B_1B_2B_3B_4<\mu_4^{(23)}=(1.2913392)^{19},$\\
 $~~~~~\frac{\varepsilon}{4}=\lambda_5^{(23)}<B_1\cdots B_5<\mu_5^{(23)}=(1.3821298)^{18},$\\
 $~~~~~\frac{3\varepsilon^2}{16}=\lambda_6^{(23)}<B_1\cdots B_6<\mu_6^{(23)}=(1.4820067)^{17},$\\
 $~~~~~\frac{\varepsilon^3}{8}=\lambda_{7}^{(23)}<B_1\cdots B_7<\mu_7^{(23)}=(1.5922103)^{16},$\\
 $~~~~~\frac{\varepsilon^4}{16}=\lambda_{8}^{(23)}<B_1\cdots B_8<\mu_8^{(23)}=(1.7142025)^{15},$\\
 $~~~~~\frac{3\varepsilon^{18}}{4^4B_{23}^{14}}=\lambda_{9}^{(23)}<B_1\cdots B_9<\mu_9^{(23)}=(1.8497215)^{14},$\\
 $~~~~~\frac{\varepsilon^{15}}{4^3B_{23}^{13}}=\lambda_{10}^{(23)}<B_1\cdots B_{10}<\mu_{10}^{(23)}=(2.000844)^{13},$\\
 $~~~~~\frac{\varepsilon^{12}}{4^3B_{23}^{12}}=\lambda_{11}^{(23)}<B_1\cdots B_{11}<\mu_{11}^{(23)}=(2.1700714)^{12},$\\
 $~~~~~\frac{\varepsilon^{10}}{2\times4^2B_{23}^{11}}=\lambda_{12}^{(23)}<B_1\cdots B_{12}<\mu_{12}^{(23)}=(2.3604401)^{11}.$\vspace{2mm}\\
We find that $max\{\phi_{s,n-s}(\lambda_{s}^{(23)}),\phi_{s,n-s}(\mu_{s}^{(23)}), {\mbox {for}}~s=1,2,3,4,5,6,7,8,9,10\}=\phi_{1,22}(\mu_{1}^{(23)})$, which is $<\omega_{23}$. Therefore using Lemma 10 we have
 \begin{equation}\label{4.6.5}B_i<max\{B_{i+1},B_{i+2},\cdots,B_{23}\},~ {\mbox {for}}~ i=2,3,4,5,6,7,8,9,10,11.\end{equation}
From \eqref{2.0a},\eqref{2.0b} and Claims (i),(ii) we find that, $max\{B_{13},B_{14},\cdots,B_{23}\}<B_{13}\leq\frac{B_{21}}{\varepsilon^2}<1.7114$. So using \eqref{4.6.5} we get each of $B_2,\cdots,B_{11}$ is $<max\{B_{12},\frac{B_{21}}{\varepsilon^2}\}$. \vspace{2mm}\\
\noindent{\bf Claim(iv)} $B_{12}<\frac{B_{21}}{\varepsilon^2}$

We find that $\phi_{11,12}(\lambda_{11}^{(23)})<\omega_{23}$, but $\phi_{11,12}(\mu_{11}^{(23)})>\omega_{23}$, so we apply Lemma 11 with $\sigma_{11}^{(23)}=(2.135)^{12}$. Here $\phi_{11,12}(\sigma_{11}^{(23)})<\omega_{23}$.\\
In Case(i), when $B_1B_2\cdots B_{11}<(2.135)^{12}$, then we have $B_{12}<max\{B_{13},B_{14},\cdots,B_{23}\}$, which is $\leq\frac{B_{21}}{\varepsilon^2}<1.7114$.\\
In Case(ii), when $B_1B_2\cdots B_{11}\geq(2.135)^{12}$, then we have $B_{12}<\frac{\mu_{12}^{(23)}}{\sigma_{11}^{(23)}}<\frac{(2.3604)^{11}}{(2.135)^{12}}<1.43$. \\
So we have $B_{12}<\frac{B_{21}}{\varepsilon^2}$.\vspace{1mm}\\ Using Claims(iii), (iv) we get, each of $B_2, B_3,\cdots,B_{13}$ is $<\frac{B_{21}}{\varepsilon^2}$.\vspace{2mm}\\
 \noindent{\bf Final Contradiction}

 Now $2B_2+B_3+2B_5+2B_7+\cdots+2B_{23}<3(\frac{B_{21}}{\varepsilon^2})+2\{5(\frac{1}{\varepsilon^2})+\frac{3/2}{\varepsilon}+\frac{1}{\varepsilon}+\frac{3}{2}+1\}
B_{21}+2B_{23}<32.68$ for $B_{21}<0.376$ and $B_{23}<0.293$, giving thereby a contradiction to the weak inequality $(2,1,2,2,\cdots,2,2)_w$.~$\Box$
\subsection{$n=24$}
Here we have $\omega_{24}=34.78$, $B_1\leq\gamma_{24}<4.0065998$. Using \eqref{2.3}, we have $l_{24}=0.1464<B_{24}<3.6718429=m_{24}$.\vspace{2mm}\\
\noindent{\bf Claim(i)} $B_{24}<0.29$

Suppose $B_{24}\geq0.29$. The inequality $(23^*,1)$ gives $32.68(B_{24})^{\frac{-1}{23}}+B_{24}>34.78$. But this is not true for $0.29\leq B_{24}<3.6718429$. So we must have $B_{24}<0.29$.\vspace{2mm}\\
\noindent{\bf Claim(ii)} $B_{22}<0.374$

Suppose $B_{22}\geq0.374$. The inequality $(21^*,3)$ gives $28.605({B_{22}B_{23}B_{24}})^{\frac{-1}{21}}+4B_{22}-\frac{B_{22}^3}{B_{23}B_{24}}>34.78$. But this is not true for $0.374\leq B_{22}\leq \frac{3}{2}B_{24}<0.435$ and $\frac{1}{2}B_{22}^2\leq B_{23}B_{24}\leq\frac{4}{3}B_{24}^2<\frac{4}{3}(0.29)^2$. So we must have $B_{22}<0.374$.\vspace{2mm}\\
\noindent{\bf Claim(iii)} $B_i<max\{B_{12},B_{13},1.703\},~ {\mbox {for}}~ i=2,3,4,5,6,7,8,9,10,11$

Using \eqref{2.2} and Lemmas 8, 9 we have\vspace{2mm}\\
 $~~~~~1=\lambda_1^{(24)}<B_1<\mu_1^{(24)}=4.0065998,$\\
 $~~~~~\frac{3}{4}=\lambda_2^{(24)}<B_1B_2<\mu_2^{(24)}=(1.1298027)^{22},$\\
 $~~~~~\frac{1}{2}=\lambda_3^{(24)}<B_1B_2B_3<\mu_3^{(24)}=(1.2034139)^{21},$\\
  $~~~~~\frac{1}{4}=\lambda_4^{(24)}<B_1B_2B_3B_4<\mu_4^{(24)}=(1.2837599)^{20},$\\
 $~~~~~\frac{\varepsilon}{4}=\lambda_5^{(24)}<B_1\cdots B_5<\mu_5^{(24)}=(1.3716707)^{19},$\\
 $~~~~~\frac{3\varepsilon^2}{16}=\lambda_6^{(24)}<B_1\cdots B_6<\mu_6^{(24)}=(1.4681092)^{18},$\\
 $~~~~~\frac{\varepsilon^3}{8}=\lambda_{7}^{(24)}<B_1\cdots B_7<\mu_7^{(24)}=(1.5741992)^{17},$\\
 $~~~~~\frac{\varepsilon^4}{16}=\lambda_{8}^{(24)}<B_1\cdots B_8<\mu_8^{(24)}=(1.6912584)^{16},$\\
 $~~~~~\frac{\varepsilon^{21}}{2\times4^3B_{24}^{15}}=\lambda_{9}^{(24)}<B_1\cdots B_9<\mu_9^{(24)}=(1.8208394)^{15},$\\
 $~~~~~\frac{3\varepsilon^{18}}{4^4B_{24}^{14}}=\lambda_{10}^{(24)}<B_1\cdots B_{10}<\mu_{10}^{(24)}=(1.9647888)^{14},$\\
 $~~~~~\frac{\varepsilon^{15}}{4^3B_{24}^{13}}=\lambda_{11}^{(24)}<B_1\cdots B_{11}<\mu_{11}^{(24)}=(2.1253121)^{13},$\\
 $~~~~~\frac{\varepsilon^{12}}{4^3B_{24}^{12}}=\lambda_{12}^{(24)}<B_1\cdots B_{12}<\mu_{12}^{(24)}=(2.3050671)^{12},$\\
  $~~~~~\frac{\varepsilon^{10}}{2\times4^2B_{24}^{11}}=\lambda_{13}^{(24)}<B_1\cdots B_{13}<\mu_{13}^{(24)}=(2.5072781)^{11}.$\vspace{2mm}\\
We find that $max\{\phi_{s,n-s}(\lambda_{s}^{(24)}),\phi_{s,n-s}(\mu_{s}^{(24)}), {\mbox {for}}~s=1,2,3,4,5,6,7,8,9,10\}=\phi_{1,23}(\mu_{1}^{(24)})$, which is $<\omega_{24}$. Therefore using Lemma 10 we have
 \begin{equation}\label{4.6.6}B_i<max\{B_{i+1},B_{i+2},\cdots,B_{24}\},~ {\mbox {for}}~ i=2,3,4,5,6,7,8,9,10,11.\end{equation}
From \eqref{2.0a},\eqref{2.0b} and Claims (i),(ii), we find that  $max\{B_{14},B_{15},\cdots,B_{24}\}<B_{14}\leq\frac{B_{22}}{\varepsilon^2}<1.703$. So using \eqref{4.6.6} we get each of $B_2, \cdots,B_{11}$ is $<max\{B_{12},B_{13},\frac{B_{22}}{\varepsilon^2}\}<max\{B_{12},B_{13},1.703\}$.\vspace{2mm} \\
\noindent{\bf Claim(iv)} $B_{12},~B_{13}<1.72$

We find that for $s=12,11$, $\phi_{s,n-s}(\lambda_{s}^{(24)})<\omega_{24}$, but $\phi_{s,n-s}(\mu_{s}^{(24)})>\omega_{24}$, so we apply Lemma 11 respectively with $\sigma_{12}^{(24)}=(2.22)^{12}$ and $\sigma_{11}^{(24)}=(2.09)^{13}$. Here $\phi_{12,12}(\sigma_{12}^{(24)})<\omega_{24}$ and $\phi_{11,13}(\sigma_{11}^{(24)})<\omega_{24}$.\\
First consider Lemma 11 for $s=12$ and with $\sigma_{12}^{(24)}=(2.22)^{12}$.\\
In Case(i), when $B_1B_2\cdots B_{12}<(2.22)^{12}$, then we have $B_{13}<max\{B_{14},B_{15},\cdots,B_{24}\}$, which is $\leq\frac{B_{22}}{\varepsilon^2}<1.703$.\\
In Case(ii), when $B_1B_2\cdots B_{12}\geq(2.22)^{12}$, then we have $B_{13}<\frac{\mu_{13}^{(24)}}{\sigma_{12}^{(24)}}<\frac{(2.50718)^{11}}{(2.22)^{12}}<1.72$. \\
So we have $B_{13}<1.72$.\\
Now consider Lemma 11 for $s=11$ and with $\sigma_{11}^{(24)}=(2.09)^{13}$.\\
In Case(i), when $B_1B_2\cdots B_{11}<(2.09)^{13}$, then we have $B_{12}<max\{B_{13},B_{14},\cdots,B_{24}\}$, which is $<B_{13}<1.72$.\\
In Case(ii), when $B_1B_2\cdots B_{11}\geq(2.09)^{13}$, then we have $B_{12}<\frac{\mu_{12}^{(24)}}{\sigma_{11}^{(24)}}<\frac{(2.305)^{12}}{(2.09)^{13}}<1.56$.\vspace{2mm} \\
Using Claim(iii) and (iv) we get $B_2, B_3,\cdots,B_{13}$ is $<1.72$.\vspace{2mm} \\
\noindent{\bf Final Contradiction}

Now $2B_2+2B_4+\cdots+2B_{24}<2(6\times1.72)+2\{\frac{1}{\varepsilon^2}+\frac{3/2}{\varepsilon}+\frac{1}{\varepsilon}+\frac{3}{2}+1\}B_{22}
++2B_{24}<34.78$ for $B_{22}<0.374$ and $B_{24}<0.29$, giving thereby a contradiction to the weak inequality $(2,2,\cdots,2,2)_w$.~$\Box$
\subsection{$n=25$}
Here we have $\omega_{25}=37.05$, $B_1\leq\gamma_{25}<4.1274438$. Using \eqref{2.3}, we have $l_{25}=0.1380<B_{25}<3.7902246=m_{25}$.\vspace{2mm}\\
\noindent{\bf Claim(i)} $B_{25}<0.26$

Suppose $B_{25}\geq0.26$. The inequality $(24^*,1)$ gives $34.78(B_{25})^{\frac{-1}{24}}+B_{25}>37.05$. But this is not true for $0.26\leq B_{25}\leq 3.7902246$. So we must have $B_{25}<0.26$.\vspace{2mm}\\
\noindent{\bf Claim(ii)} $B_{24}<0.311$

Suppose $B_{24}\geq0.311$. The inequality $(23^*,2)$ gives $32.68({B_{24}B_{25}})^{\frac{-1}{23}}+4B_{24}-\frac{2B_{24}^2}{B_{25}}>37.05$. But this is not true for $0.311\leq B_{24}\leq \frac{4}{3}B_{25}<0.347$ and $\frac{3}{4}(B_{24})\leq B_{25}<0.26$. So we must have $B_{24}<0.311$.\vspace{2mm}\\
\noindent{\bf Claim(iii)} $B_{23}<0.3595$

Suppose $B_{23}\geq0.3595$. The inequality $(22^*,3)$ gives $30.62({B_{23}B_{24}B_{25}})^{\frac{-1}{22}}+4B_{23}-\frac{B_{23}^3}{B_{24}B_{25}}>37.05$. But this is not true for $0.3595\leq B_{23}\leq \frac{3}{2}B_{25}<0.39$ and $\frac{1}{2}B_{23}^2\leq B_{24}B_{25}<(0.311)(0.26)$. So we must have $B_{23}<0.3595$.\vspace{2mm}\\
\noindent{\bf Claim(iv)} $B_i<max\{B_{12},B_{13},2.124\},~ {\mbox {for}}~ i=2,3,4,5,6,7,8,9,10,11$

 Using \eqref{2.2} and Lemmas 8, 9 we have\vspace{2mm}\\
 $~~~~~1=\lambda_1^{(25)}<B_1<\mu_1^{(25)}=4.1274438,$\\
 $~~~~~\frac{3}{4}=\lambda_2^{(25)}<B_1B_2<\mu_2^{(25)}=(1.1268424)^{23},$\\
 $~~~~~\frac{1}{2}=\lambda_3^{(25)}<B_1B_2B_3<\mu_3^{(25)}=(1.1985502)^{22},$\\
  $~~~~~\frac{1}{4}=\lambda_4^{(25)}<B_1B_2B_3B_4<\mu_4^{(25)}=(1.2766406)^{21},$\\
 $~~~~~\frac{\varepsilon}{4}=\lambda_5^{(25)}<B_1\cdots B_5<\mu_5^{(25)}=(1.3618756)^{20},$\\
 $~~~~~\frac{3\varepsilon^2}{16}=\lambda_6^{(25)}<B_1\cdots B_6<\mu_6^{(25)}=(1.4551357)^{19},$\\
 $~~~~~\frac{\varepsilon^3}{8}=\lambda_{7}^{(25)}<B_1\cdots B_7<\mu_7^{(25)}=(1.5574424)^{18},$\\
  $~~~~~\frac{\varepsilon^{28}}{4^4B_{25}^{17}}=\lambda_{8}^{(25)}<B_1\cdots B_8<\mu_8^{(25)}=(1.669988)^{17},$\\
 $~~~~~\frac{\varepsilon^{24}}{4^4B_{25}^{16}}=\lambda_{9}^{(25)}<B_1\cdots B_9<\mu_9^{(25)}=(1.794170)^{16},$\\
 $~~~~~\frac{\varepsilon^{21}}{2\times4^3B_{25}^{15}}=\lambda_{10}^{(25)}<B_1\cdots B_{10}<\mu_{10}^{(25)}=(1.9316359)^{15},$\\
 $~~~~~\frac{3\varepsilon^{18}}{4^4B_{25}^{14}}=\lambda_{11}^{(25)}<B_1\cdots B_{11}<\mu_{11}^{(25)}=(2.0843445)^{14},$\\
 $~~~~~\frac{\varepsilon^{15}}{4^3B_{25}^{13}}=\lambda_{12}^{(25)}<B_1\cdots B_{12}<\mu_{12}^{(25)}=(2.2546355)^{13},$\\
  $~~~~~\frac{\varepsilon^{12}}{4^3B_{25}^{12}}=\lambda_{13}^{(25)}<B_1\cdots B_{13}<\mu_{13}^{(25)}=(2.4453284)^{12}.$\vspace{2mm}\\
We find that $max\{\phi_{s,n-s}(\lambda_{s}^{(25)}),\phi_{s,n-s}(\mu_{s}^{(25)}), {\mbox {for}}~s=1,2,3,4,5,6,7,8,9,10\}=\phi_{1,24}(\mu_{1}^{(25)})$, which is $<\omega_{25}$. Therefore using Lemma 10 we have
 \begin{equation}\label{4.6.7}B_i<max\{B_{i+1},B_{i+2},\cdots,B_{25}\},~ {\mbox {for}}~ i=2,3,4,5,6,7,8,9,10,11.\end{equation}
From \eqref{2.0a},\eqref{2.0b} and Claims(i),(ii) and (iii) we find that, $max\{B_{14},B_{15},\cdots,B_{25}\}<B_{14}\leq\frac{3}{2}\frac{B_{24}}{\varepsilon^2}<2.124$. So using \eqref{4.6.7} we get each of $B_2, \cdots,B_{11}$ is $<max\{B_{12},B_{13},\frac{3}{2}\frac{B_{24}}{\varepsilon^2}\}<max\{B_{12},B_{13},2.124\}$.\vspace{2mm} \\
\noindent{\bf Claim(v)} $B_{12}, B_{13}<2.124$

Now we find that for $s=12,11$, $\phi_{s,n-s}(\lambda_{s}^{(25)})<\omega_{25}$, but $\phi_{s,n-s}(\mu_{s}^{(25)})>\omega_{25}$, so we apply Lemma 11 respectively with $\sigma_{12}^{(25)}=(2.193)^{13}$ and $\sigma_{11}^{(25)}=(2.07)^{14}$. Here $\phi_{12,13}(\sigma_{12}^{(25)})<\omega_{25}$ and $\phi_{11,14}(\sigma_{11}^{(25)})<\omega_{25}$.\\
First consider Lemma 11 for $s=12$ and with $\sigma_{12}^{(25)}=(2.193)^{13}$.\\
In Case(i), when $B_1B_2\cdots B_{12}<(2.193)^{13}$, then we have $B_{13}<max\{B_{14},B_{15},\cdots,B_{25}\}$, which is $\leq\frac{3}{2}\frac{B_{24}}{\varepsilon^2}<2.124$.\\
In Case(ii), when $B_1B_2\cdots B_{12}\geq(2.193)^{13}$, then we have $B_{13}<\frac{\mu_{13}^{(25)}}{\sigma_{12}^{(25)}}<\frac{(2.4453)^{12}}{(2.193)^{13}}<1.686$. \\
So we have $B_{13}\leq\frac{3}{2}\frac{B_{24}}{\varepsilon^2}<2.124$.\\
Now consider Lemma 11 for $s=11$ and with $\sigma_{11}^{(25)}=(2.07)^{14}$.\\
In Case(i), when $B_1B_2\cdots B_{11}<(2.07)^{14}$, then we have $B_{12}<max\{B_{13},B_{14},\cdots,B_{25}\}$, which is $<B_{13}\leq\frac{3}{2}\frac{B_{24}}{\varepsilon^2}<2.124$.\\
In Case(ii), when $B_1B_2\cdots B_{11}\geq(2.07)^{14}$, then we have $B_{12}<\frac{\mu_{12}^{(25)}}{\sigma_{11}^{(25)}}<\frac{(2.255)^{13}}{(2.07)^{14}}<1.471$. \\
So we have $B_{12}, B_{13}<2.124$.\vspace{2mm}\\
\noindent Using Claim(iii) and (iv) we get each of $B_2, B_3,\cdots,B_{14}$ is $<B_{13}<2.124$.\vspace{2mm}\\
 \noindent{\bf Final Contradiction}

 Now $2B_2+B_3+2B_5+\cdots+2B_{25}<3(2.124)+2(5\times2.124)+2\{(\frac{1}{\varepsilon^2})+\frac{3/2}{\varepsilon}+\frac{1}{\varepsilon}+\frac{3}{2}+1\}B_{23}
+2B_{25}<37.05$ for $B_{25}<0.26$ and $B_{23}<0.3595$, giving thereby a contradiction to the weak inequality $(2,1,2,\cdots,2,2)_w$.~$\Box$
\subsection{$n=26$}
Here we have $\omega_{26}=39.24$, $B_1\leq\gamma_{26}<4.2480446$. Using \eqref{2.3}, we have $l_{26}=0.1303<B_{26}<3.9084192=m_{26}$.\vspace{2mm}\\
\noindent{\bf Claim(i)} $B_{26}<0.29$

Suppose $B_{26}\geq 0.29$.
The inequality $(25^*,1)$ gives $37.05(B_{26})^{\frac{-1}{25}}+B_{26}>39.24$. But this is not true for $0.29\leq B_{26}<3.9084192$. \vspace{2mm}\\
\noindent{\bf Claim(ii)} $B_{25}<0.31$

Suppose $B_{25}\geq0.31$. Then $2B_{25}>B_{26}$. Therefore the inequality $(24^*,2)$ holds, i.e. $34.78({B_{25}B_{26}})^{\frac{-1}{24}}+4B_{25}-\frac{2B_{25}^2}{B_{26}}>39.24$. But this is not true for $0.31\leq B_{25}\leq \frac{4}{3}B_{26}<0.387$ and $\frac{3}{4}(B_{25})\leq B_{26}<0.29$. So we must have $B_{25}<0.31$.\vspace{2mm}\\
\noindent{\bf Claim(iii)} $B_{24}<0.358$

Suppose $B_{24}\geq0.358$. Then $B_{24}^2>B_{25}B_{26}$. Therefore the inequality $(23^*,3)$ holds, i.e.  $32.68({B_{24}B_{25}B_{26}})^{\frac{-1}{23}}+4B_{24}-\frac{B_{24}^3}{B_{25}B_{26}}>39.24$. But this is not true for $0.358\leq B_{24}\leq \frac{3}{2}B_{26}<0.436$ and $\frac{1}{2}B_{24}^2\leq B_{25}B_{26}<(0.31)(0.29)$. So we must have $B_{24}<0.358$.\vspace{2mm}\\
\noindent{\bf Claim(iv)} $B_i<max\{B_{12},B_{13},B_{14},2.1165\},~ {\mbox {for}}~ i=2,3,4,5,6,7,8,9,10,11$

 Using \eqref{2.2} and Lemmas 8, 9 we have\vspace{2mm}\\
 $~~~~~1=\lambda_1^{(26)}<B_1<\mu_1^{(26)}=4.2480446,$\\
 $~~~~~\frac{3}{4}=\lambda_2^{(26)}<B_1B_2<\mu_2^{(26)}=(1.1240391)^{24},$\\
 $~~~~~\frac{1}{2}=\lambda_3^{(26)}<B_1B_2B_3<\mu_3^{(26)}=(1.1939633)^{23},$\\
  $~~~~~\frac{1}{4}=\lambda_4^{(26)}<B_1B_2B_3B_4<\mu_4^{(26)}=(1.2699424)^{22},$\\
 $~~~~~\frac{\varepsilon}{4}=\lambda_5^{(26)}<B_1\cdots B_5<\mu_5^{(26)}=(1.3526843)^{21},$\\
 $~~~~~\frac{3\varepsilon^2}{16}=\lambda_6^{(26)}<B_1\cdots B_6<\mu_6^{(26)}=(1.4429963)^{20},$\\
 $~~~~~\frac{\varepsilon^3}{8}=\lambda_{7}^{(26)}<B_1\cdots B_7<\mu_7^{(26)}=(1.5418115)^{19},$\\
 $~~~~~\frac{\varepsilon^4}{16}=\lambda_{8}^{(26)}<B_1\cdots B_8<\mu_8^{(26)}=(1.650213)^{18},$\\
 $~~~~~\frac{\varepsilon^{28}}{4^4B_{26}^{17}}=\lambda_{9}^{(26)}<B_1\cdots B_9<\mu_9^{(24)}=(1.7694615)^{17},$\\
 $~~~~~\frac{\varepsilon^{24}}{4^4B_{26}^{16}}=\lambda_{10}^{(26)}<B_1\cdots B_{10}<\mu_{10}^{(26)}=(1.9010405)^{16},$\\
 $~~~~~\frac{\varepsilon^{21}}{2\times4^3B_{26}^{15}}=\lambda_{11}^{(26)}<B_1\cdots B_{11}<\mu_{11}^{(26)}=(2.0466947)^{15},$\\
 $~~~~~\frac{3\varepsilon^{18}}{4^4B_{26}^{14}}=\lambda_{12}^{(26)}<B_1\cdots B_{12}<\mu_{12}^{(26)}=(2.2084995)^{14},$\\
 $~~~~~\frac{\varepsilon^{15}}{4^3B_{26}^{13}}=\lambda_{13}^{(26)}<B_1\cdots B_{13}<\mu_{13}^{(26)}=(2.3889339)^{13},$\\
 $~~~~~\frac{\varepsilon^{12}}{4^3B_{26}^{12}}=\lambda_{14}^{(26)}<B_1\cdots B_{14}<\mu_{14}^{(26)}=(2.5909855)^{12}.$\vspace{2mm}\\
We find that $max\{\phi_{s,n-s}(\lambda_{s}^{(26)}),\phi_{s,n-s}(\mu_{s}^{(26)}), {\mbox {for}}~s=1,2,3,4,5,6,7,8,9,10\}=\phi_{1,25}(\mu_{1}^{(26)})$, which is $<\omega_{26}$. Therefore using Lemma 10 we have
 \begin{equation}\label{4.6.8}B_i<max\{B_{i+1},B_{i+2},\cdots,B_{26}\},~ {\mbox {for}}~ i=2,3,4,5,6,7,8,9,10,11.\end{equation}
From \eqref{2.0a},\eqref{2.0b} and Claims(i), (ii), (iii), we find $max\{B_{15},B_{16},\cdots,B_{26}\}<B_{15}\leq\frac{3}{2}\frac{B_{25}}{\varepsilon^2}<2.1165$. So using \eqref{4.6.8} we get that each of $B_2, \cdots,B_{11}$ is $<max\{B_{12},B_{13},B_{14},\frac{3}{2}\frac{B_{25}}{\varepsilon^2}\}<max\{B_{12},B_{13},B_{14},2.1165\}$.\vspace{2mm} \\
\noindent{\bf Claim(v)} $B_{12},B_{13},B_{14}<2.1165$

We find that for $s=13,12,11$, $\phi_{s,n-s}(\lambda_{s}^{(26)})<\omega_{26}$, but $\phi_{s,n-s}(\mu_{s}^{(26)})>\omega_{26}$, so we apply Lemma 11 respectively with $\sigma_{13}^{(26)}=(2.274)^{13}$, $\sigma_{12}^{(26)}=(2.15)^{14}$ and $\sigma_{11}^{(26)}=(2.036)^{15}$. Here $\phi_{13,13}(\sigma_{13}^{(26)})<\omega_{26}$, $\phi_{12,14}(\sigma_{12}^{(26)})<\omega_{26}$ and $\phi_{11,15}(\sigma_{11}^{(26)})<\omega_{26}$.\\
First consider Lemma 11 for $s=13$ and with $\sigma_{13}^{(26)}=(2.274)^{13}$.\\
In Case(i), when $B_1B_2\cdots B_{13}<(2.274)^{13}$, then we have $B_{14}<max\{B_{15},B_{16},\cdots,B_{26}\}$, which is $<\frac{3}{2}\frac{B_{25}}{\varepsilon^2}<2.1165$.\\
In Case(ii), when $B_1B_2\cdots B_{13}\geq(2.274)^{13}$, then we have $B_{14}<\frac{\mu_{14}^{(26)}}{\sigma_{13}^{(26)}}<\frac{(2.591)^{12}}{(2.274)^{13}}<2.1055$. \\
So we have $B_{14}<\frac{3}{2}\frac{B_{25}}{\varepsilon^2}<2.1165$.\\
Now consider Lemma 11 for $s=12$ and with $\sigma_{12}^{(26)}=(2.15)^{14}$.\\
In Case(i), when $B_1B_2\cdots B_{12}<(2.15)^{14}$, then we have $B_{13}<max\{B_{14},B_{15},\cdots,B_{26}\}$, which is $<\frac{3}{2}\frac{B_{25}}{\varepsilon^2}<2.1165$.\\
In Case(ii), when $B_1B_2\cdots B_{12}\geq(2.15)^{14}$, then we have $B_{13}<\frac{\mu_{13}^{(26)}}{\sigma_{12}^{(26)}}<\frac{(2.389)^{13}}{(2.15)^{14}}<1.831$. \\
So we have $B_{13}<\frac{3}{2}\frac{B_{25}}{\varepsilon^2}<2.1165$.\\
Now consider Lemma 11 for $s=11$ and with $\sigma_{11}^{(26)}=(2.036)^{15}$.\\
In Case(i), when $B_1B_2\cdots B_{11}<(2.036)^{15}$, then we have $B_{12}<max\{B_{13},B_{14},\cdots,B_{26}\}$, which is $<\frac{3}{2}\frac{B_{25}}{\varepsilon^2}<2.1165$.\\
In Case(ii), when $B_1B_2\cdots B_{11}\geq(2.036)^{15}$, then we have $B_{12}<\frac{\mu_{12}^{(26)}}{\sigma_{11}^{(26)}}<\frac{(2.2085)^{14}}{(2.036)^{15}}<1.54$. \\
So we have $B_{12}<B_{13}<\frac{3}{2}\frac{B_{25}}{\varepsilon^2}<2.1165$.\vspace{2mm}\\
 Using Claims(iv) and (v) we get each of $B_2, B_3,\cdots,B_{15}$ is $<2.1165$. \vspace{2mm}\\
 \noindent{\bf Final Contradiction}

 Now $2B_2+2B_4+\cdots+2B_{26}<2(7\times2.1165)+2\{(\frac{1}{\varepsilon^2})+\frac{3/2}{\varepsilon}+\frac{1}{\varepsilon}+\frac{3}{2}+1\}B_{24}
+2B_{26}<39.24$ for $B_{26}<0.29$ and $B_{24}<0.358$, giving thereby a contradiction to the weak inequality $(2,2,\cdots,2,2)_w$.~$\Box$
\subsection{$n=27$}
Here we have $\omega_{27}=41.78$, $B_1\leq\gamma_{27}<4.3684312$. Using \eqref{2.3}, we have $l_{27}=0.1231<B_{27}<4.0264547=m_{27}$. \vspace{2mm}\\
\noindent{\bf Claim(i)} $B_{27}<0.226$

Suppose $B_{27}\geq0.226$. The inequality $(26^*,1)$ gives $39.24(B_{27})^{\frac{-1}{26}}+B_{27}>41.78$. But this is not true for $0.226\leq B_{27}\leq 4.0264547$. So we must have $B_{27}<0.226$.\vspace{2mm}\\
\noindent{\bf Claim(ii)} $B_i<max\{B_{13},B_{14},2.195\},~ {\mbox {for}}~ i=2,3,4,5,6,7,8,9,10,11,12$

 Using \eqref{2.2} and Lemmas 8, 9 we have\vspace{2mm}\\
 $~~~~~1=\lambda_1^{(27)}<B_1<\mu_1^{(27)}=4.3684312,$\\
 $~~~~~\frac{3}{4}=\lambda_2^{(27)}<B_1B_2<\mu_2^{(27)}=(1.1213882)^{25},$\\
 $~~~~~\frac{1}{2}=\lambda_3^{(27)}<B_1B_2B_3<\mu_3^{(27)}=(1.1896237)^{24},$\\
  $~~~~~\frac{1}{4}=\lambda_4^{(27)}<B_1B_2B_3B_4<\mu_4^{(27)}=(1.2636277)^{23},$\\
 $~~~~~\frac{\varepsilon}{4}=\lambda_5^{(27)}<B_1\cdots B_5<\mu_5^{(27)}=(1.34403991)^{22},$\\
 $~~~~~\frac{3\varepsilon^2}{16}=\lambda_6^{(27)}<B_1\cdots B_6<\mu_6^{(27)}=(1.4316096)^{21},$\\
 $~~~~~\frac{\varepsilon^3}{8}=\lambda_{7}^{(27)}<B_1\cdots B_7<\mu_7^{(27)}=(1.5271911)^{20},$\\
$~~~~~\frac{\varepsilon^{36}}{2\times4^4B_{27}^{19}}=\lambda_{8}^{(27)}<B_1\cdots B_8<\mu_8^{(27)}=(1.6317719)^{19},$\\
 $~~~~~\frac{3\varepsilon^{32}}{4^5B_{27}^{18}}=\lambda_{9}^{(27)}<B_1\cdots B_9<\mu_9^{(24)}=(1.7464974)^{18},$\\
 $~~~~~\frac{\varepsilon^{28}}{4^4B_{27}^{17}}=\lambda_{10}^{(27)}<B_1\cdots B_{10}<\mu_{10}^{(27)}=(1.8727046)^{17},$\\
 $~~~~~\frac{\varepsilon^{24}}{4^4B_{27}^{16}}=\lambda_{11}^{(27)}<B_1\cdots B_{11}<\mu_{11}^{(27)}=(2.0119609)^{16},$\\
 $~~~~~\frac{\varepsilon^{21}}{2\times4^3B_{27}^{15}}=\lambda_{12}^{(27)}<B_1\cdots B_{12}<\mu_{12}^{(27)}=(2.1661136)^{15},$\\
 $~~~~~\frac{3\varepsilon^{18}}{4^4B_{27}^{14}}=\lambda_{13}^{(27)}<B_1\cdots B_{13}<\mu_{13}^{(27)}=(2.3373592)^{14},$\\
 $~~~~~\frac{\varepsilon^{15}}{4^3B_{27}^{13}}=\lambda_{14}^{(27)}<B_1\cdots B_{14}<\mu_{14}^{(27)}=(2.5283215)^{13}.$\vspace{2mm}\\
We find that $max\{\phi_{s,n-s}(\lambda_{s}^{(27)}),\phi_{s,n-s}(\mu_{s}^{(27)}), {\mbox {for}}~s=1,2,3,4,5,6,7,8,9,10,11\}=\phi_{11,16}(\mu_{11}^{(27)})$, which is $<\omega_{27}$. Therefore using Lemma 10 we have
 \begin{equation}\label{4.6.9}B_i<max\{B_{i+1},B_{i+2},\cdots,B_{27}\},~ {\mbox {for}}~ i=2,3,4,5,6,7,8,9,10,11,12.\end{equation}
From \eqref{2.0a},\eqref{2.0b} and Claim(i) we find that, $max\{B_{15},B_{16},\cdots,B_{27}\}<B_{15}\leq\frac{B_{27}}{\varepsilon^3}<2.195$. So using \eqref{4.6.9} we get that each of $B_2, \cdots,B_{12}$ is $<max\{B_{13},B_{14},\frac{B_{27}}{\varepsilon^3}\}<max\{B_{13},B_{14},2.195\}$.\vspace{2mm} \\
\noindent{\bf Claim(iii)} $B_{13},~B_{14}<2.195$

We find that for $s=13,12$, $\phi_{s,n-s}(\lambda_{s}^{(27)})<\omega_{27}$, but $\phi_{s,n-s}(\mu_{s}^{(27)})>\omega_{27}$, so we apply Lemma 11 respectively with $\sigma_{13}^{(27)}=(2.25)^{14}$ and $\sigma_{12}^{(27)}=(2.13)^{15}$. Here $\phi_{13,14}(\sigma_{13}^{(27)})<\omega_{27}$ and $\phi_{12,15}(\sigma_{12}^{(27)})<\omega_{27}$.\\
First consider Lemma 11 for $s=13$ and with $\sigma_{13}^{(27)}=(2.25)^{14}$.\\
In Case(i), when $B_1B_2\cdots B_{13}<(2.25)^{14}$, then we have $B_{14}<max\{B_{15},B_{16},\cdots,B_{27}\}$, which is $<\frac{B_{27}}{\varepsilon^3}<2.195$.\\
In Case(ii), when $B_1B_2\cdots B_{13}\geq(2.25)^{14}$, then we have $B_{14}<\frac{\mu_{14}^{(27)}}{\sigma_{13}^{(27)}}<\frac{(2.5284)^{13}}{(2.25)^{14}}<2.026$. \\
So we have $B_{14}<\frac{B_{27}}{\varepsilon^3}<2.195$.\\
Now consider Lemma 11 for $s=12$ and with $\sigma_{12}^{(27)}=(2.13)^{15}$.\\
In Case(i), when $B_1B_2\cdots B_{12}<(2.13)^{15}$, then we have $B_{13}<max\{B_{14},B_{15},\cdots,B_{27}\}$, which is $<\frac{B_{27}}{\varepsilon^3}<2.195$.\\
In Case(ii), when $B_1B_2\cdots B_{12}\geq(2.13)^{15}$, then we have $B_{13}<\frac{\mu_{13}^{(27)}}{\sigma_{12}^{(27)}}<\frac{(2.3374)^{14}}{(2.13)^{15}}<1.725$. \\
So we have $B_{13}<\frac{B_{27}}{\varepsilon^3}<2.195$.\vspace{2mm}\\
Using Claim(ii) and (iii) we get each of $B_2, B_3,\cdots,B_{15}$ is $<2.195$.\vspace{2mm}\\
 \noindent{\bf Final Contradiction}

 Now
$2B_2+B_3+2B_5+\cdots+2B_{27}<
3(\frac{B_{27}}{\varepsilon^3})+2\{
(6\times\frac{1}{\varepsilon^3}+\frac{3/2}{\varepsilon^2}+\frac{1}{\varepsilon^2}+\frac{3/2}{\varepsilon}+\frac{1}{\varepsilon}+\frac{3}{2}+1)\}
B_{27}<41.78$ ~for $B_{27}<0.226$, giving thereby a contradiction to the weak inequality $(2,1,2,\cdots,2,2)_w$.~$\Box$
\subsection{$n=28$}
Here we have $\omega_{28}=44.36$, $B_1\leq\gamma_{28}<4.488631$. Using \eqref{2.3}, we have $l_{28}=0.1164<B_{28}<4.144353=m_{28}$. \vspace{2mm}\\
 \noindent{\bf Claim(i)} $B_{28}<0.228$

Suppose $B_{28}\geq0.228$. The inequality $(27^*,1)$ gives $41.78(B_{28})^{\frac{-1}{27}}+B_{28}>44.36$. But this is not true for $0.228\leq B_{28}\leq 4.144353$. So we must have $B_{28}<0.228$.\vspace{2mm}\\
 \noindent{\bf Claim(ii)} $B_i<max\{B_{13},B_{14},B_{15},2.215\},~ {\mbox {for}}~ i=2,3,4,5,6,7,8,9,10,11,12$
 Using \eqref{2.2} and Lemmas 8, 9 we have\vspace{2mm}\\
 $~~~~~1=\lambda_1^{(28)}<B_1<\mu_1^{(28)}=4.488631,$\\
 $~~~~~\frac{3}{4}=\lambda_2^{(28)}<B_1B_2<\mu_2^{(28)}=(1.118873)^{26},$\\
 $~~~~~\frac{1}{2}=\lambda_3^{(28)}<B_1B_2B_3<\mu_3^{(28)}=(1.1855192)^{25},$\\
  $~~~~~\frac{1}{4}=\lambda_4^{(28)}<B_1B_2B_3B_4<\mu_4^{(28)}=(1.257657)^{24},$\\
 $~~~~~\frac{\varepsilon}{4}=\lambda_5^{(28)}<B_1\cdots B_5<\mu_5^{(28)}=(1.3358932)^{23},$\\
 $~~~~~\frac{3\varepsilon^2}{16}=\lambda_6^{(28)}<B_1\cdots B_6<\mu_6^{(28)}=(1.4209042)^{22},$\\
 $~~~~~\frac{\varepsilon^3}{8}=\lambda_{7}^{(28)}<B_1\cdots B_7<\mu_7^{(28)}=(1.5134819)^{21},$\\
 $~~~~~\frac{\varepsilon^{40}}{4^5 B_{28}^{20}}=\lambda_{8}^{(28)}<B_1\cdots B_8<\mu_8^{(28)}=(1.6145296)^{20},$\\
 $~~~~~\frac{\varepsilon^{36}}{2\times 4^4 B_{28}^{19}}=\lambda_{9}^{(28)}<B_1\cdots B_9<\mu_9^{(24)}=(1.7250912)^{19},$\\
 $~~~~~\frac{3\varepsilon^{32}}{4^{5}B_{28}^{18}}=\lambda_{10}^{(28)}<B_1\cdots B_{10}<\mu_{10}^{(28)}=(1.8463778)^{18},$\\
 $~~~~~\frac{\varepsilon^{28}}{4^4B_{28}^{17}}=\lambda_{11}^{(28)}<B_1\cdots B_{11}<\mu_{11}^{(28)}=(1.9798026)^{17},$\\
 $~~~~~\frac{\varepsilon^{24}}{4^4B_{28}^{16}}=\lambda_{12}^{(28)}<B_1\cdots B_{12}<\mu_{12}^{(28)}=(2.1270228)^{16},$\\
 $~~~~~\frac{\varepsilon^{21}}{2\times 4^3B_{28}^{15}}=\lambda_{13}^{(28)}<B_1\cdots B_{13}<\mu_{13}^{(28)}=(2.2899914)^{15},$\\
 $~~~~~\frac{3\varepsilon^{18}}{4^4B_{28}^{14}}=\lambda_{14}^{(28)}<B_1\cdots B_{14}<\mu_{14}^{(28)}=(2.4710303)^{14},$\\
 $~~~~~\frac{\varepsilon^{15}}{4^3B_{28}^{13}}=\lambda_{15}^{(28)}<B_1\cdots B_{15}<\mu_{15}^{(28)}=(2.6729135)^{13}.$\vspace{2mm}\\
We find that $max\{\phi_{s,n-s}(\lambda_{s}^{(28)}),\phi_{s,n-s}(\mu_{s}^{(28)}), {\mbox {for}}~s=1,2,3,4,5,6,7,8,9,10,11\}=\phi_{11,17}(\mu_{11}^{(28)})$, which is $<\omega_{28}$. Therefore using Lemma 10 we have
 \begin{equation}\label{4.6.10}B_i<max\{B_{i+1},B_{i+2},\cdots,B_{28}\},~ {\mbox {for}}~ i=2,3,4,5,6,7,8,9,10,11,12.\end{equation}
From \eqref{2.0a},\eqref{2.0b} and Claim(i) we find that, $max\{B_{16},B_{17},\cdots,B_{28}\}<B_{16}\leq\frac{B_{28}}{\varepsilon^3}<2.215$. So using \eqref{4.6.10} we get each of $B_2, \cdots,B_{12}$ is $<max\{B_{13},B_{14},B_{15},\frac{B_{28}}{\varepsilon^3}\}<max\{B_{13},B_{14},B_{15},2.215\}$.\vspace{2mm} \\
 \noindent{\bf Claim(iii)} $B_{13},B_{14},B_{15}<2.215$

We find that for $s=14,13,12$, $\phi_{s,n-s}(\lambda_{s}^{(28)})<\omega_{28}$, but $\phi_{s,n-s}(\mu_{s}^{(28)})>\omega_{28}$, so we apply Lemma 11 respectively with $\sigma_{14}^{(28)}=(2.355)^{14}$, $\sigma_{13}^{(28)}=(2.23)^{15}$ and $\sigma_{12}^{(28)}=(2.11)^{16}$. Here $\phi_{14,14}(\sigma_{14}^{(28)})<\omega_{28}$, $\phi_{13,15}(\sigma_{13}^{(28)})<\omega_{28}$ and $\phi_{12,16}(\sigma_{12}^{(28)})<\omega_{28}$.\\
First consider Lemma 11 for $s=14$ and with $\sigma_{14}^{(28)}=(2.355)^{14}$.\\
In Case(i), when $B_1B_2\cdots B_{14}<(2.355)^{14}$, then we have $B_{15}<max\{B_{16},B_{17},\cdots,B_{28}\}$, which is $<\frac{B_{28}}{\varepsilon^3}<2.215$.\\
In Case(ii), when $B_1B_2\cdots B_{14}\geq(2.355)^{14}$, then we have $B_{15}<\frac{\mu_{15}^{(28)}}{\sigma_{14}^{(28)}}<\frac{(2.673)^{13}}{(2.355)^{14}}<2.204$. \\
So we have $B_{15}<\frac{B_{28}}{\varepsilon^3}<2.215$.\\
Now consider Lemma 11 for $s=13$ and with $\sigma_{13}^{(28)}=(2.23)^{15}$.\\
In Case(i), when $B_1B_2\cdots B_{13}<(2.23)^{15}$, then we have $B_{14}<max\{B_{15},B_{16},\cdots,B_{28}\}$, which is $<\frac{B_{28}}{\varepsilon^3}<2.215$.\\
In Case(ii), when $B_1B_2\cdots B_{13}\geq(2.23)^{15}$, then we have $B_{14}<\frac{\mu_{14}^{(28)}}{\sigma_{13}^{(28)}}<\frac{(2.472)^{14}}{(2.23)^{15}}<1.898$. \\
So we have $B_{14}<\frac{B_{28}}{\varepsilon^3}<2.215$.\\
Now consider Lemma 11 for $s=12$ and with $\sigma_{12}^{(28)}=(2.11)^{16}$.\\
In Case(i), when $B_1B_2\cdots B_{12}<(2.11)^{16}$, then we have $B_{13}<max\{B_{14},B_{15},\cdots,B_{28}\}$, which is $<\frac{B_{28}}{\varepsilon^3}<2.215$.\\
In Case(ii), when $B_1B_2\cdots B_{12}\geq(2.11)^{16}$, then we have $B_{13}<\frac{\mu_{13}^{(28)}}{\sigma_{12}^{(28)}}<\frac{(2.2899)^{15}}{(2.11)^{16}}<1.62$. \\
So we have $B_{13}<\frac{B_{28}}{\varepsilon^3}<2.215$.\\
Using Claim(ii) and (iii) we get each of $B_2, B_3,\cdots,B_{16}$ is $<2.215$.\vspace{2mm}\\
 \noindent{\bf Final Contradiction}

Now
$2B_2+2B_4+\cdots+2B_{28}<
2\{8\times\frac{1}{\varepsilon^3}+\frac{3/2}{\varepsilon^2}+\frac{1}{\varepsilon^2}+\frac{3/2}{\varepsilon}+\frac{1}{\varepsilon}+\frac{3}{2}+1\}
B_{28}<44.36$ for $B_{28}<0.228$, giving thereby a contradiction to the weak inequality $(2,2,\cdots,2,2)_w$.~$\Box$
\subsection{$n=29$}
Here $\omega_{29}=47.18$, $B_1\leq\gamma_{29}<4.6086676$. Using \eqref{2.3}, we have $l_{29}=0.1102<B_{29}<4.2621353=m_{29}$.\vspace{2mm}\\
 \noindent{\bf Claim(i)} $B_{29}<0.201$

The inequality $(28^*,1)$ gives $44.36(B_{29})^{\frac{-1}{28}}+B_{29}>47.18$. But this is not true for $0.201\leq B_{29}\leq4.2621353$. So we must have $B_{29}<0.201$.\vspace{2mm}\\
 \noindent{\bf Claim(ii)} $B_{28}<0.246$

The inequality $(27^*,2)$ gives $41.78({B_{28}B_{29}})^{\frac{-1}{27}}+4B_{28}-\frac{2B_{28}^2}{B_{29}}>47.18$. But this is not true for $0.246\leq B_{28}\leq \frac{4}{3}B_{29}<0.268$ and $\frac{3}{4}(B_{28})\leq B_{29}<0.201$. So we must have $B_{28}<0.246$.\vspace{2mm}\\
 \noindent{\bf Claim(iii)} $B_{27}<0.2835$

The inequality $(26^*,3)$ gives $39.24({B_{27}B_{28}B_{29}})^{\frac{-1}{26}}+4B_{27}-\frac{B_{27}^3}{B_{28}B_{29}}>47.18$. But this is not true for $0.2835\leq B_{27}\leq \frac{3}{2}B_{29}<0.302$ and $\frac{1}{2}B_{27}^2\leq B_{28}B_{29}<(0.246)(0.201)$. So we must have $B_{27}<0.2835$.\vspace{2mm}\\
 \noindent{\bf Claim(iv)} $B_i<max\{B_{14},B_{15},2.3888\},~ {\mbox {for}}~ i=2,3,4,5,6,7,8,9,10,11,12,13$

 Using \eqref{2.2} and Lemmas 8, 9 we have\vspace{2mm}\\
 $~~~~~1=\lambda_1^{(29)}<B_1<\mu_1^{(29)}=4.6086676,$\\
 $~~~~~\frac{3}{4}=\lambda_2^{(29)}<B_1B_2<\mu_2^{(29)}=(1.1164824)^{27},$\\
 $~~~~~\frac{1}{2}=\lambda_3^{(29)}<B_1B_2B_3<\mu_3^{(29)}=(1.18162593)^{26},$\\
  $~~~~~\frac{1}{4}=\lambda_4^{(29)}<B_1B_2B_3B_4<\mu_4^{(29)}=(1.2520101)^{25},$\\
 $~~~~~\frac{\varepsilon}{4}=\lambda_5^{(29)}<B_1\cdots B_5<\mu_5^{(29)}=(1.3281938)^{24},$\\
 $~~~~~\frac{3\varepsilon^2}{16}=\lambda_6^{(29)}<B_1\cdots B_6<\mu_6^{(29)}=(1.41081792)^{23},$\\
 $~~~~~\frac{\varepsilon^3}{8}=\lambda_{7}^{(29)}<B_1\cdots B_7<\mu_7^{(29)}=(1.5005968)^{22},$\\
 $~~~~~\frac{\varepsilon^4}{16}=\lambda_{8}^{(29)}<B_1\cdots B_8<\mu_8^{(29)}=(1.5983668)^{21},$\\
 $~~~~~\frac{\varepsilon^{40}}{4^5B_{29}^{20}}=\lambda_{9}^{(29)}<B_1\cdots B_9<\mu_9^{(24)}=(1.7050819)^{20},$\\
 $~~~~~\frac{\varepsilon^{36}}{2\times4^{4}B_{29}^{19}}=\lambda_{10}^{(29)}<B_1\cdots B_{10}<\mu_{10}^{(29)}=(1.8218444)^{19},$\\
 $~~~~~\frac{3\varepsilon^{32}}{4^5B_{29}^{18}}=\lambda_{11}^{(29)}<B_1\cdots B_{11}<\mu_{11}^{(29)}=(1.9499335)^{18},$\\
 $~~~~~\frac{\varepsilon^{28}}{4^4B_{29}^{17}}=\lambda_{12}^{(29)}<B_1\cdots B_{12}<\mu_{12}^{(29)}=(2.0908415)^{17},$\\
 $~~~~~\frac{\varepsilon^{24}}{4^4B_{29}^{16}}=\lambda_{13}^{(29)}<B_1\cdots B_{13}<\mu_{13}^{(29)}=(2.2463188)^{16},$\\
 $~~~~~\frac{\varepsilon^{21}}{2\times4^3B_{29}^{15}}=\lambda_{14}^{(29)}<B_1\cdots B_{14}<\mu_{14}^{(29)}=(2.4184275)^{15},$\\
 $~~~~~\frac{3\varepsilon^{18}}{4^4B_{29}^{14}}=\lambda_{15}^{(29)}<B_1\cdots B_{15}<\mu_{15}^{(29)}=(2.6096202)^{14}.$\vspace{2mm}\\
We find that $max\{\phi_{s,n-s}(\lambda_{s}^{(29)}),\phi_{s,n-s}(\mu_{s}^{(29)}), {\mbox {for}}~s=1,2,3,4,5,6,7,8,9,10,11,12\}=\phi_{12,17}(\mu_{12}^{(29)})$, which is $<\omega_{29}$. Therefore using Lemma 10 we have
 \begin{equation}\label{4.6.11}B_i<max\{B_{i+1},B_{i+2},\cdots,B_{29}\},~ {\mbox {for}}~ i=2,3,4,5,6,7,8,9,10,11,12,13.\end{equation}
From \eqref{2.0a},\eqref{2.0b} and Claims (i),(ii) and (iii) we find that , $max\{B_{16},B_{17},\cdots,B_{29}\}<B_{16}\leq\frac{B_{28}}{\varepsilon^3}<2.3888$. So using \eqref{4.6.11} we get that each of $B_2, \cdots,B_{13}$ is $<max\{B_{14},B_{15},\frac{B_{28}}{\varepsilon^3}\}<max\{B_{14},B_{15},2.3888\}$. \\
 \noindent{\bf Claim(v)} $B_{14},B_{15}<2.3888$

We find that for $s=14,13$, $\phi_{s,n-s}(\lambda_{s}^{(29)})<\omega_{29}$, but $\phi_{s,n-s}(\mu_{s}^{(29)})>\omega_{29}$, so we apply Lemma 11 respectively with $\sigma_{14}^{(29)}=(2.34)^{15}$ and $\sigma_{13}^{(29)}=(2.21)^{16}$. Here $\phi_{14,15}(\sigma_{14}^{(29)})<\omega_{29}$ and $\phi_{13,16}(\sigma_{13}^{(29)})<\omega_{29}$.\\
First consider Lemma 11 for $s=14$ and with $\sigma_{14}^{(29)}=(2.34)^{15}$.\\
In Case(i), when $B_1B_2\cdots B_{14}<(2.34)^{15}$, then we have $B_{15}<max\{B_{16},B_{17},\cdots,B_{29}\}$, which is $<\frac{B_{28}}{\varepsilon^3}<2.3888$.\\
In Case(ii), when $B_1B_2\cdots B_{14}\geq(2.34)^{15}$, then we have $B_{15}<\frac{\mu_{15}^{(29)}}{\sigma_{14}^{(29)}}<\frac{(2.6097)^{14}}{(2.34)^{15}}<1.969$. \\
So we have $B_{15}<\frac{B_{28}}{\varepsilon^3}<2.3888$.\\
Now consider Lemma 11 for $s=13$ and with $\sigma_{13}^{(29)}=(2.21)^{16}$.\\
In Case(i), when $B_1B_2\cdots B_{13}<(2.21)^{16}$, then we have $B_{14}<max\{B_{15},B_{16},\cdots,B_{29}\}$, which is $<\frac{B_{28}}{\varepsilon^3}<2.3888$.\\
In Case(ii), when $B_1B_2\cdots B_{13}\geq(2.21)^{16}$, then we have $B_{14}<\frac{\mu_{14}^{(29)}}{\sigma_{13}^{(29)}}<\frac{(2.4185)^{15}}{(2.21)^{16}}<1.75$. \\
So we have $B_{14}<\frac{B_{28}}{\varepsilon^3}<2.3888$.\\
Using Claims(iv) and (v) we get each of $B_2, B_3,\cdots,B_{16}$ is $<2.3888$.\vspace{2mm}\\
 \noindent{\bf Final Contradiction}

Now
$2B_2+B_3+2B_5+\cdots+2B_{29}<
3(\frac{B_{28}}{\varepsilon^3})+2
(6\times\frac{B_{28}}{\varepsilon^3})+2(\frac{3/2}{\varepsilon^2}+\frac{1}{\varepsilon^2}+\frac{3/2}{\varepsilon}+\frac{1}{\varepsilon}+\frac{3}{2}+1)B_{27}+2
B_{29}<47.18$ for $B_{29}<0.201$, $B_{28}<0.246$ and $B_{27}<0.2835$, giving thereby a contradiction to the weak inequality $(2,1,2,\cdots,2,2)_w$.~$\Box$
\subsection{$n=30$}
Here $\omega_{30}=49.86$. We have $B_1\leq\gamma_{30}<4.7285667$. Using \eqref{2.3}, we have $l_{30}=0.1045<B_{30}<4.3798196=m_{30}$.\vspace{2mm}\\
\noindent{\bf Claim(i)} $B_{30}<0.231$

Suppose $B_{30}\geq 0.231$.
The inequality $(29^*,1)$ gives $47.18(B_{30})^{\frac{-1}{29}}+B_{30}>49.86$. But this is not true for $0.231\leq B_{30}\leq 4.3798196$. \vspace{2mm}\\
\noindent{\bf Claim(ii)} $B_{29}<0.247$

Suppose $B_{29}\geq0.247$, then $2B_{29}>B_{30}$. Therefore the inequality $(28^*,2)$ holds i.e. $44.36({B_{29}B_{30}})^{\frac{-1}{28}}+4B_{29}-\frac{2B_{29}^2}{B_{30}}>49.86$. But this is not true for $0.247\leq B_{29}\leq \frac{4}{3}B_{30}<0.308$ and $\frac{3}{4}(B_{29})\leq B_{30}<0.231$. So we must have $B_{29}<0.247$.\vspace{2mm}\\
\noindent{\bf Claim(iii)} $B_{28}<0.285$

Suppose $B_{28}\geq0.285$, then $B_{28}^2>B_{29}B_{30}$. Therefore the inequality $(27^*,3)$ holds, i.e. $41.78({B_{28}B_{29}B_{30}})^{\frac{-1}{27}}+4B_{28}-\frac{B_{28}^3}{B_{29}B_{30}}>49.86$. But this is not true for $0.285\leq B_{28}\leq \frac{3}{2}B_{29}<0.347$ and $\frac{1}{2}B_{28}^2<B_{29}B_{30}<(0.247)(0.231)$. So we must have $B_{28}<0.285$.\vspace{2mm}\\
\noindent{\bf Claim(iv)} $B_i<max\{B_{14},B_{15},B_{16},2.3985\},~ {\mbox {for}}~ i=2,3,4,5,6,7,8,9,10,11,12,13$

 Using \eqref{2.2} and Lemmas 8, 9 we have\vspace{2mm}\\
 $~~~~~1=\lambda_1^{(30)}<B_1<\mu_1^{(30)}=4.7285667,$\\
 $~~~~~\frac{3}{4}=\lambda_2^{(30)}<B_1B_2<\mu_2^{(30)}=(1.114207)^{28},$\\
 $~~~~~\frac{1}{2}=\lambda_3^{(30)}<B_1B_2B_3<\mu_3^{(30)}=(1.1779273)^{27},$\\
  $~~~~~\frac{1}{4}=\lambda_4^{(30)}<B_1B_2B_3B_4<\mu_4^{(30)}=(1.2466561)^{26},$\\
 $~~~~~\frac{\varepsilon}{4}=\lambda_5^{(30)}<B_1\cdots B_5<\mu_5^{(30)}=(1.3209137)^{25},$\\
 $~~~~~\frac{3\varepsilon^2}{16}=\lambda_6^{(30)}<B_1\cdots B_6<\mu_6^{(30)}=(1.4012902)^{24},$\\
 $~~~~~\frac{\varepsilon^3}{8}=\lambda_{7}^{(30)}<B_1\cdots B_7<\mu_7^{(30)}=(1.4884616)^{23},$\\
 $~~~~~\frac{\varepsilon^4}{16}=\lambda_{8}^{(30)}<B_1\cdots B_8<\mu_8^{(30)}=(1.5831813)^{22},$\\
 $~~~~~\frac{\varepsilon^6}{16}=\lambda_{9}^{(30)}<B_1\cdots B_9<\mu_9^{(24)}=(1.6863321)^{21},$\\
 $~~~~~\frac{\varepsilon^{40}}{4^{5}B_{30}^{20}}=\lambda_{10}^{(30)}<B_1\cdots B_{10}<\mu_{10}^{(30)}=(1.7989201)^{20},$\\
 $~~~~~\frac{\varepsilon^{36}}{2\times4^4B_{30}^{19}}=\lambda_{11}^{(30)}<B_1\cdots B_{11}<\mu_{11}^{(30)}=(1.9221086)^{19},$\\
 $~~~~~\frac{3\varepsilon^{32}}{4^5B_{30}^{18}}=\lambda_{12}^{(30)}<B_1\cdots B_{12}<\mu_{12}^{(30)}=(2.057247)^{18},$\\
 $~~~~~\frac{\varepsilon^{28}}{4^4B_{30}^{17}}=\lambda_{13}^{(30)}<B_1\cdots B_{13}<\mu_{13}^{(30)}=(2.2059099)^{17},$\\
 $~~~~~\frac{\varepsilon^{24}}{4^4B_{30}^{16}}=\lambda_{14}^{(30)}<B_1\cdots B_{14}<\mu_{14}^{(30)}=(2.3699437)^{16},$\\
 $~~~~~\frac{\varepsilon^{21}}{2\times4^3B_{30}^{15}}=\lambda_{15}^{(30)}<B_1\cdots B_{15}<\mu_{15}^{(30)}=(2.551525)^{15},$\\
  $~~~~~\frac{3\varepsilon^{18}}{4^4B_{30}^{14}}=\lambda_{16}^{(30)}<B_1\cdots B_{16}<\mu_{16}^{(30)}=(2.7532392)^{14}.$\vspace{2mm}\\
We find that $max\{\phi_{s,n-s}(\lambda_{s}^{(30)}),\phi_{s,n-s}(\mu_{s}^{(30)}), {\mbox {for}}~s=1,2,3,4,5,6,7,8,9,10,11,12\}=\phi_{12,18}(\mu_{12}^{(30)})$, which is $<\omega_{30}$. Therefore using Lemma 10 we have
 \begin{equation}\label{4.6.12}B_i<max\{B_{i+1},B_{i+2},\cdots,B_{29}\},~ {\mbox {for}}~ i=2,3,4,5,6,7,8,9,10,11,12,13.\end{equation}
From \eqref{2.0a},\eqref{2.0b} and Claims (i),(ii) and (iii) we find that, $max\{B_{17},B_{18},\cdots,B_{29}\}<B_{17}\leq\frac{B_{29}}{\varepsilon^3}<2.3985$. So using \eqref{4.6.12} we get that  each of $B_2, \cdots,B_{13}$ is $<max\{B_{14},B_{15},B_{16},\frac{B_{29}}{\varepsilon^3}\}<max\{B_{14},B_{15},B_{16},2.3985\}$. \\
\noindent{\bf Claim(v)} $B_{14},B_{15},B_{16}<2.3985$

We find that for $s=15,14,13$, $\phi_{s,n-s}(\lambda_{s}^{(30)})<\omega_{30}$, but $\phi_{s,n-s}(\mu_{s}^{(30)})>\omega_{30}$, so we apply Lemma 11 respectively with $\sigma_{15}^{(30)}=(2.438)^{15}$, $\sigma_{14}^{(30)}=(2.305)^{16}$ and $\sigma_{13}^{(30)}=(2.185)^{17}$. Here $\phi_{15,15}(\sigma_{15}^{(30)})<\omega_{30}$, $\phi_{14,16}(\sigma_{14}^{(30)})<\omega_{30}$ and $\phi_{13,17}(\sigma_{13}^{(30)})<\omega_{30}$.\\
First consider Lemma 11 for $s=15$ and with $\sigma_{15}^{(30)}=(2.438)^{15}$.\\
In Case(i), when $B_1B_2\cdots B_{15}<(2.438)^{15}$, then we have $B_{16}<max\{B_{17},B_{18},\cdots,B_{30}\}$, which is $<\frac{B_{29}}{\varepsilon^3}<2.3985$.\\
In Case(ii), when $B_1B_2\cdots B_{15}\geq(2.438)^{15}$, then we have $B_{16}<\frac{\mu_{16}^{(30)}}{\sigma_{15}^{(30)}}<\frac{(2.7533)^{14}}{(2.438)^{15}}<2.26$. \\
So we have $B_{16}<\frac{B_{29}}{\varepsilon^3}<2.3985$.\\
Now consider Lemma 11 for $s=14$ and with $\sigma_{14}^{(30)}=(2.305)^{16}$.\\
In Case(i), when $B_1B_2\cdots B_{14}<(2.305)^{16}$, then we have $B_{15}<max\{B_{16},B_{17},\cdots,B_{30}\}$, which is $<\frac{B_{29}}{\varepsilon^3}<2.3985$.\\
In Case(ii), when $B_1B_2\cdots B_{14}\geq(2.305)^{16}$, then we have $B_{15}<\frac{\mu_{15}^{(30)}}{\sigma_{14}^{(30)}}<\frac{(2.5516)^{15}}{(2.305)^{16}}<1.993$. \\
So we have $B_{15}<\frac{B_{29}}{\varepsilon^3}<2.3985$.\\
Further consider Lemma 11 for $s=13$ and with $\sigma_{13}^{(30)}=(2.185)^{17}$.\\
In Case(i), when $B_1B_2\cdots B_{13}<(2.185)^{17}$, then we have $B_{14}<max\{B_{15},B_{16},\cdots,B_{30}\}$, which is $<\frac{B_{29}}{\varepsilon^3}<2.3985$.\\
In Case(ii), when $B_1B_2\cdots B_{13}\geq(2.185)^{17}$, then we have $B_{14}<\frac{\mu_{14}^{(30)}}{\sigma_{13}^{(30)}}<\frac{(2.36995)^{16}}{(2.185)^{17}}<1.68$. \\
So we have $B_{14}<\frac{B_{29}}{\varepsilon^3}<2.3985$.\\
Using Claims(iv) and (v) we get each of $B_2, B_3,\cdots,B_{17}$ is $<\frac{B_{29}}{\varepsilon^3}<2.3985$.\vspace{2mm}\\
\noindent{\bf Final Contradiction}

 Now
$2B_2+2B_4+\cdots+2B_{30}<
2
(8\times\frac{B_{29}}{\varepsilon^3})+2(\frac{3/2}{\varepsilon^2}+\frac{1}{\varepsilon^2}+\frac{3/2}{\varepsilon}+\frac{1}{\varepsilon}+\frac{3}{2}+1)B_{28}+2
B_{30}<49.86$ for $B_{30}<0.231$, $B_{29}<0.247$ and $B_{28}<0.285$, giving thereby a contradiction to the weak inequality $(2,2,\cdots,2,2)_w$.~$\Box$

\subsection{$n=31$}
Here $\omega_{31}=53.04$, $B_1\leq\gamma_{31}<4.8483483$. Using \eqref{2.3}, we have $l_{31}=0.0991<B_{31}<4.4974263=m_{31}$. \vspace{2mm}\\
\noindent{\bf Claim(i)} $B_{31}<0.173$

Suppose $B_{31}\geq0.173$. The inequality $(30^*,1)$ gives $49.86(B_{31})^{\frac{-1}{30}}+B_{31}>53.04$. But this is not true for $0.173\leq B_{31}\leq 4.4974263$. So we must have $B_{31}<0.173$.\\
\noindent{\bf Claim(ii)} $B_i<max\{B_{15},B_{16},2.5199\},~ {\mbox {for}}~ i=2,3,4,5,6,7,8,9,10,11,12,13,14$

 Using \eqref{2.2} and Lemmas 8, 9 we have\vspace{2mm}\\
 $~~~~~1=\lambda_1^{(31)}<B_1<\mu_1^{(31)}=4.8483483,$\\
 $~~~~~\frac{3}{4}=\lambda_2^{(31)}<B_1B_2<\mu_2^{(31)}=(1.1120388)^{29},$\\
 $~~~~~\frac{1}{2}=\lambda_3^{(31)}<B_1B_2B_3<\mu_3^{(31)}=(1.1744085)^{28},$\\
  $~~~~~\frac{1}{4}=\lambda_4^{(31)}<B_1B_2B_3B_4<\mu_4^{(31)}=(1.2415716)^{27},$\\
 $~~~~~\frac{\varepsilon}{4}=\lambda_5^{(31)}<B_1\cdots B_5<\mu_5^{(31)}=(1.3140139)^{26},$\\
 $~~~~~\frac{3\varepsilon^2}{16}=\lambda_6^{(31)}<B_1\cdots B_6<\mu_6^{(31)}=(1.3922837)^{25},$\\
 $~~~~~\frac{\varepsilon^3}{8}=\lambda_{7}^{(31)}<B_1\cdots B_7<\mu_7^{(31)}=(1.47700291)^{24},$\\
 $~~~~~\frac{\varepsilon^4}{16}=\lambda_{8}^{(31)}<B_1\cdots B_8<\mu_8^{(31)}=(1.5688842)^{23},$\\
 $~~~~~\frac{\varepsilon^6}{16}=\lambda_{9}^{(31)}<B_1\cdots B_9<\mu_9^{(31)}=(1.6687218)^{22},$\\
 $~~~~~\frac{\varepsilon^{45}}{4^{5}B_{31}^{21}}=\lambda_{10}^{(31)}<B_1\cdots B_{10}<\mu_{10}^{(31)}=(1.7774458)^{21},$\\
 $~~~~~\frac{\varepsilon^{40}}{4^5B_{31}^{20}}=\lambda_{11}^{(31)}<B_1\cdots B_{11}<\mu_{11}^{(31)}=(1.8961171)^{20},$\\
 $~~~~~\frac{\varepsilon^{36}}{2\times4^4B_{31}^{19}}=\lambda_{12}^{(31)}<B_1\cdots B_{12}<\mu_{12}^{(31)}=(2.0259616)^{19},$\\
 $~~~~~\frac{3\varepsilon^{32}}{4^5B_{31}^{18}}=\lambda_{13}^{(31)}<B_1\cdots B_{13}<\mu_{13}^{(31)}=(2.1684016)^{18},$\\
 $~~~~~\frac{\varepsilon^{28}}{4^4B_{31}^{17}}=\lambda_{14}^{(31)}<B_1\cdots B_{14}<\mu_{14}^{(31)}=(2.3250968)^{17},$\\
 $~~~~~\frac{\varepsilon^{24}}{4^4B_{31}^{16}}=\lambda_{15}^{(31)}<B_1\cdots B_{15}<\mu_{15}^{(31)}=(2.497994)^{16},$\\
  $~~~~~\frac{\varepsilon^{21}}{2\times4^3B_{31}^{15}}=\lambda_{16}^{(31)}<B_1\cdots B_{16}<\mu_{16}^{(31)}=(2.6893851)^{15}.$\vspace{2mm}\\
We find that $max\{\phi_{s,n-s}(\lambda_{s}^{(31)}),\phi_{s,n-s}(\mu_{s}^{(31)}), {\mbox {for}}~s=1,2,3,4,5,6,7,8,9,10,11,12,13\}=\phi_{13,18}(\mu_{1}^{(31)})$, which is $<\omega_{31}$. Therefore using Lemma 10 we have
\begin{equation}\label{4.6.13}B_i<max\{B_{i+1},B_{i+2},\cdots,B_{31}\},~ {\mbox {for}}~ i=2,3,4,5,6,7,8,9,10,11,12,13,14.\end{equation}
From \eqref{2.0a},\eqref{2.0b} and Claim (i), we find that $max\{B_{17},B_{18},\cdots,B_{31}\}<B_{17}\leq\frac{3}{2}\frac{B_{31}}{\varepsilon^3}<2.5199$. So using \eqref{4.6.13} we get that each of $B_2, \cdots,B_{14}$ is $<max\{B_{15},B_{16},\frac{3}{2}\frac{B_{31}}{\varepsilon^3}\}<max\{B_{15},B_{16},2.5199\}$. \\
\noindent{\bf Claim(iii)} $B_{15},B_{16}<2.5199$

We find that for $s=15,14$, $\phi_{s,n-s}(\lambda_{s}^{(31)})<\omega_{31}$, but $\phi_{s,n-s}(\mu_{s}^{(31)})>\omega_{31}$, so we apply Lemma 11 respectively with $\sigma_{15}^{(31)}=(2.42)^{16}$ and $\sigma_{14}^{(31)}=(2.299)^{17}$. Here $\phi_{15,16}(\sigma_{15}^{(31)})<\omega_{31}$ and $\phi_{14,17}(\sigma_{14}^{(31)})<\omega_{31}$.\\
First consider Lemma 11 for $s=15$ and with $\sigma_{15}^{(31)}=(2.42)^{16}$.\\
In Case(i), when $B_1B_2\cdots B_{15}<(2.42)^{16}$, then we have $B_{16}<max\{B_{17},B_{18},\cdots,B_{31}\}$, which is $<\frac{3}{2}\frac{B_{31}}{\varepsilon^3}<2.5199$.\\
In Case(ii), when $B_1B_2\cdots B_{15}\geq(2.42)^{16}$, then we have $B_{16}<\frac{\mu_{16}^{(31)}}{\sigma_{15}^{(31)}}<\frac{(2.6894)^{15}}{(2.42)^{16}}<2.0128$. \\
So we have $B_{16}\leq\frac{3}{2}\frac{B_{31}}{\varepsilon^3}<2.5199$.\\
Now consider Lemma 11 for $s=14$ and with $\sigma_{14}^{(31)}=(2.299)^{17}$.\\
In Case(i), when $B_1B_2\cdots B_{14}<(2.299)^{17}$, then we have $B_{15}<max\{B_{16},B_{17},\cdots,B_{31}\}$, which is $\leq\frac{3}{2}\frac{B_{31}}{\varepsilon^3}<2.5199$.\\
In Case(ii), when $B_1B_2\cdots B_{14}\geq(2.299)^{17}$, then we have $B_{15}<\frac{\mu_{15}^{(31)}}{\sigma_{14}^{(31)}}<\frac{(2.498)^{16}}{(2.299)^{17}}<1.642$. \\
So we have $B_{15}\leq\frac{3}{2}\frac{B_{31}}{\varepsilon^3}<2.5199$.\\
Using Claims(ii) and (iii) we get each of $B_2, B_3,\cdots,B_{17}$ is $\leq\frac{3}{2}\frac{B_{31}}{\varepsilon^3}<2.5199$.\vspace{2mm}\\
\noindent{\bf Final Contradiction }

Now
$2B_2+B_3+2B_5+\cdots+2B_{31}<
3(\frac{3}{2}\frac{B_{31}}{\varepsilon^3})+2
(7\times\frac{3/2}{\varepsilon^3}+\frac{1}{\varepsilon^3}+\frac{3/2}{\varepsilon^2}+\frac{1}{\varepsilon^2}+\frac{3/2}{\varepsilon}+\frac{1}{\varepsilon}+\frac{3}{2}+1)B_{31}<53.04$ for $B_{31}<0.173$, giving thereby a contradiction to the weak inequality $(2,1,2,\cdots,2,2)_w$.~$\Box$
\subsection{$n=32$}
Here $\omega_{32}=56.06$, $B_1\leq\gamma_{32}<4.9680344$. Using \eqref{2.3} we have $l_{32}=0.0942<B_{32}<4.6149714=m_{32}$. \vspace{2mm}\\
\noindent{\bf Claim(i)} $B_{32}<0.203$

Suppose $B_{32}\geq0.203$. The inequality $(31^*,1)$ gives $53.04(B_{32})^{\frac{-1}{31}}+B_{32}>56.06$. But this is not true for $0.203\leq B_{32}\leq 4.6149714$. So we must have $B_{32}<0.203$.\vspace{2mm}\\
\noindent{\bf Claim(ii)} $B_{30}<0.261$

Suppose $B_{30}\geq0.261$. The inequality $(29^*,3)$ gives $47.18({B_{30}B_{31}B_{32}})^{\frac{-1}{29}}+4B_{30}-\frac{B_{30}^3}{B_{31}B_{32}}>56.06$. But this is not true for $0.261\leq B_{30}\leq \frac{3}{2}B_{32}<0.305$ and $\frac{1}{2}B_{30}^2<B_{31}B_{32}<\frac{4}{3}(0.203)^2$. So we must have $B_{30}<0.261$.\vspace{2mm}\\
\noindent{\bf Claim(iii)} $B_i<max\{B_{15},B_{16},B_{17},2.5344\},~ {\mbox {for}}~ i=2,3,4,5,6,7,8,9,10,11,12,13,14$
 Using \eqref{2.2} and Lemmas 8, 9 we have\vspace{2mm}\\
 $~~~~~1=\lambda_1^{(32)}<B_1<\mu_1^{(32)}=4.9680344,$\\
 $~~~~~\frac{3}{4}=\lambda_2^{(32)}<B_1B_2<\mu_2^{(32)}=(1.1099695)^{30},$\\
 $~~~~~\frac{1}{2}=\lambda_3^{(32)}<B_1B_2B_3<\mu_3^{(32)}=(1.1710561)^{29},$\\
  $~~~~~\frac{1}{4}=\lambda_4^{(32)}<B_1B_2B_3B_4<\mu_4^{(32)}=(1.2367358)^{28},$\\
 $~~~~~\frac{\varepsilon}{4}=\lambda_5^{(32)}<B_1\cdots B_5<\mu_5^{(32)}=(1.3074634)^{27},$\\
 $~~~~~\frac{3\varepsilon^2}{16}=\lambda_6^{(32)}<B_1\cdots B_6<\mu_6^{(32)}=(1.3837502)^{26},$\\
 $~~~~~\frac{\varepsilon^3}{8}=\lambda_{7}^{(32)}<B_1\cdots B_7<\mu_7^{(32)}=(1.4661739)^{25},$\\
 $~~~~~\frac{\varepsilon^4}{16}=\lambda_{8}^{(32)}<B_1\cdots B_8<\mu_8^{(32)}=(1.5553893)^{24},$\\
 $~~~~~\frac{\varepsilon^6}{16}=\lambda_{9}^{(32)}<B_1\cdots B_9<\mu_9^{(24)}=(1.6521469)^{23},$\\
 $~~~~~\frac{3\varepsilon^{50}}{4^{6}B_{32}^{22}}=\lambda_{10}^{(32)}<B_1\cdots B_{10}<\mu_{10}^{(32)}=(1.7572829)^{22},$\\
 $~~~~~\frac{\varepsilon^{45}}{4^5B_{32}^{21}}=\lambda_{11}^{(32)}<B_1\cdots B_{11}<\mu_{11}^{(32)}=(1.8717771)^{21},$\\
 $~~~~~\frac{\varepsilon^{40}}{4^5B_{32}^{20}}=\lambda_{12}^{(32)}<B_1\cdots B_{12}<\mu_{12}^{(32)}=(1.9967464)^{20},$\\
 $~~~~~\frac{\varepsilon^{36}}{2\times4^4B_{32}^{19}}=\lambda_{13}^{(32)}<B_1\cdots B_{13}<\mu_{13}^{(32)}=(2.1334819)^{19},$\\
 $~~~~~\frac{3\varepsilon^{32}}{4^5B_{32}^{18}}=\lambda_{14}^{(32)}<B_1\cdots B_{14}<\mu_{14}^{(32)}=(2.2834814)^{18},$\\
 $~~~~~\frac{\varepsilon^{28}}{4^4B_{32}^{17}}=\lambda_{15}^{(32)}<B_1\cdots B_{15}<\mu_{15}^{(32)}=(2.4484926)^{17},$\\
  $~~~~~\frac{\varepsilon^{24}}{4^4B_{32}^{16}}=\lambda_{16}^{(32)}<B_1\cdots B_{16}<\mu_{16}^{(32)}=(2.6305651)^{16},$\\
   $~~~~~\frac{\varepsilon^{21}}{2\times4^3B_{32}^{15}}=\lambda_{17}^{(32)}<B_1\cdots B_{17}<\mu_{17}^{(32)}=(2.8321145)^{15}.$\vspace{2mm}\\
We find that $max\{\phi_{s,n-s}(\lambda_{s}^{(32)}),\phi_{s,n-s}(\mu_{s}^{(32)}), {\mbox {for}}~s=1,2,3,4,5,6,7,8,9,10,11,12,13\}=\phi_{1,31}(\mu_{1}^{(32)})$, which is $<\omega_{32}$. Therefore using Lemma 10 we have
\begin{equation}\label{4.6.14}B_i<max\{B_{i+1},B_{i+2},\cdots,B_{32}\},~ {\mbox {for}}~ i=2,3,4,5,6,7,8,9,10,11,12,13,14.\end{equation}
From \eqref{2.0a},\eqref{2.0b} and Claims (i), (ii), we find that $max\{B_{18},B_{19},\cdots,B_{32}\}<B_{18}\leq\frac{B_{30}}{\varepsilon^3}<2.5344$. So using \eqref{4.6.14} we get that each of $B_2, \cdots,B_{14}$ is $<max\{B_{15},B_{16},B_{17},\frac{B_{30}}{\varepsilon^3}\}<max\{B_{15},B_{16},B_{17},2.5344\}$.\vspace{2mm} \\
\noindent{\bf Claim(iv)} $B_{15},B_{16},B_{17}<2.5344$

We find that for $s=16,15,14$, $\phi_{s,n-s}(\lambda_{s}^{(32)})<\omega_{32}$, but $\phi_{s,n-s}(\mu_{s}^{(32)})>\omega_{32}$, so we apply Lemma 11 respectively with $\sigma_{16}^{(32)}=(2.506)^{16}$, $\sigma_{15}^{(32)}=(2.4)^{17}$ and $\sigma_{14}^{(32)}=(2.27)^{18}$. Here
$\phi_{16,16}(\sigma_{16}^{(32)})<\omega_{32}$, $\phi_{15,17}(\sigma_{15}^{(32)})<\omega_{32}$ and $\phi_{14,18}(\sigma_{14}^{(32)})<\omega_{32}$.\\
First consider Lemma 11 for $s=16$ and with $\sigma_{16}^{(32)}=(2.506)^{16}$.\\
In Case(i), when $B_1B_2\cdots B_{16}<(2.506)^{16}$, then we have $B_{17}<max\{B_{18},B_{19},\cdots,B_{32}\}$, which is $\leq\frac{B_{30}}{\varepsilon^3}<2.5344$.\\
In Case(ii), when $B_1B_2\cdots B_{16}\geq(2.506)^{16}$, then we have $B_{17}<\frac{\mu_{17}^{(32)}}{\sigma_{16}^{(32)}}<\frac{(2.833)^{15}}{(2.506)^{16}}<2.512$. \\
So we have $B_{17}<\frac{B_{30}}{\varepsilon^3}<2.5344$.\\
Now consider Lemma 11 for $s=15$ and with $\sigma_{15}^{(32)}=(2.4)^{17}$.\\
In Case(i), when $B_1B_2\cdots B_{15}<(2.4)^{17}$, then we have $B_{16}<max\{B_{17},B_{18},\cdots,B_{32}\}$, which is $<\frac{B_{30}}{\varepsilon^3}<2.5344$.\\
In Case(ii), when $B_1B_2\cdots B_{15}\geq(2.4)^{17}$, then we have $B_{16}<\frac{\mu_{16}^{(32)}}{\sigma_{15}^{(32)}}<\frac{(2.6306)^{16}}{(2.4)^{17}}<1.809$. \\
So we have $B_{16}<\frac{B_{30}}{\varepsilon^3}<2.5344$.\\
Next consider Lemma 11 for $s=14$ and with $\sigma_{14}^{(32)}=(2.27)^{18}$.\\
In Case(i), when $B_1B_2\cdots B_{14}<(2.27)^{18}$, then we have $B_{15}<max\{B_{16},B_{17},\cdots,B_{32}\}$, which is $<\frac{B_{30}}{\varepsilon^3}<2.5344$.\\
In Case(ii), when $B_1B_2\cdots B_{14}\geq(2.27)^{18}$, then we have $B_{15}<\frac{\mu_{15}^{(32)}}{\sigma_{14}^{(32)}}<\frac{(2.4485)^{17}}{(2.27)^{18}}<1.596$. \\
So we have $B_{15}<\frac{B_{30}}{\varepsilon^3}<2.5344$.\\
Using Claims(iii) and (iv) we get each of $B_2, B_3,\cdots,B_{18}$ is $<\frac{B_{30}}{\varepsilon^3}<2.5344$.\vspace{2mm}\\
\noindent{\bf Final Contradiction}

 Now
$2B_2+2B_4+2B_6+\cdots+2B_{32}<
2(9\times\frac{1}{\varepsilon^3}+\frac{3/2}{\varepsilon^2}+\frac{1}{\varepsilon^2}+
\frac{3/2}{\varepsilon}+\frac{1}{\varepsilon}+\frac{3}{2}+1)B_{30}+2
B_{32}<56.06$ for $B_{32}<0.203$ and $B_{30}<0.261$, giving thereby a contradiction to the weak inequality $(2,2,2,\cdots,2,2)_w$.~$\Box$
\subsection{$n=33$}
Here $\omega_{33}=59.58$, $B_1\leq\gamma_{33}<5.0876409$. Using \eqref{2.3}, we have $l_{33}=0.0896<B_{33}<4.7324725=m_{33}$. \vspace{2mm}\\
\noindent{\bf Claim(i)} $B_{33}<0.155$

Suppose $B_{33}\geq 0.155$.
The inequality $(32^*,1)$ gives $56.06(B_{33})^{\frac{-1}{32}}+B_{33}>59.58$. But this is not true for $0.155
\leq B_{33}< 4.7324725$. \\
 \noindent{\bf Claim(ii)} $B_i<max\{B_{16},B_{17},B_{18},2.2577\},~ {\mbox {for}}~ i=2,3,4,5,6,7,8,9,10,11,12,13,14,15$

 Using \eqref{2.2} and Lemmas 8, 9 we have\vspace{2mm}\\
 $~~~~~1=\lambda_1^{(33)}<B_1<\mu_1^{(33)}=5.0876409,$\\
 $~~~~~\frac{3}{4}=\lambda_2^{(33)}<B_1B_2<\mu_2^{(33)}=(1.1079939)^{31},$\\
 $~~~~~\frac{1}{2}=\lambda_3^{(33)}<B_1B_2B_3<\mu_3^{(33)}=(1.1678596)^{30},$\\
  $~~~~~\frac{1}{4}=\lambda_4^{(33)}<B_1B_2B_3B_4<\mu_4^{(33)}=(1.2321321)^{29},$\\
 $~~~~~\frac{\varepsilon}{4}=\lambda_5^{(33)}<B_1B_2B_3B_4B_5<\mu_5^{(33)}=(1.3012373)^{28},$\\
 $~~~~~\frac{3\varepsilon^2}{16}=\lambda_6^{(33)}<B_1B_2B_3B_4B_5B_6<\mu_6^{(33)}=(1.3756537)^{27},$\\
 $~~~~~\frac{\varepsilon^3}{8}=\lambda_{7}^{(33)}<B_1\cdots B_7<\mu_7^{(33)}=(1.4559193)^{26},$\\
 $~~~~~\frac{\varepsilon^4}{16}=\lambda_{8}^{(33)}<B_1\cdots B_8<\mu_8^{(33)}=(1.5426418)^{25},$\\
 $~~~~~\frac{\varepsilon^6}{16}=\lambda_{9}^{(33)}<B_1\cdots B_9<\mu_9^{(24)}=(1.6365102)^{24},$\\
 $~~~~~\frac{\varepsilon^{55}}{2\times4^{5}B_{33}^{23}}=\lambda_{10}^{(33)}<B_1\cdots B_{10}<\mu_{10}^{(33)}=(1.7383141)^{23},$\\
 $~~~~~\frac{3\varepsilon^{50}}{4^6B_{33}^{22}}=\lambda_{11}^{(33)}<B_1\cdots B_{11}<\mu_{11}^{(33)}=(1.8489335)^{22},$\\
 $~~~~~\frac{\varepsilon^{45}}{4^5B_{33}^{21}}=\lambda_{12}^{(33)}<B_1\cdots B_{12}<\mu_{12}^{(33)}=(1.9693991)^{21},$\\
 $~~~~~\frac{\varepsilon^{40}}{4^5B_{33}^{20}}=\lambda_{13}^{(33)}<B_1\cdots B_{13}<\mu_{13}^{(33)}=(2.1008861)^{20},$\\
 $~~~~~\frac{\varepsilon^{36}}{2\times4^4B_{33}^{19}}=\lambda_{14}^{(33)}<B_1\cdots B_{14}<\mu_{14}^{(33)}=(2.24475296)^{19},$\\
 $~~~~~\frac{3\varepsilon^{32}}{4^5B_{33}^{18}}=\lambda_{15}^{(33)}<B_1\cdots B_{15}<\mu_{15}^{(33)}=(2.4025756)^{18},$\\
  $~~~~~\frac{\varepsilon^{28}}{4^4B_{33}^{17}}=\lambda_{16}^{(33)}<B_1\cdots B_{16}<\mu_{16}^{(33)}=(2.5761882)^{17},$\\
   $~~~~~\frac{\varepsilon^{24}}{4^4B_{33}^{16}}=\lambda_{17}^{(33)}<B_1\cdots B_{17}<\mu_{17}^{(33)}=(2.7677614)^{16},$\\
      $~~~~~\frac{\varepsilon^{21}}{2\times4^3B_{33}^{15}}=\lambda_{18}^{(33)}<B_1\cdots B_{18}<\mu_{18}^{(33)}=(2.9798226)^{15}.$\vspace{2mm}\\
We find that $max\{\phi_{s,n-s}(\lambda_{s}^{(33)}),\phi_{s,n-s}(\mu_{s}^{(33)}), {\mbox {for}}~s=1,2,3,4,5,6,7,8,9,10,11,12,13,14\}=\phi_{14,19}(\mu_{14}^{(33)})$, which is $<\omega_{33}$. Therefore using Lemma 10 we have
\begin{equation}\label{4.6.15}B_i<max\{B_{i+1},B_{i+2},\cdots,B_{33}\},~ {\mbox {for}}~ i=2,3,4,5,6,7,8,9,10,11,12,13,14,15.\end{equation}
From \eqref{2.0a},\eqref{2.0b} and Claim(i) we find that, $max\{B_{19},B_{20},\cdots,B_{33}\}<B_{19}\leq\frac{3}{2}\frac{B_{33}}{\varepsilon^3}<2.2577$. So using \eqref{4.6.15} we get that each of $B_2, \cdots,B_{15}$ is $<max\{B_{16},B_{17},B_{18},\frac{3}{2}\frac{B_{33}}{\varepsilon^3}\}<max\{B_{16},B_{17},B_{18},2.2577\}$.\vspace{2mm} \\
\noindent{\bf Claim(iii)} $B_{16},B_{17},B_{18}<2.7$

We find that for $s=17,16,15$, $\phi_{s,n-s}(\lambda_{s}^{(33)})<\omega_{33}$, but $\phi_{s,n-s}(\mu_{s}^{(33)})>\omega_{33}$, so we apply Lemma 11 respectively with $\sigma_{17}^{(33)}=(2.616)^{16}$, $\sigma_{16}^{(33)}=(2.5)^{17}$ and $\sigma_{15}^{(33)}=(2.39)^{18}$. Here
$\phi_{17,16}(\sigma_{17}^{(33)})<\omega_{33}$, $\phi_{16,17}(\sigma_{16}^{(33)})<\omega_{33}$ and $\phi_{15,18}(\sigma_{15}^{(33)})<\omega_{33}$.\\
First consider Lemma 11 for $s=17$ and with $\sigma_{17}^{(33)}=(2.616)^{16}$.\\
In Case(i), when $B_1B_2\cdots B_{17}<(2.616)^{16}$, then we have $B_{18}<max\{B_{19},B_{20},\cdots,B_{33}\}$, which is $<\frac{3}{2}\frac{B_{33}}{\varepsilon^3}<2.2577$.\\
In Case(ii), when $B_1B_2\cdots B_{17}\geq(2.616)^{16}$, then we have $B_{18}<\frac{\mu_{18}^{(33)}}{\sigma_{17}^{(33)}}<\frac{(2.9797)^{15}}{(2.616)^{16}}<2.7$. \\
So we have $B_{18}<2.7$.\\
Now consider Lemma 11 for $s=16$ and with $\sigma_{16}^{(33)}=(2.5)^{17}$.\\
In Case(i), when $B_1B_2\cdots B_{16}<(2.5)^{17}$, then we have $B_{17}<max\{B_{18},B_{19},\cdots,B_{33}\}$, which is $<2.7$.\\
In Case(ii), when $B_1B_2\cdots B_{16}\geq(2.5)^{17}$, then we have $B_{17}<\frac{\mu_{17}^{(33)}}{\sigma_{16}^{(33)}}<\frac{(2.7678)^{16}}{(2.5)^{17}}<2.038$. \\
So we have $B_{17}<2.7$.\\
Next consider Lemma 11 for $s=15$ and with $\sigma_{15}^{(33)}=(2.39)^{18}$.\\
In Case(i), when $B_1B_2\cdots B_{15}<(2.39)^{18}$, then we have $B_{16}<max\{B_{17},B_{18},\cdots,B_{33}\}$, which is $<2.7$.\\
In Case(ii), when $B_1B_2\cdots B_{15}\geq(2.39)^{18}$, then we have $B_{16}<\frac{\mu_{16}^{(33)}}{\sigma_{15}^{(33)}}<\frac{(2.5762)^{17}}{(2.39)^{18}}<1.498$. \\
So we have $B_{16}<2.7$.\\
Using Claims(ii) and (iii) we get each of $B_2, B_3,\cdots,B_{18}$ is $<2.7$. \vspace{2mm}\\
\noindent{\bf Final Contradiction}

Now
$2B_2+B_3+2B_5+2B_7+\cdots+2B_{33}<
3\times2.7+2(7\times2.7)+2(\frac{3/2}{\varepsilon^3}+\frac{1}{\varepsilon^3}+\frac{3/2}{\varepsilon^2}+\frac{1}{\varepsilon^2}+
\frac{3/2}{\varepsilon}+\frac{1}{\varepsilon}+\frac{3}{2}+1)B_{33}<59.58$ for $B_{33}<0.155$, giving thereby a contradiction to the weak inequality $(2,1,2,\cdots,2,2)_w$.~$\Box$

\end{document}